Oscar Sheynin

**Theory of Probability and Statistics
As Exemplified in Short Dictums**

Second revised and enlarged edition

Berlin, 2009

**Contents**

> Pour prévoir des mathématiques
> la vraie méthode est d'étudier
> leur histoire et leur état present.
> Poincaré (1909, p. 167)

> I do feel how wrongful it was to
> work for so many years at
> statistics and neglect its history.
> K. Pearson (1978, p. 1)

## 1. Introduction

I am presenting a first-ever scientific collection of short sayings on probability and statistics expressed by most various men of science, many classics included, from antiquity to Kepler to our time. Quite understandably, the reader will find here no mathematical formulas and in some instances he will miss a worthy subject. Markov chains provide a good example: their inventor had not said anything about them suitable for my goal. Nevertheless, the scope of the collected dictums is amazingly broad which reflects both the discussions concerning the lack of a solid foundation of probability theory until the 1930s and the great extent of applications of probability and statistics. And I have also included two related and most important topics, the treatment of observations and randomness.

The present arrangement of the collected material is the result of several attempts, but some further improvement might have perhaps been possible here. In any case, I have discovered that a chronological arrangement, even of each separate subject, is not good enough. Having been unable to check several sayings included in a mostly popular collection Gaither et al (1996), I copied them with proper references.

I am reasonably sure that this work will prove useful both for scientific and pedagogic purposes, – for providing a background even in unexpected cases and for greatly stimulating the study of probability and statistics. As to a general description of the history of these disciplines, I refer readers to Hald (1990; 1998) and Sheynin (2005/2009), also available at www. sheynin.de

This second edition is much extended, several mistakes have been corrected and the material is more properly arranged.

## 2. Theory of Probability

### 2.1. Origin and Aims

**1.** It is in basic attitudes towards the phenomenal world, in religious and moral teachings and barriers, that I incline to seek for an explanation of the delay [of the origin of probability theory]. Mathematics never leads thought, but only expresses it. Kendall (1956, § 31).

**Comment.** In particular, Kendall mentioned the gamblers' psychology.

**2.** Two persons claim a cloak [or aught else. The first says] The whole of it is mine. [The second says] Half of it is mine. [They share it in the ratio 3:1.] *Mishnah*, Baba Metria 1[1].

**Comment.** A thousand years later Pascal would have given the same answer and stochastically justified it. *Mishnah* (Blackman 1951 – 1955) was compiled in the 4th century of our era. It consists of more than 60 treatises and is the first part of the Talmud.

**3.** Neither joint-stock societies, nor banks, nor stock exchanges stood in need of probability; their demands on probability only appeared in the 19th century, when methods of scientific gain superseded downright robbery. Mrochek (1934, p. 50).

**Comment.** Mrochek could have also mentioned marine and life insurance, cf. Ostrogradsky (1847): "Do not pay more; … try to pay less so as to gain something. Do not worry about the insurance society: it will not incur losses".

**4.** Joignant la rigueur des demonstrations de la science [matheseos] à l'incertitude du hasard, et coinciliant ces choses en apparence contraire, elle peut, tirant son nom des deux, s'arroger à bon droit ce titre stupéfiant: La Géomètrie du hasard. Pascal (1654/1998, p. 172).

**5.** Toutefois je veux croire qu'en considérant ces choses [games of chance] plus attentivement, le lecteur apercevra bientôt qu'il ne s'agit pas ici d'un simple jeu d'esprit, mais qu'on y jette les fondements d'une spéculation fort intéressante et profonde. Huygens (1657/1920, p. 58).

**6.** Ich wünschte, dass jemand verschiedenen Arten von Glücksspielen (bei denen es schöne Beispiele gibt für diese Theorie) mathematisch bearbeitete. Das wäre zugleich anziehend und nützlich und nicht unwürdig Deiner oder eines anderen grossen Mathematikers. Leibniz, letter of 1703 to Jakob Bernoulli, in Latin; Gini (1946, p. 404).

**7.** There is a given Number of each of several Sorts of Things … put promiscuously together; out of which a given Number … is to be taken as it happens; To find the probability that there shall come out precisely a given Number of each Sort … T. Simpson (1740, Problem 6).

**Comment.** Ostrogradsky (1848) considered an equivalent problem which he reasonably thought to be useful for, as we say now, statistical control of supplies, but his main formula (p. 342) was extremely involved and hardly ever checked by anyone else. My main point here is that Huygens (1657, Additional problem No. 4) solved the same problem as Simpson did which goes to show that the emerging theory of probability was indeed not a "simple jeu d'esprit" (see No. 5). See also No. 8. And here is Poisson (1837, p. 1):

*Les géomètres du XVII$^e$ siècle qui se sont occupés du calcul des probabilités, ne l'ont employé qu'à déterminer les chances de différents jeux de cette époque; et ce n'est que dans le siècle suivant qu'il a pris toute son extension, et qu'il est devenu une des principales branches des mathématiques, soit par le nombre et l'utilité de ses applications, soit par le genre d'analyse auquel il a donné naissance.*

**8.** An insect is placed on a table; find the chance that it will be found at time *t* at distance *r* from the point where it has started; if the insect be supposed to hop, we have the simpler case of which that of the crawling insect is the limit. Glaisher (1873).

**Comment.** Random walks with applications to natural science had begun to be studied in earnest from the mid-19$^{th}$ century (Dutka 1985), but games of chance can be interpreted as such.

**9.** The Reader may here observe the Force of Numbers, which can be successfully applied, even to those things, which one would imagine are subject to no Rules. There are very few things which we know, which are not capable of being reduc'd to a Mathematical Reasoning, and when they cannot, it's a Sign our Knowledge of them is very small and confus'd. And where mathematical reasoning can be had, it's a great folly to make use of any other, as to grope for a thing in the dark, when you have a Candle standing by you. I believe the Calculation of the Quantity of probability might be improved to a very useful and pleasant Speculation, and applied to a great many Events which are accidental, besides those of Games … Arbuthnot (1692, Introduction) as quoted by Gaither et al (1996).

**Comment.** Todhunter (1865, p. 49) believed that that work was a translation of the Huygens tract (1657).

**10.** We should not take into account the number of the doubts but rather consider how great is their incongruity and what is their disagreement with what exists. Sometimes a single doubt is more powerful than a thousand other doubts. Maimonides (1963, II-23).

**Comment.** Cf. J. Bernoulli (1713, pt. 4, chapter 3 and beginning of chapter 2). In chapter 2 he stated: "Probabilities are estimated both by the *number* and *the weight of the arguments* …"

**11.** To make *conjectures* about something is the same as to measure its probability. Therefore, the art of *conjecturing* or *stochastics* (*ars conjectandi* sive *stochastice*) is defined as the art of measuring the probability of things as exactly as possible, to be able always to choose what will be found the best, the more satisfactory, serene and reasonable for our judgements and actions. This alone supports all the wisdom of the philosopher and the prudence of the politician. Jakob Bernoulli (1713/2005, p. 17).

**12.** Under uncertain and dubious circumstances we ought to suspend our actions until more light is thrown. If, however, the necessity of action brooks no delay, we must always choose from among two possibilities that one which seems more suitable, safe, reasonable, or at least more probable, even if none of them is actually such. Ibidem, p. 19.

**13.** Since complete certitude can only seldom be attained, necessity and custom desire that that, which is only morally certain, be considered as absolutely certain. Ibidem, p. 20.

**14.** We must hold fast to the most probable account. Plato (1929, § 44d, p. 99).

**Comment.** Descartes (No. 19) and later scholars adhered to the same opinion.

**15.** Es gibt also drei Grundhaltungen: zwei fehlerhafte, durch Übermaß und Unzulänglichkeit gekennzeichnet, und eine richtige: die Mitte. Aristotle 1979, 1108b, p. 40).

**Comment.** Aristotle (1107b – 1109b) described this subject in detail.

**16.** Many things are probable and … though these are not demonstrably true, they guide the life of the wise man because they are so significant and clear-cut. Cicero (1997, Book 1, § 12, p. 7).

**17.** Le mieux que nous puissions faire quand nous sommes engagez à prendre parti, est d'embrasser le plus probable, puisque ce seroit vu renversement de la raison d'embrasser le moins probable. Arnauld & Nicole (1662/1992, p. 327).

**18.** [Certitude] morale, c'est-à-dire suffisante pour régler nos mœurs, ou aussi grand que celle des choses dont nous n'avons point coustume de douter touchant la conduite de la vie, bien que nous sçachions qu'il se peut faire, absolument parlant, qu'elles soient fausses. Descartes (1644/1978, pt. 4, No. 205, 483°, p. 323).

**19.** Ainsi, les actions de la vie ne souffrant souvent aucun delay, c'est une vérité très certaine que, lorsqu'il n'est pas en notre pouvoir de discerner les plus vrayes opinions, nous devons suivre les plus probables. Descartes (1637/1982, pt. 3, p. 25).

**20.** Was aber den Ausdruck *Stochastik* anlangt, so bedarf er keiner Rechtfertigung. Denn er findet sich – und zwar in dem ihm von mir beigelegten Sinne – schon in Jakob Bernoullis *Ars Conjectandi*. Bortkiewicz (1917, p. x).

**Comment.** For more detail about the history of that term see Sheynin (2009a, Note 1 to Chapter 3). Plato and Socrates applied it and it had been in usage in England at least from 1662 onward. For Wallis, in 1685, *stochastic* meant *iterative* (process). Lastly, Prevost & Lhuilier (1799, p. 3) forestalled Bortkiewicz.

**21.** Einer … mir bete, ich sollte jhm sagen, ob sein Freundt in fernen Landen lebendt oder todt were … Und ich … sagte jhm, ja oder nein, so were ich ein Ariolus und er ein Verbrecher an Gottes Gebott … Kepler (1610/1941, p. 238).

**Comment.** J. Bernoulli (1713, pt. 4, chapter 4) was prepared to declare a missing person alive or dead by weighing the possible arguments against each other.

**22.** There is no probability so great as not to allow of a contrary possibility. Hume (1739/1969, Book 1, pt. 3, § 12, p. 135).

**23.** La théorie des hasards a pour objet de déterminer ces fractions [the probabilities], et l'on voit par là que c'est le supplément le plus heureux que l'on puisse imaginer à l'incertitude de nos connaissances. Laplace (1776/1891, p. 146).

**24.** La théorie des hasards consiste … à réduire tous les événements qui peuvent avoir lieu relativement à un objet, dans un certain nombre de cas également possibles, c'est-à-dire, tels que nous soyons également indécis sur leur existence, & à déterminer le nombre de cas favorables à l'événement dont on cherche la probabilité. Le rapport de ce nombre à celui de tous les cas possibles, est la mesure de cette probabilité. Laplace (1786/1893, p. 296).

**25.** La théorie des probabilités n'est, au fond, que le bon sens réduit au calcul. Laplace (1814/1886, p. CLIII).

**Comment.** At the time, the same definition could have described mathematics as a whole.

**26.** There is no more remarkable feature in the mathematical theory of probability than the manner in which it has been found to harmonize with, and justify, the conclusions to which mankind have been led, not by reasoning, but by instinct and experience, both of the individual and of the race. At the same time it has corrected, extended, and invested them with a definiteness and precision of which these crude, though sound, appreciations of common sense were till then devoid. Crofton (1885), as quoted by Gaither et al (1996).

**27.** The science of probabilities that goes under the name of the theory of probability has as its subject the determination of the probability of an event given its connection with events whose probabilities are known. Chebyshev (1845/1951, p. 29).

**28.** The object of the theory of probabilities may be thus stated: Given the separate probabilities of any propositions, to find the probability of another proposition. Boole (1851/1952, p. 251). Alternatively (Boole 1854b/2003, p. 246): … might be thus defined. Given the probabilities of any events, of whatever kind, to find the probability of some other event connected with them.

**Comment.** De Moivre described four aims of his doctrine of chances. The main goal was to separate chance from Divine Design, see No. 29. For Laplace, the theory of probability was an applied mathematical discipline.

**29.** Sir, The greatest help I have received in writing upon this subject [the doctrine of chances] having been from your incomparable Works, especially your Method of Series; I think it my duty publickly to acknowledge, that the Improvements I have made, in the matter here treated of, are principally derived from yourself. The great benefit which has accrued to me in this respect, requires my share in the general Tribute of Thanks due to you from the learned World: But one Advantage which is more particularly my own, is the Honour I have frequently had of being admitted to your private Conversation; wherein the Doubts I have had upon any Subject relating to *Mathematics*, have been resolved by you with the greatest Humanity and Condescension. Those marks of your Favour are the more valuable to me, because I had no other pretence to them but the earnest desire of understanding your sublime and universally useful Speculations. I should think my self very happy, if having given my readers a Method of calculating the Effects of Chance, as they are the result of Play, and thereby fixing certain Rules, for estimating how far some sort of Events may rather be owing to Design than Chance, I could by this small Essay excite in others a desire of prosecuting these Studies, and of learning from your Philosophy how to collect, by a just Calculation, the Evidences of exquisite Wisdom and Design, which appear in the *Phenomena* of Nature throughout the Universe. De Moivre (1718).

**Comment.** This is the Dedication of the first edition of the *Doctrine of Chances* to Newton reprinted in 1756 on p. 329.

**30.** Altho' Chance produces Irregularities, still the Odds will be infinitely great, that in the process of Time, those Irregularities will bear no proportion to the recurrency of that Order which naturally results from ORIGINAL DESIGN. De Moivre (1733/1756, p. 251).

**31.** Newton's idea of an omnipresent activating deity, who maintains mean statistical values, formed the foundation of statistical development through Derham [a religious philosopher], Süssmilch, Niewentyt [a statistician], Price to Quetelet and Florence Nightingale. K. Pearson (1926).

**Comment.** In a private communication, E. S. Pearson explained his understanding of "mean values": K. P. actually thought about maintaining stability.

**32.** The theory of probabilities is simply the science of logic quantitatively treated. Peirce (1878/1958, p. 278).

**33.** The subject-matter of calculation in the mathematical theory of probabilities is *quantity of belief*. Donkin (1851, p. 353).

**34.** [The theory of probability is] The art of judging in cases where only probable evidence can be obtained. Newcomb (1884).

**35.** Ich habe mehr als einmal gesagt, daß man eine neue Art Logik braucht, die die Grade der Wahrscheinlichkeit behandelt. Leibniz (1765/1961, p. 515).

**36.** [The theory of probability] may be considered as a branch of logic studying all the methods which the human mind applies for acquiring new truths. Vasiliev (1892, p. 644).

**37.** We may say quite definitely that without applying probability theory the further systematization of human knowledge or the development of science are impossible (p. 218/7). *The new contemporary stage in the development of scientific thought is characterized by the need to introduce the notion of probability into the statements of the elementary laws of nature* (p. 222/p. 11). Bernstein (1928/1964, pp. 218, 222; translation, 2005, pp. 7, 11).

**38.** The mathematical theory of probability is a science which aims at reducing to calculation, where possible, the amount of credence due to propositions or statements, or to the occurrence of events, future or past, more especially as contingent or dependent upon other propositions or events the probability of which is known. Crofton (1885) as quoted by Gaither et al (1996).

**39.** The study of inductive inference belongs to the theory of probability, since observational facts can make a theory only probable but will never make it absolutely certain. Reichenbach (1951, p. 231).

**40.** The theory of probability aims at determining the chances for the occurrence of some event. The word *event* means, in general, everything whose probability is being determined. In mathematics, the word *probability* thus serves to denote some magnitude subject to measurement. Chebyshev (1879 – 1880/1936, p. 148; translation 2004, p. 141).

**41.** The theory of probability studies mathematical models of random events, and, given the probabilities of some random events, makes it possible to determine the probabilities of other random events somehow connected with the first ones. Prokhorov & Sevastianov (1999, p. 77).

**42.** The actual science of logic is conversant at present only with things either certain, impossible, or entirely doubtful, none of which (fortunately) we have to reason on. Therefore the true Logic for this world is the Calculus of Probabilities (which is, or ought to be in a reasonable man's mind). This branch of Math. … is the only "Mathematics for Practical Men", as we ought to be. Maxwell, letter of 1850; Campbell & Garnett (1882/1884, p. 97).

**43.** The smallness of probability is compensated by the greatness of the evil; and the sensation is equally lively, as if the evil were more probable. Hume (1739/1969, Book 2, pt 3, § 9, p. 444).

**44.** We should, perhaps, spare a few moments to consider this stiff drab Victorian figure [Todhunter], so unlike the colourful authors of whom he wrote, so meticulous in his attention to detail and so blind to the broad currents of his subject; for his *History* [1865] has stood for nearly a hundred years, without an imitator or a rival, and we are all indebted to it. Kendall (1963, p. 205).

### 2.2. Main Notions

**45.** Probabilitas est gradus possibilitatis. Leibniz, manuscript of 1678, published 1901, as quoted by Biermann & Faak (1957).

**46.** As to *probability*, this is the degree of certainty, and it differs from the latter as a part from the whole. Namely, if the integral and absolute certainty, which we designate by letter α or by unity 1, will be thought to consist, for example, of five probabilities, as though of five parts, three of which favour the existence or realization of some event, with the other ones, however, being against it, we will say that this event has 3/5α, or 3/5 of certainty. Jakob Bernoulli (1713/2005, pp. 14 – 15).

**Comment.** This definition is not quite formal, and the equality of the five probabilities is not stipulated. In his further exposition of the unfinished *Ars Conjectandi* Bernoulli had not applied it.

**47.** If $p$ is the number of chances by which a certain event may happen, and $q$ is the number of chances by which it may fail, the happenings as much as the failings have their degree of probability; but if all the chances by which the event may happen or fail were equally easy, the probability of happening will be to the probability of failing as $p$ to $q$. De Moivre (1712/1984, p. 237).

**48.** The Probability of an Event is greater, or less, according to the number of chances by which it may Happen, compar'd with the number of all chances, by which it may either Happen or Fail. De Moivre (1718, beginning of Intro.) as quoted by Sylla (2006, p. 112).

It is the comparative magnitude of the number of Chances to happen, in respect to the whole number of Chances either to happen or to fail, which is the true measure of Probability. De Moivre (1738) as quoted by Schneider (1968, p. 279).

**49.** If we constitute a Fraction whereof the Numerator be the number of Chances whereby an Event may happen, and the Denominator the number of all the Chances whereby it may either happen or fail, that Fraction will be a proper designation of the Probability of happening. De Moivre (1756, pp. 1 – 2).

**50.** In all cases, the Expectation of obtaining any Sum is estimated by multiplying the value of the Sum expected by the Fraction which represents the Probability of obtaining it. Ibidem, p. 3.

**Comment.** This definition is possibly contained in the previous editions of the *Doctrine* as well. The concepts of probability and expectation had likely been made use of from the second half of the 16[th] century onward: the gains in the Genoise lottery had always been much lower than their expectations, see for example Biermann (1957). The Huygens treatise (1657), or, more precisely, its original Dutch text first published in 1660 was based on the concept of expectation (called *value of chance*). The Latin text of 1657 (translated by van Schooten from the Dutch) mentions lots or expectations (*sors sive expectatio*). De Moivre (1712/1984, p. 237) applied the term *sors* (fortune, lot) and Arbuthnot (1712/1977, p. 32) wrote *Lot*, or *Value of Expectation*.

**51.** *The probability of any event* is the ratio between the value at which an expectation depending on the happening of the event ought to be computed, and the value of the thing expected upon its happening. Bayes (1764/1970, p. 136).

**52.** La probabilité de l'existence d'un événement n'est ainsi que le rapport du nombre des cas favorables à celui de tous les cas possibles, lorsque nous ne voyons d'ailleurs aucune raison pour laquelle l'un de ces cas arriverait plutôt que l'autre. Elle peut être conséquemment représentée par une fraction dont le numérateur est le nombre des cas favorables, et le dénominateur celui de tous les cas possibles. Laplace (1776/1891, p. 146).

**Comment.** Laplace (1814/1886, p. XI) largely repeated this definition much later, see also NNo. 23 and 24. And Laplace (1786/1893, p. 296; 1814/1886, p. VIII) had also stated that "La probabilité est relative en partie à cette ignorance, en partie à nos connaissances".

In 1777, the first edition of the *Encyclopaedia Britannica* (vol. 3, p. 513) carried a short anonymous note, *Probability*, defining it as a logical concept.

**53.** Two Events are independent, when they have no connection one with the other, and that the happening of one neither forwards nor obstructs the happening of the other.

Two Events are dependent, when they are so connected together as that the Probability of either's happening is altered by the happening of the other. De Moivre (1718/1738, p. 6).
**Comment.** I did not see the first edition of 1718. Two great French mathematicians had thought of translating De Moivre's *Doctrine*, but abandoned that idea:

*Il est vrai que j'ai eu autrefois l'idée de donner une traduction de l'Ouvrage de Moivre, accompagnée de notes et d'additions de ma façon, et j'avais même déjà traduit une partie de cet Ouvrage; mais j'ai depuis longtemps renoncé à ce projet, et je suis enchanté d'apprendre que vous en avez entrepris l'exécution, persuadé qu'elle répondra à la haute idée qu'on a de tout ce qui sort de votre plume. Je vous exhorte donc aussi de mon côté à continuer ce travail, et j'applaudis d'avance à vos succès de tout mon cœur.*
Lagrange, letter to Laplace; Lagrange (1776/1892).

On the Continent, the *Doctrine* remained unnoticed until the end of the 19$^{th}$ century.

**54.** Vague ideas on probability and an inexact distinction between subjective and objective probabilities are among the main obstacles against the speedy development of practical medicine. Davidov (1854, p. 66).
**Comment.** Perhaps vague ideas about probability theory in general?

**55.** Will man … eine exakte Theorie … in Angriff nehmen, so muss vor allem die Wahrscheinlichkeit der verschiedene Zustände bestimmt werden welche an einem und demselben Moleküle im Verlaufe einer sehr lange Zeit und an den verschiedenen Molekülen gleichzeitig vorkommen. Boltzmann (1872/1909, p. 317).
**Comment.** The equivalence of these two definitions for gas considered as a whole is established by the ergodic hypothesis.

**56.** I have found it convenient, instead of considering one system of … particles, to consider a large number of systems similar to each other in all respects except in the initial circumstances of the motion, which are supposed to vary from system to system, the total energy being the same in all. In the statistical investigation of the motion, we confine our attention to the number of these systems which at a given time are in a phase such that the variables which define it lie within given limits (p. 715). Boltzmann (1868/1909, § 3) defines the probability of the system being in [a certain] phase … as the ratio of the aggregate time during which it is in that phase to the whole time of the motion, the whole time being supposed to be very great. I prefer to suppose that there are a great many systems the properties of which are the same and that each of these is set in motion with a different set of values for [the coordinates and momenta] the value of the total energy … being the same in all, and to consider the number of these systems, which, at a given instant, are in this phase … (p. 721). Maxwell (1879/1890, pp. 715, 721).
**Comment.** Boltzmann's definition demanded the application of geometric probabilities. Concerning Boltzmann's work in general, Brush (1976, p. 243) remarked that he "often succumbed to unnecessary verbosity". And here is Maxwell, letter of 1873; Knott (1911, p. 114): "By the study of Boltzmann I have been unable to understand him. He could not understand me on account of my shortness, and his length was and is an equal stumbling block to me …". Indeed, "Maxwell's Deduktion [1867/1890] wegen ihre großen Kürze schwer verständlich ist". Boltzmann (1868/1909, p. 49).

**57.** Um also eine zuverlässige Grundlage für die Behandlung solcher Probleme [of the kinetic theory] zu gewinnen, scheint mir die Forderung unabweisbar, eine bestimmte Definition der Wahrscheinlichkeit an die Spitze zu stellen, die, wenn auch in gewissen Grade willkürlich, doch im Verlaufe der Untersuchung nicht mehr geändert oder durch neue Annahmen ergänzt werden darf, und man wird ferner unbedingt festhalten müssen an dem Laplaceschen Wahrscheinlichkeitssatze, nach welchem zwei notwendig wie Ursache und Wirkung mit einander verbundene Ereignisse, … sodass das Eintreten des eines von ihnen das des anderen bedingt, auch immer gleich wahrscheinlich sein müssen. Zermelo (1900, p. 318).

**Comment.** Laplace hardly ever formulated the indicated theorem. Zermelo apparently had in mind that if $P(A/B) = P(B/A) = 1$, then

$$P(AB) = P(A)P(B/A) = P(A),$$
$$P(AB) = P(B)P(A/B) = P(B), P(A) = P(B).$$

**58.** Bei Berechnung der Entropie auf molekulartheoretischem Wege wird häufig das Wort "Wahrscheinlichkeit" in einer Bedeutung angewendet, die sich nicht mit der Definition der Wahrscheinlichkeit deckt, wie sie in der Wahrscheinlichkeits-rechnung gegeben wird. Insbesondere werden die "Fälle gleicher Wahrscheinlichkeit" häufig hypothetisch festgesetzt in Fällen, wo die angewendeten theoretischen Bilder bestimmt genug sind, um statt jener hypothetischen Festsetzung eine Deduktion zu geben. Ich will in einer besonderen Arbeit zeigen, dass man bei Betrachtungen über thermische Vorgänge mit der so genannten "statistischen Wahrscheinlichkeit" vollkommen auskommt und hoffe dadurch eine logische Schwierigkeit zu beseitigen, welche der Durchführung des Boltzmannschen Prinzips noch im Wege steht. Einstein (1905/1989, p. 158).

**59.** Wenn man sich auf den Standpunkt stellt, dass die Nichtumkehrbarkeit der Naturvorgänge nur eine *scheinbare* ist, und dass der nichtumkehrbare Vorgang in einem Übergang zu einem wahrscheinlicheren Zustand bestehe, so muss man zunächst eine Definition der Wahrscheinlichkeit … eines Zustand geben. Die einzige solche Definition, die nach meiner Meinung in Betracht kommen kann, wäre die folgende: Es seien $A_1, A_2, …, A_l$ alle Zustände, welche ein nach außen abgeschlossenes System bei bestimmtem Energieeinhalt anzunehmen vermag, bzw. genauer gesagt, alle Zustände, welche wir an einem solchen System mit gewissen Hilfsmitteln zu unterscheiden vermögen. Nach der klassischen Theorie nimmt das System nach einer bestimmten Zeit einen bestimmten dieser Zustände … an, und verharrt darauf in diesem Zustand (thermodynamisches Gleichgewicht). Nach der statistischen Theorie nimmt aber das System in regelmäßiger Folge alle Zustände … immer wieder an. Beobachtet man das System eine sehr lange Zeit θ hindurch, so wird es einen gewissen Teil $τ_v$ dieser Zeit geben, so dass das System während $τ_v$ und zwar nur während $τ_v$ den Zustand $A_v$ inne hat. Es wird $τ_v/θ$ einen bestimmten Grenzwert besitzen, den wir die Wahrscheinlichkeit … des betreffenden Zustandes $A_v$ nennen. Einstein (1909/1989, p. 544).

**60.** Verallgemeinernd können wir den Satz aussprechen: Die Zustandsänderungen eines isolierten Systems erfolgen derart, dass (im Durchschnitt) wahrscheinlichere Zustände auf unwahrscheinlichere folgen. Man sieht, dass die Wahrscheinlichkeit eines Zustandes eine fundamentale thermodynamische Bedeutung haben muss. Einstein (1925/1995, p. 532).

**61.** Von Wahrscheinlichkeit kann erst gesprochen werden, wenn ein wohlbestimmtes, genau umgrenztes Kollektiv vorliegt. Kollektiv ist eine Massenerscheinung oder ein Wiederholungsvorgang, bei Erfüllung von zwei Forderungen, nämlich: es müssen die relativen Häufigkeiten der einzelnen Merkmale bestimmte Grenzwerte besitzen und diese müssen ungeändert bleiben, wenn man durch willkürliche Stellenauswahl einen Teil der Elemente aus der Gesamtheit heraushebt (p. 29). Die Aufgabe der Wahrscheinlichkeitsrechnung, und zwar die einzige, die ihr in dem vorliegenden Falle zukommt, ist es nun, aus den Wahrscheinlichkeiten, die innerhalb der beiden erstgenannten Kollektivs bestehen, die Wahrscheinlichkeiten innerhalb des abgeleiteten, dritten Kollektivs zu berechnen (pp. 30 – 31). Die Aufgabe der Wahrscheinlichkeitsrechnung ist es, aus den als gegeben vorausgesetzten Verteilungen innerhalb der Ausgangskollektivs die Verteilung innerhalb der abgeleiteten Kollektivs zu berechnen (p. 36). Mises (1928/1972, pp. 29, 30 – 31, 36).

**62.** Dans toutes ces questions [of the kinetic theory] la difficulté principale est, comme nous le verrons, de donner une définition correcte et claire de la probabilité. Langevin (1913/1914, p. 3).
**Comment.** Recall the phase and the time average probabilities and the ergodic hypothesis.

**63.** Um das Jahr 1910 konnte man in Göttingen das Bonmot hören: "Die mathematische Wahrscheinlichkeit ist ein Zahl, die zwischen Null und Eins liegt und über die man sonst nichts weis". Kamke (1933, p. 14).

**64.** We [mathematicians] are content with notions concerning force, space, time, probability, etc. Philosophers have written full volumes about these same subjects of no use for physicists or mathematicians. From our point of view, Mill, Kant and others are not better, but worse than Aristotle, Plato, Descartes, Leibniz … Nekrasov. Marginal note on p. 4 of Chuprov (1896).
**Comment.** At that time, Nekrasov was still a normal person, cf. No. 728, Comment. His viewpoint would hardly be generally shared, and in any case it did not acknowledge the lack of a real definition of probability.

**65.** The notion of equiprobable is comprehensible to all, and, like for example, the concept of real space, does not need a logical definition. We discuss the equiprobability of two or several events … just like we mention … the linearity of a ray of light. Glivenko (1939, p. 18).

**66.** The idea of equiprobability appears not as a formal logical base of the doctrine of mass phenomena; in single concrete situations, it rather is our sole method of theoretically forecasting the probabilities of events. Khinchin (1961/2004, pp. 420 – 421).

**67.** Comment expliquer la vitalité, la puissance de cette théorie quand on la voit s'édifier sur une base si contestée, si fragile, où l'on ne s'accorde même pas sur la signification de ce qu'on calcule? Au début du Calcul des Probabilités (comme à ceux des diverses sciences mathématiques), il n'y a eu aucune démarcation nette entre le concret et l'abstrait. Jusqu'à une époque assez récente, les plus importants ouvrages de Calcul des Probabilités sont restés sur ce terrain intermédiaire où les raisonnements semblent s'appliquer directement à des phénomènes réels mais où l'on attribue à ceux-ci un caractère de simplicité qu'ils ne possèdent pas réellement. Fréchet (1946, pp. 129 – 130).

**68.** Probability, in mathematics, a numerical characteristic of the degree of possibility of the occurrence of some specific event under one or another specific condition repeatable an unlimited number of times. As a category of scientific knowledge, the concept of probability reflects a particular type of connection between phenomena that are characteristic of large-scale processes. The category of probability forms the basis of a particular class of regularities – probability or statistical regularities. … Mathematical probability is an expression of the qualitatively distinctive connection between the random and the necessary. Kolmogorov (1951/1974, p. 423).

**69.** We may fairly lament that intuitive probability is insufficient for scientific purposes, but it is a historical fact. Feller (1950/1957, p. 5).

**Comment.** According to the context, this was inevitable.

**70.** It must be stressed that many of the most essential results of mathematical probability have been suggested by the non-mathematical context of real world probability, which has never even had a universally acceptable definition. In fact the relation between real world probability and mathematical probability has been simultaneously the bane of and inspiration for the development of mathematical probability. Doob (1994/1996, p. 587).

**71.** Each author … invariably reasoned about equipossible and favourable chances, attempting, however, to leave this unpleasant subject as soon as possible. Khinchin (1961/2004, p. 396).

**72.** In the first of his works, Mises some twenty years ago raised the alarm in connection with the scandalously unhappy situation concerning the foundations of the theory of probability. This single fact already constitutes an historical merit of such importance that one can forgive much. Ibidem.

**73.** Die prinzipiellen Probleme des Wahrscheinlichkeitsbegriffes eng mit der Informationstheorie zusammenhängen. (Zum Beispiel spielt in der Diskussion um die Frage der Objektivität oder

Subjektivität die so genannte Bayessche Methode eine zentrale Rolle.) Rényi (1969/1969, p. 87).

**74.** La probabilité d'un événement est la raison que nous avons de croire qu'il aura ou qu'il a eu lieu (p. 30). On rapportera … le mot *chance* aux événements en eux-mêmes et indépendamment de la connaissance que nous en avons, et l'on conservera au mot *probabilité* sa définition précédente (p. 31). Poisson (1837, pp. 30, 31).

**75.** J'ai mis entre les mots *chance* et *probabilité* la même différence que vous, et j'ai beaucoup insisté sur cette différence. Poisson, letter of 1836 to Cournot; Cournot (1843/1984, p. 6).

**76.** Rien n'est plus important que de distinguer soigneusement la double acception du terme de *probabilité*, pris tantôt dans un sens objectif, et tantôt dans un sens subjectif, si l'on veut éviter la confusion et l'erreur. Ibidem, § 240-4.

**77.** By *chance* I mean the same as probability. Bayes (1764/1970, p. 137).

**78.** Jene Unterscheidung zwischen "objektiven" und "subjektiven" Wahrscheinlichkeiten anerkanntermassen nicht stichhaltig ist, weil jede gegebene Wahrscheinlichkeit einen bestimmten Wissens- oder Unwissenheitszustand voraussetzt und in diesem Sinn notwendig subjektive ist. Bortkiewicz (1894 – 1896, p. 661).

**Comment.** A difference, and not a small one, does nevertheless exist. Chuprov's marginal note written (in German or Russian?) on his copy of Bortkiewicz' paper and now adduced to the Russian translation of the latter, see Chetverikov (1968, p. 74).

**79.** La probabilité d'un événement futur ne change pas lorsque les causes dont il dépend subissent des modifications inconnues. Catalan (1877).

**Comment.** A subjective probability is really meant here. Poisson (1837, p. 37) considered an extraction of a ball from an urn with white and black balls in an unknown proportion and proved that the probability of the ball being white was 1/2 which corresponded with "la perfaite perplexité de notre esprit" (and conformed to the principles of information theory). He also meant the subjective probability.

**80.** Probability may be described, agreeably to general usage, as importing partial incomplete belief. Edgeworth (1884, p. 223/1996, p. 6).

**81.** Incomplete knowledge must be considered as perfectly normal in probability theory; we might even say that, if we knew all the circumstances of the phenomena, there would be no place for probability, and we would know the outcome with certainty. Borel (1950, p. 16).

**82.** Who even before battle gains victory by military estimations has many chances. … Who has many chances gains victory; who has few chances does not gain victory; all the less he who has no chances at all. Burov et al (1972, p. 203).

**Comment.** A statement formulated in China in the 4$^{th}$ century BC. The chances had apparently been partly subjective and in any case not calculated statistically.

**83.** The basic object of the calculus of probability … is the probability of events on separate trials; without this probability there is no law of large numbers. Markov (1911/1981, p. 150).

**84.** Aus der Gesamtheit der Fälle, in denen die verschiedenen Wahrscheinlichkeitswendungen der Umgangssprache gebraucht werden, lässt sich eine Gruppe herauslösen, in der der Wahrscheinlichkeitsbegriff einer Präzisierung fähig ist. So entsteht die sog. Wahrscheinlichkeitsrechnung als eine exakte <u>Theorie der Massenerscheinungen und Wiederholungsvorgänge</u> in demselben Sinn wie die rationelle Mechanik Theorie der Bewegungsvorgänge, die Geometrie Theorie der Raumerscheinungen ist. Mises (1939, p. 181).

**Comment.** This utterance explains, to a certain extent, the author's definition of probability as the limit of frequency.

**85.** The idea of probability essentially belongs to a kind of inference which is repeated indefinitely. An individual inference must be either true or false, and can show no effect of probability; and, therefore, in reference to a single case considered in itself, probability can gave no meaning. Peirce (1878/1958, p. 281).

**86.** 'Twould perhaps be more convenient, in order at once to preserve the common signification of words, and mark the several degrees of evidence, to distinguish human reason into three kinds, viz. that from knowledge, from proofs, and from probabilities. … By probability [I mean] that evidence, which is still attended with uncertainty (p. 124). Our reason must be consider'd as a kind of cause, of which truth is the natural effect; but such-a-one as by irruption of other causes, and by the inconstancy of our mental powers, may frequently be prevented. By this means all knowledge degenerates into probability (p. 180). But knowledge and probability are of such contrary and disagreeing natures, that they cannot well run insensibly into each other, and that because they will not divide, but must be either entirely present, or entirely absent, therefore all knowledge resolves itself into probability (p. 181). Hume (1739/1969, Book 1, pt. 3, §11, p. 124; pt. 4, §1, pp. 180, 181).

**87.** Avoir $p$ chances d'obtenir $a$ et $q$ d'obtenir $b$, les chances étant équivalentes, me vaut $(pa + qb)/(p + q)$. Huygens (1657, Proposition 3).

**Comment.** He had not introduced any special term here. Both he and later Jakob Bernoulli (1713, pt. 1, Comment on that Proposition) justified that formula by reasonable arguments. Laplace (1814/1886, p. XVIII) considered such formulas as a definition of the expectation, and so it has remained since.

**88.** The Greek Chronologers … have made the kings of their several Cities … to reign about 35 or 40 years a-piece, one with another; which is a length so much beyond the course of nature, as is not to be credited. For by the ordinary course of nature Kings Reign, one with another, about 18 or 20 years a-piece; and if in some instances they reign, one with another, five or six years longer, in others they Reign as much shorter: 18 or 20 years is a medium. Newton (1728, p. 52).

**Comment.** The length of the reign of a king depends on that of his predecessor. An estimate of the appropriate expected value is therefore difficult the more so as the appropriate variance is unknown. Newton based the interval of 18 or 20 years on trustworthy sources.

**89.** Man lege jedem Fehler ein von seiner Größe abhängendes Moment bei, multiplicire das Moment jedes möglichen Fehlers in dessen Wahrscheinlichkeit und addire die Produkte: der Fehler, dessen Moment diesem Aggregat gleich ist, wird als mittlerer betrachtet werden müssen. Allein welche Funktion der Größe des Fehlers wir für dessen Moment wählen wollen, bleibt wieder unserer *Willkür* überlassen, wenn nur der Werth derselben immer positiv ist, und für größere Fehler größer als für kleinere. Der Verf. hat die einfachste Funktion dieser Art gewählt, nämlich das Quadrat; diese Wahl ist aber noch mit manchen anderen höchst wesentlichen Vortheilen verknüpft, die bei keiner anderen stattfinden. Gauss (1821/1887, p. 192).

**Comment.** The variance appeared in the theory of errors.

**90.** L'analyse est le seul instrument dont on se soit servi jusqu'à ce jour dans la science des probabilités, pour déterminer & fixer les rapports du hasard; la Géométrie paroissoit peu propre à un ouvrage aussi délié; cependant si l'on y regarde de près, il sera facile de reconnoître que cet avantage de l'Analyse sur la Géométrie, est tout-à-fait accidentel, & que le hasard selon qu'il est modifié & conditionné, se trouve du ressort de la géométrie aussi bien que de celui de l'analyse; pour s'en assurer, il suffira de faire attention que les jeux & les questions de conjecture ne roulent ordinairement que sur des rapports de quantités discrètes; l'esprit humain plus familier avec les nombres qu'avec les mesures de l'étendue les a toujours préférés; les jeux en sont une preuve, car leurs loix sont une arithmétique continuelle; pour mettre donc la Géométrie en possession de ses droits sur la science du hasard, il ne s'agit que d'inventer des jeux qui roulent sur l'étendue & sur ses rapports, ou calculer le petit nombre de ceux de cette nature qui sont déjà trouvés. Buffon (1777/1954, p. 471).

**Comment.** Already in 1735 an anonymous author (Buffon himself?) described Buffon's successful attempt. In essence, the first to apply geometric probability was Newton. See Sheynin (2003b) for the history of this notion.

**91.** On pourrait faire usage de Calcul des Probabilités pour rectifier les courbes ou carrer leurs surfaces. Sans doute, les géomètres n'emploieront pas ce moyen … Laplace (1812/1886, p. 365).

**Comment.** Laplace (Ibidem, p. 342) indirectly made known how to understand a surface of a curve: "L'ordonnée qui divise l'aire de la courbe en parties égales …", and he several times separated himself from "les géomètres". He devoted his next pages to geometrical probability, and although he did not say so directly, it followed that the quotation above concerned what is now called statistical simulation.

**92.** Let us examine what is probable would [would probably] have been the least apparent distance of any two or more stars, anywhere in the whole heavens, upon the supposition that they have been scattered by mere chance, as it might happen. Now it is manifest, upon this supposition, that every star being as likely to be in any one situation as another, the probability, that any particular star should happen to be within a certain distance … of any other given star, would be represented … by a fraction … Michell (1767) as quoted by Hardin (1966, p. 38).

**Comment.** Michell's problem continued to be discussed until the 20th century. In particular, some considerations about what we would now call possible deviations of an empirical distribution from its corresponding theoretical law had been voiced, see Proctor (1874, p. 99), Kleiber (1887, p. 439) and Sheynin (2009a, § 6.1.6). W. Herschel, in 1802, attempted to solve a similar problem. His calculations contained a mistake, but he correctly concluded that "casual situations will not account for the multiplied phenomena of double stars… Their existence must be owing to the influence of some general law of nature", see Sheynin (1984a, p. 162).

**93.** Les probabilités relatives à la distribution des étoiles, en les supposant semées au hasard sur la sphère céleste, sont impossibles à assigner si la question n'est pas précisée davantage (p. 7). L'ingénieux argument de Mitchell (!) ne peut cependant fournir d'évaluation numérique (p. 170). Bertrand (1888a, pp. 7, 170).

**Comment.** Other peculiarities of the sidereal system rather than the short distance between stars could have been evaluated as well or instead.

**94.** L'infini n'est pas un nombre; on ne doit pas, sans explication, l'introduire dans les raisonnements. La précision illusoire des mots pourrait faire naître des contradictions. Choisir *au hasard*, entre un nombre infini de cas possibles, n'est pas une indication suffisante. … On trace *au hasard* une corde dans un cercle. Quelle est la probabilité pour qu'elle soit plus petite que le côté du triangle équilatéral inscrit? Bertrand (1888a, p. 4).

**Comment.** Bertrand specified his problem in three differing ways by natural conditions and obtained three different answers. It was proved later that the problem had an uncountable number of solutions. See Sheynin (2003b) for its history and the application of geometric probability by mathematicians and natural scientists.

**95.** Par des raisonnements qui peuvent paraître également plausibles, il [Bertrand] trouve pour la probabilité cherchée deux valeurs différentes, tantôt 1/2, tantôt 1/3. Cette question l'a préoccupé; il en avait trouvé la solution, mais il la laisse à chercher à son lecteur. Darboux (1902/1912, p. 50).

**Comment.** Darboux mistakenly left out the third solution and did not elaborate.

**96.** [Probability is] le rapport de l'étendue des chances favorables à un événement, à l'étendue totale des chances. Cournot (1843, §18).

**Comment.** Cournot introduced this definition (in which *étendue* ought to be now replaced by *measure*) both for the discrete and the continuous cases. With regard to geometric probability, every author beginning with Newton had effectively applied it.

**97.** Supposons qu'une chose quelconque, que nous nommerons A pour abréger le discours, soit susceptible, par sa nature, de toutes les valeurs comprises entre des limites données … Soit *x* une de ces valeurs; si l'on fait une série d'expériences pour déterminer A, la probabilité que la valeur qu'on trouvera par une de ces observations n'excédera pas *x*, variera, en général, d'une expérience à l'autre. Nous la représenterons par $F_n x$ pour la $n^{ième}$ expérience. La probabilité que cette valeur sera précisément *x* ne pourra être qu'infiniment petite, puisque le nombre des valeurs possibles est infini; en faisant $dF_n x/dx = f_n x$, elle aura pour expression $f_n x dx$. Poisson (1829, p. 3).

**Comment.** To *determine* A is of course not quite proper. The first to introduce directly (not only effectively) a density curve was Simpson (see No. 655). However, Huygens (1669, correspondence 1888 – 1950, 1895, plate between pp. 530 and 531) drew a graph of an arbitrary smooth curve passing through the empirical points given by Graunt's mortality table and featuring the function $y = 1 - F(x)$ where, in modern notation, $F(x)$ was a distribution function and $0 \leq x \leq 100$.

**98.** Ich fange mit der Aufsuchung der Wahrscheinlichkeit = $\varphi x dx$ an, mit welcher erwartet werden kann, dass ein Beobachtungsfehler zwischen *x* und $x + dx$ falle, wenn er, auf gegebene Art, von einer Ursache ξ abhängt, für welche jeder zwischen zwei Grenzen – α und α liegende Werth gleich möglich ist. Bessel (1838/1876, p. 374).

**Comment.** Bessel calculates the density of a random variable functionally connected with another such variable with a given distribution. This problem is now standard, and it was repeatedly solved even before Bessel, but he was the first to formulate it in such a manner that the densities became mathematical objects.

**99.** Everywhere possible, I exclude the completely undefined expressions *random* and *at random*; where it is necessary to use them, I introduce an explanation corresponding to the pertinent case. Markov, letter of 1912 to Chuprov (Ondar 1977/1981, p. 65).

**Comment.** Recalling the Bertrand paradox, it is easy to understand Markov's apprehension. However, he could have applied the Poisson heuristic definition of random variable and moreover made the next step by introducing an appropriate notation for it. And in spite of his statement, Markov sometimes wrote *indefinite* without any explanation (which the translators of his correspondence with Chuprov, no doubt attempting to modernize the exposition, rendered as *random*) and this was definitely bad. Note that Markov never applied the then already usual terms *normal distribution* or even *coefficient of correlation*.

### 2.3. Logical Difficulties and Mistaken Opinions

**100.** Ich weiß so gut wie du, dass die Argumente, die aufgrund von Plausibilitäten etwas beweisen, ihren Mund zu voll nehmen, und dass sie einen, wenn man sich nicht in acht nimmt, ganz schön in die Irre führen, in der Geometrie ebenso wie in allen anderen Gebieten. Plato (2004, 92d, p. 56).

**101.** La théorie des Probabilités tient à des considérations si délicates, qu'il n'est pas surprenant qu'avec les mêmes données deux personnes trouvent des résultats différents, surtout dans les questions très compliquées. Laplace (1814/1886, p. XI).

**102.** Le calcul des probabilités donne lieu à de singulières illusions, auxquelles les meilleurs esprits n'ont pu toujours se soustraire. Bienaymé (1853, p. 310).

**103.** Man verwechsle nicht einen unvollständig bewiesenen Satz, dessen Richtigkeit infolgedessen problematisch ist, mit einem vollständig erwiesenen Satze der Wahrscheinlichkeitsrechnung; letzter stellt, wie das Resultat jedes anderen Kalküls, eine notwendige Konsequenz gewisser Prämissen dar, und bestätigt sich, sobald diese richtig sind, ebenso in der Erfahrung. … Nur scheint es hier doppelt geboten, den Schlüssen mit der größten Strenge zu verfahren. Boltzmann (1872/1909, p. 317).

**104.** No part of mathematical science requires more careful handling than that which treats of probabilities and averages. … Even the validity of certain methods of proof is still apparently an open question. Maxwell (1877, p. 242).

**105.** This branch of mathematics is the only one, I believe, in which good writers frequently get results entirely erroneous. … It may be doubted if there is a single extensive treatise on probabilities in existence which does not contain solutions absolutely indefensible. Peirce (1878/1958, p. 279).

**106.** The analysis of probabilities considers and quantitatively estimates even such phenomena … which, due to our ignorance, are not subject to any suppositions. Buniakovsky (1846, p. I).

**Comment.** In actual fact Buniakovsky never considered such phenomena and on p. 364 and elsewhere (1866, p. 24) he went back on his opinion (lacking in Laplace's works).

In Russia, he was the main predecessor of Chebyshev. Here are opinions about his treatise (1846). F. K. Derschau, obviously the author of Anonymous (1847, p. 44):

*Buniakovsky's book is new in subject matter, complete in its contents and scientifically up-to date. What else shall we demand of an author* [who does not claim originality]*?*

Markov (1914/1981, p. 162): "A beautiful work". Steklov (1924, p. 177): "For his time, [he had compiled] a complete and outstanding treatise". Struve (1918) called him a "Russian disciple of the French mathematical school" and Kolmogorov (1947, p. 56/2005, p. 72) mentioned the "excellent for his time treatise".

**107.** Mere ignorance is no ground for any inference whatever. … Ex nihilo nihil. Ellis (1850/1863, p. 57).

**108.** Il est aussi faisable de jetter douze points, que d'en jetter onze, car l'un et l'autre ne se peut faire que d'une seul manière; mais il est trois fois plus faisable d'en jetter sept, parce que cela se peut faire en jettant 6 et 1, 5 et 2, et 4 et 3; et une combinaison icy est aussi faisable que l'autre. … Feu Monsieur Bernoulli a cultivé cette matière sur mes exhortations. Leibniz (1714/1887, pp. 569 – 570).

**Comment.** Todhunter (1865, p. 48) had not failed to notice that both conclusions were wrong. Then, Leibniz' statement about Jakob Bernoulli was misleading: the latter had begun to *cultivér* probability theory long before Leibniz formulated his appropriate wish (in 1703, see No. 6), and not later than in 1690 had sketched the proof of the law of large numbers in his Diary (*Meditationes*).

**109.** On demande combien il y a à parier qu'on amener *croix* en jouant deux coups consécutifs. La réponse qu'on trouvera dans tout les auteurs, & suivant les principes ordinaires, est celle-ci: … Cependant cela est-il bien exact? Car pour ne prendre ici que le cas de deux coups, ne faut-il pas réduire à une les deux combinaisons qui donnent *croix* au premier coup? Car dès qu'une fois *croix* est venu, le jeu est fini, & le second coup est compté pour rien. Ainsi il n'y a proprement que trois combinaisons de possibles … Donc il n'y a que 2 contre 1 à parier. D'Alembert (1754/1966).

**Comment.** The game ends after the first throw with probability 1/2, and before it begins the success of the second attempt has probability 1/4. The total probability is 3/4 and the correct odds are 3:1.

**110.** Vous me demanderez peut-être quels sont les principes qu'il faut, selon moi, substituer à ceux dont je révoque en doute l'exactitude? Ma réponse sera celle que j'ai déjà faite; je n'en sais rien, et je suis même très-porté à croire que la matière dont il s'agit [the theory of probability], ne peut être soumise, au moins à plusieurs égards, à un calcul exact et précis, également net dans ses principes et dans ses résultats. D'Alembert (1767b, pp. 309 – 310).

**Comment.** Here are D'Alembert's hardly satisfactory pertinent remarks (Todhunter 1865, pp. 265, 280, 286) made in 1761, 1768 and 1780 respectively.

   1) His solution of the problem described in No. 109, *faute d'en connôitre une meilleure*, was not at all *éxacte à la rigueur*, but the *régle ordinaire* is still worse. Then, *pour parvenir à une théorie satisfaisante* it is necessary to know how to find the true ratio of unequal probabilities [?]; to understand when a probability might be neglected; and to estimate expectations *selon que la probabilité est plus ou moins grande*.

   2) He does not wish to ruin the calculus of probability, he only desires *qu'il soit éclairci et modifié*.

   3) Some principles of that calculus ought to be either changed, or at least expounded *de manière à ne plus laisser aucun nuage*.

**111.** Da er [D'Alembert] seines unerträglichen Hochmuths ungeachtet längst begreiffen konnte, dass er sich dazu [to a lucrative position in Russia] gar nicht schickte. Sonsten besteht seine Philosophie nach des H. [errn?] Bernoullis Ausdruck in einer impertinente suffisance, dass er alle seine Fehler auf das unverschämteste zu vertheidigen sucht, welches ihm nur gar zu oft wiederfahren, so dass er seit vielen Jahren aus Verdruss nichts mehr mit der Mathematic zu thun haben will. In seiner Hydrodynamic hat er die meisten Sätze des H. [errn?] Bernoullis, die doch durch häuffige Erfahrung bestätigt waren, gantz cavalierement widerlegt, da doch seine eigenen der Erfahrung widersprechen. Und er hat seinen Hochmuth noch nicht so weit überwinden können, dass er seinen offenbaren Irrthum hätte bekennen wollen. Euler, letter of 1763; Juskevic et al (1959, p. 221).

**112.** Mir schrecken [bei solchen Untersuchungen] die Einwürfe D'Alembert's nicht zurück, der diesen Kalkül zu verdächten versucht hat. Zuerst nämlich hat dieser bedeutende Geometer die mathematischen Studien beiseite gelegt; jetzt scheint er sie sogar zu bekämpfen, indem er unternommen hat, eine Reihe von Grundsätzen umzustürzen, die auf das sicherste begründet sind. Den Laien mögen seine Einwürfe gewichtig erscheinen, doch die Furcht liegt fern, dass die Wissenschaft selbst Schaden durch sie erleide. Euler (1785, in Latin, p. 408). Translated by Netto (1908, pp. 226 – 227).

**113.** The negative probability may no doubt be an index of the removal from possibility of the circumstances, or of the alteration of data which must take place before possibility begins. But I have not yet seen a problem in which such interpretation was worth looking for. I have, however, stumbled upon the necessity of interpretation at the other end of the scale: as in a problem in which the chance of an event happening turns out to be 2.5; meaning that under the given hypothesis the event must happen twice, with an even chance of happening a third time. De Morgan (1864, p. 421n).

**Comment.** Rice & Seneta (2005) studied De Morgan's work, referred to that memoir misdating it and mistakenly described its contents in a single phrase. They never noticed his statement quoted above. And here is De Morgan's opinion from his letter to John Herschel of 1842 (De Morgan, Sophia 1882, p. 147): $\sin\infty = 0$, $\cos\infty = 0$, $\tan\infty = \mp\sqrt{-1}$, $\cot\infty = \mp\sqrt{-1}$. And he suspects that $\sec\infty = \csc\infty = 0$ or $\infty$. I did not find Herschel's reply.

**114.** Le calcul des probabilités ne me semble avoir été réellement, pour ses illustres inventeurs, qu'un texte commode à d'ingénieux et difficiles problèmes numériques, qui n'en conservent pas moins toute leur valeur abstraite, comme les théories analytiques dont il a été ensuite l'occasion, ou, si l'on veut, l'origine. Quant à la conception philosophique sur laquelle repose une telle doctrine, je la crois radicalement fausse et susceptible de conduire aux plus absurdes conséquences. Je ne parle pas seulement de l'application évidemment illusoire qu'on a souvent tenté d'en faire au prétendu perfectionnement des sciences sociales: ces essais, nécessairement chimériques, seront caractérisés dans la dernière partie de cet ouvrage. C'est la notion fondamentale de la probabilité évaluée, qui me semble directement irrationnelle et même sophistique: je la regarde comme essentiellement impropre à régler notre conduite en aucun cas, si ce n'est tout au plus dans les jeux de hasard. Elle nous amènerait habituellement, dans la pratique, à rejeter, comme numériquement invraisemblable, des événements qui vont pourtant s'accomplir. On s'y propose le problème insoluble de suppléer à la suspension de jugement, si nécessaire en tant d'occasions. Les applications utiles qui semblent lui être dues, le simple bon sens, dont cette doctrine a souvent faussé les aperçus, les avait toujours clairement indiquées d'avance. Comte (1830 – 1842/1877, t. 2, p. 255).

**Comment.** It is hardly necessary to refute minutely this opinion of a most prominent opponent of the applications of probability theory. Note, however, that statements as *propose le problème insoluble de suppléer à la suspension* are certainly libellous. Perhaps Comte thought about the Petersburg paradox (see No. 115) which Niklaus Bernoulli had apparently invented to strengthen by its study the logical basis of the theory of probability. Nevertheless, Vasiliev (Bazhanov 2002, p. 131) positively appraised Comte's general attitude towards mathematics.

### 2.4. Moral Expectation

**115.** A promet de donner un écu à B, si avec un dé ordinaire il amené au premier coup six points, deux écus s'il amené le six [only] au second, trois écus s'il amené ce point au troisième coup … On demande quelle est l'espérance de B. … On demande la même chose si A promet à B de lui donner des écus en cette progression 1, 2, 4, 8, 16 &c, ou 1, 3, 9, 27, &c … Niklaus Bernoulli, letter of 1713 to Montmort; Montmort (1713/1980, p. 402).

**Comment.** In 1738, Daniel Bernoulli studied this game (No. 116). He published his memoir (and many others) in Petersburg, hence the name of the game. Its conditions were insignificantly changed: the progression 1, 2, 4, … persisted, but the die was replaced by a coin. On the further history of the game see Sheynin (2009a, § 3.3.4). In particular, Freudenthal (1951) recommended to consider a series of such games with the role of the gamblers in each of them to be decided by lot.

**116.** No characteristic of the persons themselves ought to be taken into consideration; only those matters should be weighed carefully that pertain to the terms of the risk (§ 2). But anyone who considers the problem with perspicacity and interest will ascertain that the concept of *value* which we have used in this rule may be defined in a way which renders the entire procedure universally acceptable without reservation. To do this the determination of the *value* of an item must not be based on its *price*, but rather on the *utility* it yields (§ 3).

In many games, even those that are absolutely fair, both of the players may expect to suffer a loss; indeed, this is Nature's admonition to avoid the dice altogether (§ 13). D. Bernoulli (1738, §§ 2, 3, 13; translation 1954).

**Comment.** This was the author's starting point which led him to the introduction of the concept of moral expectation; the term itself was introduced in 1732 by Gabriel Cramer in his letter to Niklaus Bernoulli and quoted by Daniel. The immediate cause of Daniel's study was the Petersburg paradox, see No. 115.

**117.** Le jeu mathématiquement le plus égal est toujours désavantageux. Laplace (1812/1886, p. 443).

**118.** L'écu que le pauvre a mis a part pour payer un impôt de nécessité, & l'écu qui complète les sacs d'un financier, n'ont pour l'Avare & pour le Mathématicien que la même valeur, celui-ci les comptera par deux unités égales, l'autre se les appropriera avec un plaisir égal, au lieu que l'homme sensé comptera l'écu du pauvre pour un louis, & l'écu du financier pour un liard. Buffon (1777/1954, p. 469).

**119.** Nous nommerons cet ["usual"] avantage *espérance mathématique*, pour le distinguer de l'espérance morale qui dépend, comme lui, du bien espéré et de la probabilité de l'obtenir, mais qui se règle encore sur mille circonstances variables qu'il presque toujours impossible de définir, et plus encore d'assujettir au calcul. Laplace (1812/1886, p. 189).

**Comment.** Laplace's term apparently only took root in the French and the Russian literature.

**120.** La théorie de l'espérance morale est devenue classique, jamais le mot ne pu être plus exactement employé: on l'étudie, on l'enseigne, on la développe dans des livres justement célèbres. Le succès s'arrête là, on n'en a jamais fait en n'en pourra faire aucun usage. Bertrand (1888a, p. 66).

**Comment.** At the end of the 19th century, issuing from Daniel Bernoulli's idea, economists began to develop the theory of marginal utility. Bertrand's opinion was thus refuted.

**121.** M. Ostrogradsky n'admet point l'hypothèse de Daniel Bernoulli; il exprime la fortune morale par une fonction arbitraire de la fortune physique; et il parvient à donner à la solution des questions principales relatives à la fortune morale, toute l'étendue et l'exactitude que l'on peut désirer. Fuss (1836, pp. 24 – 25),

**Comment.** Nothing else is known on this point since Ostrogradsky had not published his study. It may be thought that he was not entirely satisfied with his results.

### 2.5. Law of Large Numbers, Central Limit Theorem

**122.** …da bringe an ein gemein Orth zusammen 1900 alter … Guldenthaler … was nun etwan ein Thaler zu schwer, das ist der ander zu leicht, daß es also auf … 100 Pfund … nichts auffträgt. Kepler (letter dated 1627/1930).

**Comment.** More correctly: the mean weight will not change. Such statements, characteristic of the prehistory of the law of large numbers (Sheynin 2009a, § 3.2.3), had been occurring in the theory of errors even in the second half of the 19th century. In an unpublished and undated manuscript Boscovich (Sheynin 1973, pp. 317 – 318) calculated the chances of sums of random errors for a most elementary example, but did not consider mean values.

**123.** When the purchaser of several life annuities comes to divide his capital … upon several young lives – upon ten, twenty, or more – this annuitant may be assured, without hazard or risk of the enjoyment of [a considerable profit]. De Witt (1671, last lines); translated by Hendriks (1852, pp. 232 – 249).

**Comment.** Condorcet (1785/1972, p. 226) mentioned such purchases.

**124.** [About irregularity in statistical figures provided in his mortality Table: it] seems rather be owing to chance, as are also the other irregularities in the series of age, which would rectify themselves were the number of years [of observation] much more considerable. Halley (1694/1942, p. 5).

**Comment.** Halley intuitively believed in the principle of large numbers, as I would call it. However, the irregularities could have well been occasioned by systematic influences.

**125.** Even the most stupid person … feels sure that the more … observations are taken into account, the less is the danger of straying from the goal. … However, it remains to investigate something that no one had perhaps until now run across even in his thoughts. It certainly remains to inquire whether, when the number of observations … increases, the probability of attaining the real ratio [the theoretical probability] continually augments so that it finally exceeds any given degree of certitude. Or [to the contrary], the problem has, so to say, an asymptote; i. e., that there exists such a degree of certainty which can never be exceeded no matter how the observations be multiplied … To tell the truth, if this [former possibility] failed to happen, it would be necessary to question our attempt at experimentally determining the number of cases [the theoretical probability] (pp. 29 – 30). [His first example: posterior estimation of the unknown ratio of white and black pebbles in an urn. His second example:] If we replace the urn … by air or by a human body which contain … sources of various changes or diseases … we will be able to determine by observation in exactly the same way how much easier can one or another event occur … The ratio between the number of cases … is accepted not as precise and strict … but … with a certain latitude … contained between two limits … (pp. 30 – 31). J. Bernoulli (1713/2005, pp. 29 – 30, 30 – 31).

**Comment.** Bernoulli assumed the possibility of estimating a non-existing probability. Statisticians of the 19$^{th}$ and early 20$^{th}$ century rejected it (Sheynin 2001a, p. 102), and in modern textbooks on probability theory that assumption is not mentioned.

**126.** His solution must be said to be from the practical standpoint a failure … he failed entirely to determine what the *adequate* number of observations were for such limits (p. 202). Bernoulli saw the importance of a certain problem; so did Ptolemy, but it would be rather absurd to call Kepler's or Newton's solution of planetary motion by Ptolemy's name (p. 210)! K. Pearson (1925, pp. 202, 210). A most fundamental principle of statistics has been attributed to Bernoulli instead of its real discoverer De Moivre. K. Pearson (1978, p. 1).

**Comment.** The first statement is wrong because existence theorems are important in themselves. Pearson's second utterance is unfair because Ptolemy issued from the wrong geocentric system of the world whereas Bernoulli did not base himself on any mistaken presumption.

**127.** [In case of need] nous devons suppléer au défaut de la théorie par le moyen de l'observation; elle nous sert à établir des lois empiriques, qui deviennent presque aussi certaines que les lois rationnelles, quand elles reposent sur des observations assez répétées … Cuvier (1812/1861, p. 67).

**128.** Les choses de toutes natures sont soumises à une loi universelle qu'on peut appeler *la loi des grands nombres*. Elle consiste en ce que, si l'on observe des nombres très considérables d'événements d'une même nature, dépendants de causes constantes et de causes qui varient irrégulièrement, tantôt dans un sens, tantôt dans l'autre, c'est-à-dire sans que leur variation soit progressive dans aucun sens déterminé, on trouvera, entre ces nombres, du rapports à très peu près constantes. Poisson (1837, p. 7).

**Comment.** This diffuse definition does not explain the essence of its subject. See Sheynin (2009a, § 8.7).

**129.** Markov (like Chebyshev) attributes [Poisson's law of large numbers] to the case when all the probabilities following one another … are known beforehand. … In concluding, Markov admitted that perhaps there did exist some kind of ambiguity in Poisson's reasoning, but he believed that it was necessary to take into account the later authors' understanding of the term … Bortkiewicz, in a letter of 1897 to Chuprov, about his recent conversation with Markov (Sheynin 1990b/1996, p. 42).

**Comment.** In addition to his reasoning, Poisson offered examples of his law, including the stability of the sea level and the interval between molecules, where no probabilities were known beforehand.

**130.** Cette proposition fondamentale de la théorie des probabilité contenant comme cas particulier la loi de J. Bernoulli, est déduite par Mr Poisson d'un formule, qu'il obtient en calculant approximativement la valeur d'une intégrale definie, assez compliquée. Toute ingénieuse que sont la méthode employée par le célèbre Géomètre, il reste à être impossible de montrer la limite de l'erreur que peut admettre son analyse approximative, et par cette incertitude de la valeur de l'erreur, sa démonstration n'est pas rigoureuse. Chebyshev (1846, p. 259).

**131.** In the very beginning of the calculus of probability there must be a law on which all [its] applications to reality rests. In all justice, this law may be called the law of large numbers. It is independent both of the Bernoulli and the Poisson theorem and is their basis. It reads: If a trial, in which an … event having probability $p$ can occur, is repeated $n$ times with $n$ being sufficiently large, then this event must happen about $np$ times. Romanovsky (1912, p. 22).

**132.** On est donc en droit d'affirmer que la *Loi des grands nombres* n'a pas d'existence réelle, et qu'il ne faut plus s'appuyer de ces mots dans le sens que M. Poisson leur avait donné. Bienaymé (1855, p. 202).

**133.** Poisson (1837, §101) arbeitet schon ganz direkt und recht exakt mit der Fouriertransformation. Wie Laplace begnügte er sich noch mit beschränkten stochastischen variablen. Aber während von De Moivre bis Laplace von den zu addierenden stochastischen Variablen immer vorausgesetzt wurde, dass sie *dieselbe* Verteilung besaßen, ließ er zum ersten Male eine ziemlich willkürliche Reihe von unabhängigen stochastischen Variablen zu. Für sie bewies er dann den zentralen Grenzwertsatz mit Hilfe der Fourier-transformierten*.

   *Kann natürlich bei Laplace, Fourier, Poisson den Mangel an Exaktheit rügen, aber wollte man die Exaktheit zum Kriterium erheben, so müsste man auch die Infinitesimalrechnung historisch

erst bei Weierstrass anfangen lassen. – Übrigens ist auch die übliche Behauptung, die Fouriertransformation von Verteilungsfunktionen rühre von Cauchy her, vollkommen abwegig. Wohl bedingt Cauchy … sich der Fouriertransformation souveräner als Poisson. Freudenthal & Steiner (1966, pp. 171 – 172).

**134.** Each new experiment is as a new stroke of the pencil which bestows an additional vivacity on the colours, without either multiplying or enlarging the figure. Hume (1739/1969, Book 1, pt 3, § 12, p. 135).

**135.** An ingenious Friend has communicated to me a Solution of the inverse [with regard to De Moivre's finding to whom the author referred a few lines above] Problem, in which he has shewn what the Expectation is, when an Event has happened $p$ times, and failed $q$ times, that the original ratio of the Causes for the Happening or Failing of an Event should deviate in any given Degree from that of $p$ to $q$. And it appears from this Solution, that where the Number of Trials is very great, the Deviation must be inconsiderable: Which shews that we may hope to determine the Proportions, and, by degrees, the whole Nature, of unknown Causes, by a sufficient Observation of their Effects. D. Hartley (1749, p. 339).

**Comment.** This rather recently noticed passage had become the starting point of a discussion about the actual author of the Bayes posthumous memoir. Apparently, however, the assumption of any other authorship of that memoir was refuted, see Sheynin (2003a).

**136.** Bayes … paraît avoir possédé à un très-haut degré les qualités du géomètre. Bienaymé (1838, p. 514).

**137.** Let us imagine … a person just brought forth into this world, and left to collect from his observation of the order and course of events what powers and causes take place in it. The Sun would, probably, be the first object that would engage his attention; but after losing it the first night he would be entirely ignorant whether he should ever see it again. [What would be for him the odds in favour of its return after regularly observing the sunrise a million times.] R. Price; Bayes (1764/1970, p. 150).

**138.** One wou'd appear ridiculous, who wou'd say, that 'tis only probable the sun will rise tomorrow. Hume (1739/1969, Book 1, pt. 3, § 11, p. 124).

**139.** Si un événement peut être produit par un nombre $n$ de causes différentes, les probabilités de l'existence de ces causes prises de l'événement sont entre elles comme les probabilités de l'événement prises de ces causes, et la probabilité de l'existence de chacune d'elles est égale à la probabilité de l'événement prise de cette cause, divisée par la somme de toutes les probabilités de l'événement prises de chacune de ces causes. Laplace (1774/1891, p. 29).

**Comment.** This is the principle of inverse probability. Cf. No. 563.

**140.** [Concerning Laplace (1774)] Les remarques que vous faites sur l'aberration de la théorie ordinaire, lorsqu'on regarde communément comme également probables, m'ont paru aussi justes qu'ingénieuses; c'est une nouvelle branche très importante que vous ajouter à la théorie des hasards, et qui était nécessaire pour mettre cette théorie à l'abri de toute atteinte … Lagrange (1775/1892, p. 58).

**141.** La probabilité de l'existence d'une quelconque de ces causes est donc une fraction, dont le numérateur est la probabilité de l'événement résultante de cette cause, et dont le dénominateur est la somme des probabilités semblables relatives à toutes les causes. Si ces diverses causes, considérées *a priori*, sont inégalement probables, il faut, au lieu de la probabilité de l'événement, résultante de chaque cause, employer le produit de cette probabilité par celle de la cause elle-même. Laplace (1814/1886, pp. XIV – XV).

**Comment.** Laplace verbally formulated the particular and general forms of the so-called Bayes theorem. He hardly realized how difficult it was for a non-mathematically minded reader to understand suchlike definitions. For a similar case see No. 139.

**142.** Wenn man aus der Wahrscheinlichkeit des Eintretens eines Ereignisses auf die Wahrscheinlichkeit der Ursachen zurückschließt, so ist dies ein schlüpfrigen Boden. Gauss, letter of 1845; Biermann (1965).

**143.** Note that the Gauss method based on the law of hypotheses does not demand that [the number of observations] be certainly very large. … We saw, however, that this method is not really reliable. Chebyshev (1879 – 1880/1936, p. 250; translation 2004, p. 229).

**Comment.** The *law of hypotheses*, as Chebyshev called it, was the principle of inverse probability. Gauss only applied it in 1809, in his first justification of the principle of least squares.

**144.** Bayes' theorem has come back from the cemetery to which it had been consigned and is today an object of research and application among a growing body of statisticians. Cornfield (1967, p. 41).

**Comment.** The main cause of the dissatisfaction with the *Bayes approach* was, and is, the lack or insufficient knowledge of the prior (law of) probability of the studied event.

**145.** In our time the Bayesian approach is experiencing its second youth. This is partly connected with the decision theory offered by Wald. Tiurin (1975, p. 56).

**146.** A Bayesian assigns probabilities to events the non-Bayesian would not include in his probability model and the statistical tests of the two may differ. The non-Bayesians assert that the Bayesians' very principles are contradictory. The Bayesian asserts that their opponents are inefficient statisticians. Doob (1976, p. 202).

**147.** [Price knew] of no person who has shown how to deduce the solution of the converse problem. … What Mr De Moivre has done therefore cannot be thought sufficient. R. Price (Bayes 1764/1970, p. 135).

**Comment.** Price communicated the memoir of Bayes, his late friend and colleague, and compiled a covering letter published as its Introduction. De Moivre's theorem estimated the statistical probability given its theoretical counterpart and an increasing number of Bernoulli trials. The random variables studied in both cases (by De Moivre and Bayes) had differing variances. It was Gauss who introduced the notion of variance (cf. No. 89), so that Bayes had a deep understanding of his subject.

**148.** Le calcul des fonctions génératrices est le fondement d'un théorie que je me propose de publier bientôt sur les probabilités. Laplace (1811/1898, p. 360).

**149.** De 1770 à 1809 [from 1774 to 1811], pendant près de quarante ans, Laplace avait donné des Mémoires nombreux sur les probabilités: mais, quelque intérêt qu'il y eût dans ces Mémoires, il n'avait pas voulu les rédiger en théorie générale. Aussitôt qu'il a reconnu la propriété des fonctions de probabilités [the central limit theorem], il voit clairement que c'est un principe qui régit presque toutes les applications, et il compose sa théorie [his treatise of 1812]. Bienaymé (1853, pp. 311 – 312).

**150.** La science, quoique Pascal l'ait créée il y a deux siècles, n'est pas bien loin de son origine. Le développement analytique a progressé presque seul. Ibidem, p. 310.

**Comment.** As usual, Fermat is somehow forgotten.

**151.** Cette supposition, la plus naturelle et la plus simple de toutes, résulte de l'emploi du cercle répétiteur dans la mesure des angles des triangles. Laplace (1818/1886, p. 536).

**Comment.** The repeating theodolite substantially reduced the effect of the error of reading. The order of the effects of the two main errors of measurement, of sighting and reading, became equal one to another, hence the supposition of the normal law.

**152.** Je montré que si l'erreur totale résulte de l'accumulation d'un grand nombre de petites erreurs partielles, si ces dernières suivent des lois quelconques, mais symétriques, l'erreur résultante totale sera soumise à la loi de Gauss. Poincaré, written in 1901 (1921, p. 343).

**Comment.** See No. 153.

**153.** [Poincaré] will remain in history as the most penetrating critic of quantitative methods, and the great proponent of qualitative ones. Ekeland (1984/1988, p. 35).

**154.** On the whole my custom of terming the curve the Gauss – Laplacian or normal curve saves us from proportioning the merit of discovery between the two great astronomer mathematicians. K. Pearson (1905, p. 189).

**Comment.** Pearson only became aware of De Moivre's pertinent achievement in 1924.

**155.** [Bessel (1838) provided a non-rigorous proof of a version of the central limit theorem, and on 28 February 1839 Gauss wrote to him: he read that memoir] mit grossem Interesse … doch bezog sich, … dieses Interesse weniger auf die Sache selbst, als auf Ihre Darstellung. Denn jene ist mir seit vielen Jahren familiär … Gauss (*Werke*, Bd. 8, pp. 146 – 147).

**156.** A twilight lasting all but a whole century [after Laplace and Poisson] fell over the European probability theory. … It might be stated that, in those times, in spite of winning ever more regions of applied knowledge, European probability not only did not develop further as a mathematical science, – it literally degraded. … ["The minor scientists"] usually inferred … that its theorems might be proved *not quite rigorously* … Khinchin (1937, p. 37/2005, p. 43).

**157.** Who can fail to see, with extreme regret, the absolute neglect, in the academic institutions, of one of the most important branches of mathematics? Only a few universities teach the rudiments of the theory of probability, and, up to now, there does not exist a single Russian work, or translation, on advanced, or even elementary probability. … We hope that Russian scientists will soon try to make up this deficiency. Brashman (1841). Translation: Gnedenko & Sheynin (1978/1992, p. 248).

**158.** [The central limit theorem that he, Chebyshev, obtained, was not derived] in a rigorous way. … we have made various assumptions but did not determine the boundary of the ensuing error. In its present state, mathematical analysis cannot derive this boundary in any satisfying fashion. Chebyshev (1879 – 1880/1936, p. 224; translation 2004, p. 206).

**Comment.** Chebyshev had in essence proved this theorem in 1887.

**159.** It was hardly possible to write down Chebyshev's lectures minutely, and it is natural that their extant record is fragmentary. Prudnikov (1964, p. 183).

**Comment.** I am inclined to believe this statement although it contradicts the utterance of A. N. Krylov, the Editor of the lectures (Chebyshev 1879 – 1889) published in 1936. Incidentally, this edition is corrupted by more than a hundred (I repeat: a hundred) mathematical misprints.

**160.** Especially towards the end of his life, as well as in his lectures, Chebyshev sometimes deviated from his own demands of clear formulas and rigour of proof in probability … Bernstein (1945/1964, p. 425; translation 2004, p. 84).

**161.** Tschebycheff's principal contribution to our subject is a simple formula for approximating, or rather finding an inferior limit, to a certain probability which is often required by the mathematical statistician. Edgeworth (1922, p. 109/1996, p. 155).

**Comment.** Edgeworth undoubtedly thought about the Bienaymé – Chebyshev inequality. His disregard of rigour in probability theory was characteristic of the time.

**162.** It [the Bienaymé – Chebyshev inequality] was first stated and proven by Chebyshev, but the main idea of [its] proof was pointed out much earlier by Bienaymé. Markov (1900/1913, p. 89; 1924, p. 92).

**Comment.** The French scientist provided the actual proof, but did not isolate the inequality from its context. True, conforming to the text of his memoir (1853) as a whole, he apparently only considered a sum of identically distributed terms.

**163.** I consider his [Chebyshev's] memoir as of *minor importance* since it contains that, which was sufficiently rigorously proved much earlier and included in generally known treatises [a reference follows: Laurent (1873, pp. 144 – 165)]. It is only interesting as being one of the successful applications of Chebyshev's great inventions to earlier exhausted problems. Nekrasov (1900, p. 384/2004, p. 33).

**Comment.** Everything is wrong here. Incidentally, Laurent had discussed that theorem on his pages 144 – 145 and I am inclined to believe that Nekrasov intentionally provided wrong page numbers to attribute to him a detailed investigation. The missing reference to Chebyshev should have been (1887; in Russian, 1890 – 1891, in French).

**164.** Nekrasov's work on the central limit theorem in general should be appraised as unsatisfactory. Having been a very powerful analyst, he chose an unfortunate, purely analytic rather than the stochastic approach which to a large extent predetermined his failure. Soloviev (1997, p. 21/2008, p. 363).

**165.** The principal upheaval due to Chebyshev consisted not only in his being the first to demand … that the limit theorems be proved with absolute rigour, … but mainly in that in each instance he strove to determine exact estimates of the deviations from limit regularities … in the form of inequalities unconditionally true for any number of trials.

Furthermore, Chebyshev was the first to appreciate clearly and use the full power of the concepts of *random variable* and its *expectation* (mean value). Kolmogorov (1947, p. 56); translation Gnedenko & Sheynin (1978/1992, pp. 254 – 255, where the last three lines were somehow mistakenly attributed to these two authors).

**Comment.** The expression *full power* seems unfortunate. Chebyshev had not introduced even a heuristic definition of (neither, therefore, a special notation for) a random variable. He was thus unable to study densities or generating functions as mathematical objects. In addition, the entire development of probability theory may be represented as an ever more complete use of the power of the concepts mentioned.

**166.** It is our specific Russian feature – an inclination to conservatism, to isolation from international science. Even Chebyshev …, in spite of his splendid analytical talent, was a pathological conservative. Novikov (2002, p. 330).

**Comment.** Cf. Liapunov's strange statement in No. 171.

**167.** One of the most distinguished [Russian] mathematicians … was in the habit of expressly telling his students that he did not advise [them] to engage in the philosophical aspect of mathematics since this was not very helpful for acquiring the knowledge of mathematics, and even rather harmful. V. A. Latyshev, 1893 as quoted by Prudnikov (1964, p. 91).

**Comment.** Prudnikov added that Latyshev had certainly meant Chebyshev. The further development of mathematics proved that this viewpoint was erroneous.

**168.** For its transformation from a mathematical amusement into a method of natural sciences the theory of probability is mainly obliged to the Petersburg school which accomplished this fundamental progress leaving West European mathematicians far behind. Bernstein (1940/1964, p. 11; translation 2004, p. 110).

The genius of Chebyshev and his associates, who, in this field [theory of probability], had left mathematicians of Western Europe far behind, have surmounted the crisis of the theory of probability that had brought to a stop its development a hundred years ago. Bernstein (1945/1964, p. 432; translation 2004, p. 92).

**Comment.** The first statement is utterly wrong. "Amusement" is an invention, pure and simple. Daniel Bernoulli studied smallpox epidemics, i. e. epidemiology, which also belongs to natural sciences, and this is just one example. And it was Markov who, in spite of his previous brilliant work, found himself "far behind" West European theory of probability. "Crisis" in the second statement means satisfaction with insufficient rigour.

**169.** Nunmehr finde ich, dass schon Liapunov … allgemeine Resultate dargelegt hat, die nicht nur über diejenigen des Herrn v. Mises hinausgehen, sondern aus denen auch die meisten der von mir bewiesen Tatsachen abgeleitet werden können. Lindeberg (1922, p. 211).

**170.** His [Chebyshev's] idea of mathematics had nothing in common with the science defined by Bertrand Russell, which knows not what it is speaking about, nor whether what it is saying is true (p. 65). The most outstanding exponent of Chebyshev's ideas [in the theory of probability], who shaped them with greatest skill in his classical course [1900] was A. A. Markov … [He made] up some deficiencies in Chebyshev's celebrated proof [of the central limit theorem] … (p. 72). Bernstein (1947/2001, pp. 65, 72).

**171.** The partisans of Riemann's extremely abstract ideas delve ever deeper into function-theoretic research and pseudo-geometric investigations in spaces of four, and of more than four dimensions*. In their investigations they sometimes go so far that any present or future possibility of seeing the meaning of their work with respect to some applications is lost, whereas Chebyshev and his followers invariably stick to reality.

   *These investigations were recently often connected, but have nothing in common with Lobachevsky's geometric research. Like Chebyshev, the great geometer always remained on solid ground and would have hardly seen these transcendental investigations as a development of his ideas. Liapunov (1895/1946, pp. 19 – 20).

**Comment.** The one-sidedness of the Chebyshev school is clearly revealed here. Liapunov possibly referred indirectly to Klein, who, in 1871, had presented a unified picture of the non-Euclidean geometry in which the findings of Lobachevsky and Riemann appeared as partial cases. Anyway, Liapunov rejected this unified approach.

**172.** In the beginning there was De Moivre, Laplace, and many Bernoullis, and they begat limit theorems, and the wise men saw that it was good and they called it by the name of Gauss. Then there were new generations and they said that it had experimental vigor but lacked in rigor. Then came Chebyshev, Liapunov and Markov and they begat a proof and Polya saw that it was momentous and he said that its name shall be called Central Limit Theorem.

   Then came Lindeberg and he said that it was elementary, for Taylor had expanded that which needed expansion and he said it twice, but Lévy had seen that Fourier transforms are characteristic functions and he said "Let them multiply and bring forth limit theorems and stable laws". And it was good, stable, and sufficient, but they asked "Is it necessary?" Lévy answered, "I shall say verily unto you that it is not necessary, but the time shall come when Gauss will have no parts except that they be in the image of Gauss himself, and then it will be necessary". It was a prophecy, and when Cramér announced that the time has come, and there was much rejoicing and Lévy said it must be recorded in the bibles and he did record it, and it came to pass that there were many limit theorems and many were central and they overflowed the chronicles and this was the history of the central limit theorem. Le Cam (1986, p. 86).

**173.** Bertrand and Poincaré wrote treatises on the calculus of probability, a subject neither of the two appeared to know. Le Cam (1986, p. 81).

**Comment.** This statement is ungenerous: in those times, no one except Markov *knew* probability theory, and even he failed in one important aspect (NNo. 190, 191).

**174.** [About the Poincaré treatise (1896).] The excessively respectful attitude towards … Bertrand is surprising. No traces of a special acquaintance with the literature on probability are seen. The course is written in such a way as though Laplace and Poisson, especially the latter, never lived. Bortkiewicz, in a letter to Chuprov of 1897 (Sheynin 1990b/1996, p. 40).

**175.** [Liapunov] understood and was able to appreciate the achievements of the West European mathematicians, made in the second half of the [19$^{th}$] century, better than the other representatives of the [Chebyshev] Petersburg school. Bernstein (1945/1964, p. 427; translation 2004, p. 87).

**176.** The work of Bernstein completes the classical direction of the Russian school of the theory of probability … and paves the way for the development of this important branch of mathematics. Kolmogorov & Sarmanov (1960, p. 219).

**177.** During the last 25 years, statistical science has made great progress, thanks to the brilliant schools of British and American statisticians, among whom the name of Professor R. A. Fisher should be mentioned in the foremost place. During the same time, largely owing to the work of French and Russian mathematicians, the classical calculus of probability has developed into a purely mathematical theory satisfying modern standards with respect to rigour. Cramér (1946, p. vii).

**178.** Limit theorems are usually rather decently formulated, but as a rule their proofs are helplessly long, difficult and entangled. Their sole raison d'être consists in obtaining comparatively simple stochastic distributions possibly describing some real phenomena. Tutubalin (1977, p. 59).

**179.** In probability theory, only very small (as compared, for example, with physics) groups of authors refer to each other. This means that the interest has narrowed which was largely caused by its unwieldy mathematical machinery and which is a typical sign of degeneration. Ibidem.

**180.** Je m'aperçus alors qu'en certain sens ce calcul [des probabilités] n'existait pas; il fallait le créer. Lévy (1970, p. 71).
**Comment.** He described the impressions of his youth.
**181.** It seemed clear that here [Lévy (1925)] was a first attempt to present the theory as a connected whole, using mathematically rigorous methods. It contained the first systematic exposition of the theory of random variables, their probability distributions and their characteristic functions. Cramér (1976, p. 516).

### 2.6. The Axiomatic Approach and the Frequentist Theory

**182.** Wer die Entwicklung der Wahrscheinlichkeitsrechnung in den letzten Jahrzehnten verfolgt, wird sich dem Eindruck nicht verschließen, dass dieser Zweig der mathematischen Wissenschaft in zweierlei Hinsicht hinter allen anderen bedeutend zurückgeblieben ist. Es fehlt einmal – nur wenige Arbeiten russischer Mathematiker bilden da eine Ausnahme – den *analytischen Sätzen* der Wahrscheinlichkeitsrechnung jene *Präzision der Formulierung und Beweisführung*, die in anderen Teilen der Analysis längst zur Selbstverständlichkeit geworden ist. Und es besteht andererseits, trotz mancher wertvoller Ansätze, über die *Grundlagen der Wahrscheinlichkeitsrechnung als einer mathematischen Disziplin* heute noch so gut wie keine Klarheit: was um so erstaunlicher erscheint, als wir nicht nur in einer Zeit lebhaften Interesses für axiomatische Fragen innerhalb der Mathematik, sondern auch in einer Epoche stetig wachsender Ausbreitung der Wahrscheinlichkeitsrechnung auf den verschiedensten Anwendungsgebiete leben. Mises (1919/1964, p. 35).
**183.** The claim to rank among the pure sciences must rest upon the degree in which it [the theory of probability] satisfies the following conditions: 1° That the principles upon which its methods are founded should be of an axiomatic nature. Boole (1854a/1952, p. 288).
**Comment.** Boole formulated two more general scientific conditions.
**184.** [In 1887 Tikhomandritsky showed Chebyshev his *Course* in probability. Chebyshev remarked that] it was necessary to transform the entire theory of probability. Tikhomandritsky (1898, p. iv).
**Comment.** Nothing more is known about this episode. In principle, apart from axiomatization, it was possible to consider random variables, densities, integral distributions, etc. as the main objects of probability theory. This, however, only happened some 30 years later.
**185.** [A] preliminary theorem, if it ought not rather be called an axiom … [Given,] mutually exclusive hypotheses … of which the probabilities … are $p_1, p_2, …$; if new information … changes the probabilities of some of them, … then the probabilities of [the other ones] have the *same ratios* to one another [as] before. The most important case … is [this. Some of the *n* hypotheses] must be rejected or … others must be admitted … Then [the other ones will] have the same mutual ratios [of their probabilities as] before. By means of these theorems, I think that we may … avoid the necessity of constructions by means of balls, dice, etc. Donkin (1851, p. 356).
**Comment.** Boole (1854b/2003, p. 163) thought that "perhaps [Donkin's principle] might … be regarded as axiomatic".
**186.** I shall not defend these basic theorems connected to the basic notions of the calculus of probability, notions of equal probability, of independence of events, and so on, since I know that one can argue endlessly on the basic principles even of a precise science such as geometry. Markov (1911/1981, p. 149).
**Comment.** It is clearly seen here and in No. 109 that probability theory could not have been rigorously constructed until the advent of the axiomatic approach.
**187.** Various concepts are defined not by words, each of which can in turn demand definition, but rather by [our] attitude towards them ascertained little by little. Markov (1900; 1908, p. 2; 1924, p. 2).
**Comment.** Also see his unfortunate attempt at reforming the theory of probability (NNo. 190, 191).
**188.** I shall not go a step out of that region where my competence is beyond any doubt. Markov, letter of 1910 to Chuprov; Ondar (1977/1981, p. 59).
**Comment.** This principle seems too rigid; anyway, it explains why Markov never said a word about the possibility of applying his *chains* in natural sciences. I have included the next item to show one of Markov's strong mathematical traits.
**189.** Many mathematicians apparently believe that going beyond the field of abstract reasoning into the sphere of effective calculations would be humiliating. Markov (1899a, p. 30).
**Comment.** The best example of his own calculations is his table of the normal distribution (1888). With regard to precision, it remained beyond compare up to the 1930s.
**190.** Axiom. If, given some information, events *p, q, r, …, u, v* are equally possible and divided, with respect to event *A*, into favourable and unfavourable, then, upon adding to this information the indication that *A* has occurred, those from among the events *p, q, r, …, u, v* which are not favourable to event *A* fall through, whereas the others remain equally possible as they were previously. Markov (1900/1924, p. 10).
**191.** The addition and the multiplication theorems along with the axiom mentioned above serve as an unshakeable base for the calculus of probability as a chapter of pure mathematics. Ibidem, p. 24.
**Comment.** Markov proved these theorems with a reference to his axiom, which, however, was not of a mathematical nature and his statement seems absolutely wrong. The axiom persisted through all the editions of Markov's treatise, and his final conclusion first appeared in 1913. No one however had ever mentioned either it, or the axiom, or the *unshakeable base*. On p. 159 Markov independently formulated what could now be called the extended axiom of addition (and multiplication) but his final declaration was certainly patently wrong. Note that his axiom was a particular case of Donkin's theorem (No. 185) and that Donkin only weakly hinted at an axiomatic theory of probability. Being Chebyshev's follower, Markov underrated the axiomatic direction of probability (and the theory of functions of complex variables), see A. A. Youshkevich (1974, p. 125).

**192.** Durch die Untersuchungen über die Grundlagen der Geometrie wird uns die Aufgabe nahe gelegt, nach diesem Vorbilde diejenigen physikalischen Disziplinen axiomatisch zu behandeln, in denen schon heute die Mathematik eine hervorragende Rolle spielt: dies sind in erster Linie die Wahrscheinlichkeitsrechnung und die Mechanik.

Was die Axiome der Wahrscheinlichkeitsrechnung angeht, so scheint es mir wünschenwert, dass mit der logischen Untersuchung derselben zugleich eine strenge und befriedigende Entwicklung der Methode der mittleren Werte in der mathematischen Physik, speziell in der kinetischen Gastheorie Hand in Hand gehe. Hilbert (1901/1970, p. 306).

**Comment.** In those days, probability theory had been an applied mathematical, but hardly physical science. Perhaps from the beginning of the 19th century onwards *mathematical physics* was the same as theoretical part of physics, see the title of Biot's fundamental treatise (1816). Then, in 1850, in accordance with Poisson's long-standing proposal, a chair of the calculus of probability and mathematical physics was established at the École Polytechnique (Bru 1981, p. 87) and even 62 years later the subtitle of Poincaré (1912) read *Cours de physique mathématique*. The term *theory of means* had been widely used beginning at least with Condorcet instead of the already existing expression *theory of errors* (No. 657). Moreover, the former theory was more general since it also studied such means that did not correspond to any object in the real world (e. g., mean price of bread), see Sheynin (1986, p. 311), where, in particular, I quote Quetelet (1846, pp. 60, 65 – 67).

**193.** Zweck des vorliegenden Heftes ist eine axiomatische Begründung der Wahrscheinlichkeitsrechnung. Der leitende Gedanke des Verfassers war dabei, die Grundbegriffe der Wahrscheinlichkeitsrechnung, welche noch unlängst für ganz eigenartig galten, natürlicherweise in die Reihe der allgemeinen Begriffsbildungen der modernen Mathematik einzuordnen (p. III).

Herrn A. Khintchine, der das ganze Manuskript sorgfältig durchgelesen und dabei mehrere Verbesserungen vorgeschlagen hat, danke ich an dieser Stelle herzlich (p. IV). Kolmogorov (1933, pp. III, IV).

**Comment.** I adduced Kolmogorov's Acknowledgement because no one had mentioned Khinchin in this connection.

**194.** Die Wahrscheinlichkeitstheorie als mathematische Disziplin soll und kann genau in demselben Sinne axiomatisiert werden wie die Geometrie oder die Algebra. Das bedeutet, dass, nachdem die Namen der zu untersuchenden Gegenstände und ihrer Grundbeziehungen sowie die Axiome, denen diese Grundbeziehungen zu gehorchen haben, angegeben sind, die ganze weitere Darstellung sich ausschließlich auf diese Axiome gründen soll und keine Rücksicht auf die jeweilige konkrete Bedeutung dieser Gegenstände und Beziehungen nehmen darf. Ibidem, p. 1.
**Comment.** I am adducing without comment a recent statement (Vovk & Shafer 2003, p. 27): [In our book Shafer & Vovk (2001)] *we show how the classical core of probability theory can be based directly on game-theoretic martingales, with no appeal to measure theory. Probability again becomes* [a] *secondary concept but is now defined in terms of martingales …*
**195.** Kolmogoroff a donné à l'axiomatique du calcul des probabilités une forme qui semble définitive. Lévy (1949, p. 55).
**196.** Das Kolmogorows Buch der entscheidende Durchbruch zur Entwicklung der Wahrscheinlichkeitstheorie als mathematische Disziplin war, lag wohl weniger daran, dass nun die Axiome explizit und klar ausgesprochen waren, sondern daran, dass es Kolmogorow gelang die theoretischen Grundlagen für die Theorie der stochastischen Prozesse zu legen, bedingte Wahrscheinlichkeiten rigoros zu behandeln, und einen einheitlichen Rahmen für viele klassische Fragen zu schaffen. Krengel (1990, p. 461).
**197.** To most mathematicians mathematical probability was to mathematics as black marketing to marketing; … The confusion between probability and the phenomena to which it is applied … still plagues the subject; [the significance of the Kolmogorov monograph] was not appreciated for years, and some mathematicians sneered that … perhaps probability needed rigor, but surely not *rigor mortis*; … The role of measure theory in probability … still embarrasses some who like to think that mathematical probability is not a part of analysis. Doob (1989).
**198.** It was some time before Kolmogorov's basis was accepted by probabilists. The idea that a (mathematical) random variable is simply a function, with no romantic connotation, seemed rather humiliating to some probabilists. Doob (1994, p. 593).
**199.** Kolmogorov's book still ranks as the basic document of modern probability theory. If in 1920 it might be said … that this theory was not a mathematical subject, it was impossible to express such an opinion after the publication of this book in 1933. Cramér (1976, p. 520).
**200.** 1933 kam A. N. Kolmogorov mit dem Ei des Kolumbus. Freudenthal & Steiner (1966, p. 190).
**Comment.** As the legend goes, Columbus cracked an egg which enabled it to stand firmly on his table.
**201.** Die Aufgabe, die Wahrscheinlichkeitsrechnung im Geiste der modernen Mathematik exakt axiomatisch zu begründen, wurde in befriedigender Weise erstmalig von A. N. Kolmogorov … gelöst … Rényi (1969/1969, p. 84).
**202.** Kolmogorov's axiomatics is the one widely accepted. However, the concept itself of practical applications largely follows the Mises idea. Tutubalin (1977, p. 15).
**Comment.** Cf. NNo. 58 and 59.
**203.** Das Schibboleth [the touchstone] in wahrscheinlichkeits-theoretischer Axiomatik ist der Produktsatz für unabhängige stochastische Größen (oder Ereignisse). Wer naiv eine halb mathematische, halb physikalische Wahrscheinlichkeitstheorie betreibt, sagt "wenn zwei Ereignisse unabhängig sind, so gilt für ihre Wahrscheinlichkeiten der Produktsatz", der Axiomatiker dagegen "wenn für sie der Produktsatz gilt, so heißen sie unabhängig". Freudenthal & Steiner (1966, p. 189).

**204.** Nach meiner Ansicht ist die Misessche Grundlegung der Wahrscheinlichkeitsrechnung total verunglückt und wir sind durch sie nicht besser gestellt, als wenn wir einfach direkt von einer Verteilungsfunktion und den zugehörigen Stieltjes-Integralen ausgehen. Hausdorff, letter of 1920; Girlich (1996, p. 42).

**205.** The readiness with which we physicists and many others who use the theory in practice adopt Frequency as the basis, does mean taking things a little too easy. There are grave objections to this, mainly that we thereby cut ourselves off from ever applying rational probability considerations to a single event. I attempt here to keep aloof from the frequency definition of probability and yet to reduce the volume of *axiomatic* statements to a minimum. Schrödinger (1945 – 1948, p. 51).

**206.** The theory of probability is a doctrine of mass phenomena (p. 400). It can be only linked to set theory … as the doctrine about totalities of the most general kind (p. 402). It is impossible on principle (as Mises admits) to construct an individual collective (p. 408). The theory of probability can be developed on the basis of the axiomatic theory in an incomparably simpler and easier way than when it is founded on the frequentist theory (p. 412). The foundation of the theory of probability offered by the frequentist theory is not therefore a logical basis in the sense of modern mathematics (p. 418). Khinchin (1961/2004, pp. 400, 402, 408, 412, 418).

**207.** For Mises, there were never two theories, one "pure", the other "applied", but one theory only, a frequency theory, mathematically rigorous and guided by an operational approach. Hilda Geiringer; Mises (1964, p. v).

**208.** Pourtant de nombreux savants ont de tout temps défendu la théorie empirique. R. von Mises a d'une manière très ingénieuse essayé de remédier à ses difficultés. Mais c'est aussi impossible que la quadrature du cercle. Lévy (1970, p. 79).

**209.** The main feature of the Mises approach consists in that, from the very beginning, it assumes everything as it actually happens in an experiment. … Its opposition to the axiomatic method is entirely based on an misunderstanding. Alimov (1980, pp. 32 – 33).

**210.** [Fechner's] Ausführungen bildeten – wenigsten für mich – die Anregung zu der neuen Betrachtungsweise. Mises (1928/1972, p. 99).

**Comment.** In the first edition (1928) of that source, Mises, when introducing his *collective*, remarked that that notion adjoined Fechner's *Kollektivgegenstand.*

**211.** Unter einem Kollektivgegenstande … verstehe ich einen Gegenstand, der aus unbestimmt vielen, nach Zufall variierenden, Exemplaren besteht, die durch einen Art- oder Gattungsbegriff zusammenhalten werden. Fechner (1897, p. 3).

**Comment.** In statistics, Fechner is indeed meritorious, mostly for his attempts to study systematically observational series with pertinent asymmetric densities.

**212.** All the leading statisticians from Poisson to Quetelet, Galton, Edgeworth and Fechner … have realized that asymmetry must be in some way described before we can advance in our theory of variation [in biology]. K. Pearson (1905, p. 189).

**Comment.** Pearson had not documented his statement and I am unable to confirm it with respect to Poisson. As to Quetelet, he wavered, to say the least, cf. No. 486C.

**213.** Ich war immer für die Ideen G. T. Fechners zugänglich und habe mich auch in wichtigen Punkten an diesen Denker angelehnt. Freud (1925/1963, p. 86).

**214.** Der babylonische Thurm wurde nicht vollendet, weil die Werkleute sich nicht verständigen konnten, wie sie ihn bauen sollten; mein psychophysisches Bauwerk dürfte bestehen bleiben, weil die Werkleute sich nicht werden verständigen können, wie sie es einreißen sollen. Fechner (1877, p. 215).

**Comment.** Fechner's investigations (directly bearing on statistics) had indeed been criticized; see however NNo. 210 – 212.

## 3. Statistics and Mathematical Statistics

### 3.1. Origin and Aims

**215.** Statistics has long had a neighbourly relation with philosophy of science in the epistemological city although statistics has usually been more modest in scope and more pragmatic in outlook. In a strict sense, statistics is part of philosophy of science, but in fact the two areas are usually studied separately. Kruskal (1978, p. 1082).

**216.** Statistics, as we now understand the term, did not commence until the seventeenth century, and then not in the field of 'statistics' [Staatswissenschaft] (§ 4). The true ancestor of modern statistics is … Political Arithmetic (§ 5). Kendall (1960, §§ 4, 5).

**217.** [According to unnamed German authors, 1806 – 1807] die wahren Staatskräfte seien geistig, nicht materiell. Die Statistik dürfte nicht zum Skelette herabgewürdigt, sondern ihr eine höhere, eine lebendigere Tendenz untergelegt werden … [The issues should be] der Nationalgeist, die Freiheitsliebe, das Genie und der Charakter der großen und der kleinen Männer an der Spitze des Staates … Knies (1850, p. 24).

**218.** Une description simple va au fond des choses, tandis que les nombres et les tableaux s'arrêtent à la surface. Mone (1824/1834) as quoted by Knies (1850, p. 80).

**Comment.** Other opinions, both similar and opposite to Mone's statement, can be quoted. Roslavsky (1839/1963, pp. 181 – 182): the "number element" was an innovation [?] to be admitted "with greatest caution". Moreau de Jonnès (1847) began his book by stating that statistics was "la science des faits sociaux exprimés par de termes numériques", but his subject index lacked such expressions as mathematics or probability theory. And Dufau (1840) inserted the words "l'étude des lois …" in the subtitle of his book.

**219.** What is a common measure of Time, Space, Weight & motion? What number of Elementall sounds or letters will … make a speech or language? How to give names to names, and how to adde and subtract sensate & to ballance the weight and power of words; which is Logick & reason? Petty (1927, vol. 2, pp. 39 – 40).

**Comment.** The cofounder of political arithmetic was a philosopher of science, congenial in some respects with Leibniz, his younger contemporary.

**220.** Petty's claim to fame as an economist lies not so much in his originality or his theoretical ability as in his analytical skill. His insistence on measurement and his clear schematic view of the economy make him the first econometrician, and he was certainly evolving and using concepts and analytical methods that were in advance of his time. Deane (1978, p. 703).

**221.** [Petty rejected the use of] comparative and superlative Words [and decided to express himself] in terms of Number, Weight, or Measure; [of using] only Arguments of Sense [and of considering] only such Causes as have visible Foundations in Nature. Petty (1690/1899, p. 244).

**Comment.** Petty was cofounder of political arithmetic, the predecessor of statistics.

**222.** Im politischen Haushalt, wie bei Erforschung von Naturerscheinigungen sind die Zahlen immer das Entscheidende; sie sind die letzten unerbittlichen Richter … Humboldt, in 1838, as quoted by Knies (1850, p. 145).

**223.** I often say that when you can measure what you are speaking about, and express it in numbers, you know something about it; but when you cannot measure it, when you cannot express it in numbers, your knowledge is of a meagre and unsatisfactory kind … Thomson, Lord Kelvin (1883/1889, p. 73).

**224.** Measurement does not necessarily mean progress. Failing the possibility of measuring that which you desire, the lust for measurement may, for example, merely result in your measuring something else – and perhaps forgetting the difference – or in your ignoring some things because they cannot be measured. G. U. Yule in his correspondence (with Kendall?) as quoted by Kendall (1952/1970, p. 423).

**Comment.** And here is a similar statement (Andreski 1972, p. 120): *The gravest kind of danger stems from the illusion that, because certain kinds of data can be quantified and processed by a computer, therefore they must be more important than those which cannot be measured.*

**225.** I found that all mentioned to dye of the *French Pox* were returned by the *Clerks* of Saint *Giles*'s and Saint *Martin*'s *in the Fields* only, in which place I understood that most of the vilest and most miserable Houses of Uncleanness were: from whence I concluded, that only *hated* persons, and such, whose very *Noses* were eaten off, were reported by the *Searchers* to have died of this too frequent *Malady*. Graunt (1662/1899, p. 356).

**Comment.** This is a classical example of revealing systematic corruption in the data.

**226.** In this Parish there were born 15 *Females* for 16 *Males*, whereas in *London* there were 13 for 14, which shews, that *London* is somewhat more apt to produce *Males* than the Country. Ibidem, p. 389.

**Comment.** Graunt rather than Arbuthnot (1712) was the first to note that more boys were being born than girls. His comparison of a tiny parish with London is not however convincing.

**227.** Now having (I know not by what accident) engaged my thoughts upon the *Bills of Mortality*, and so far succeeded therein, as to have reduced several great confused *Volumes* into a few perspicuous *Tables*, and abridged such *Observations* as naturally flowed from them, into a few succinct *Paragraphs*, without any long Series of multiloquacious *Deductions* … Ibidem, p. 320.

**Comment.** In actual fact, Graunt had thus formulated the aims of statistics, cf. NNo. 257, 392.

**228.** [After describing economic issues] it is no less necessary to know how many People there be of each Sex, State, Age, Religion, Trade, Rank, or Degree, etc. by the knowledge whereof, Trade and Government may be made more certain and Regular. Ibidem, p. 396.

**229.** I conclude, That a clear knowledge of all these particulars and many more, whereat I have shot but at rovers, is necessary, in order to good, certain and easie Government and even to balance Parties and Factions both in *Church* and *State*. But whether the knowledge thereof be necessary to many, or fit for others than the Sovereign and his chief Ministers, I leave to consideration. Ibidem, p. 397.
**230.** [Graunt is] the father of modern epidemiology. Greenwood (1932, p. 10).
**231.** The concept of a life table was an outstanding innovation and it lay ready for Halley's use (p. 13). … Graunt's work created the subject of demography. [It] contributed to statistics in general (p. 14). D. V. Glass (1963, pp. 13, 14).
**232.** Graunt is memorable mainly because he discovered … the uniformity and predictability of many biological phenomena taken in the mass … Thus he, more than any more man, was the founder of statistics. Willcox (Graunt 1662/1939, p. XIII).
**233.** Petty perhaps suggested the subject of the inquiry … probably assisted with comments upon medical and other questions here and there … proposed [some] figures … and may have revised or even written the Conclusion … W. Hull (Petty 1899, vol. 1, p. LII).
**Comment.** Hull left too little for Graunt.
**234.** I have also (like the author of those *Observations* [like Graunt]) Dedicated this *Discourse* to … the Duke of Newcastle … Petty (1674).
**235.** Es ist die Statistik ohnehin keine Disciplin, die man mit einem leeren Kopf gleich fassen kann. Es gehöret eine wohldigerirte Philosophie, eine geläufige Kenntnis der Europäischen Staats- wie auch der Natur-Geschichte nebst einer Menge von Begriffen und Grundsätzen dazu, und die vielen ganz verschiedenen Artikel der Staatsverfassung der heutigen Reiche gehörig begreifen zu können. Achenwall (1752/1756, Vorrede).
**236.** Geschichte ist eine fortlaufende Statistik, und Statistik stillstehende Geschichte. Schlözer (1804, p. 86).
**Comment.** The author provided this pithy saying as an illustration, but his followers (of the *Staatswissenschaft* direction) adopted it as a definition of statistics. A *stillstehende Geschichte* did not demand any studies of causes and effects.
**237.** Man unterschied zwischen höherer [Achenwallian] und gemeinen Statistik … bis auf unsere Tag. Knies (1850, p. 25).
**238.** Moses sent them to spy out the land of Canaan, and said to them, "Go up into the Negeb yonder, and go up into the hill country, and see what the land is, and whether the people who dwell in it are strong or weak, whether they are few or many, and whether the land that they dwell in is good or bad, and whether the cities that they dwell in are camps or strongholds, and whether the land is rich or poor, and whether there is wood in it or not. Numbers 13:17.
**Comment.** Graunt's endeavours (1662) were in tune with Moses' aims.
**239.** [Statistics] diffère beaucoup de la science de l'économie politique, qui examine et compare les effets des institutions, et recherche les causes principales de la richesse et de la prospérité des peuples. Ces considérations … ne sont point le premier objet de la statistique qui exclut presque toujours les discussions et les conjectures. L'arithmétique politique … doit aussi être distinguée de la statistique. Delambre (1819, p. LXVIII).
**Comment.** A strange statement characteristic of those times.
**240.** L'esprit de dissertation et de conjectures est, en général, oppose aux véritables progrès de la statistique, que est surtout une science d'observation. Fourier (1821, pp. iv – v).
**Comment.** Fourier was Editor of a collection of statistical tables in four volumes (1821 – 1829) describing Paris and the Département de la Seine. Nevertheless, bearing in mind his general standing, his viewpoint seems strange.
**241.** [The statistician is] not obliged to treat causes or consequences. Poroshin (1838/1963, p. 101).
**242.** [Statistics] does not discuss causes, nor reason upon probable effects; it seeks only to collect, arrange and compare, that class of facts which alone can form the basis of correct conclusions … with respect to social and political government. …
All conclusions shall be drawn from well-attested data and shall admit of mathematical demonstration (p. 1). Statistics … are closely allied to the other sciences, and receive contributions from all of them … (p. 2). Like other sciences, that of Statistics seeks to deduce from well-established facts certain general principles … does not admit of any kind of speculation … It is not, however, true that Statistics consist merely of columns of figures; it is simply required that all conclusions shall be drawn from well-attested data, and shall admit of mathematical demonstration (p. 3). Anonymous (1839, pp. 1, 2, 3).
**Comment.** From the editorial in the first issue of the official periodical of the just established London (later Royal) Statistical Society.
**243.** These absurd restrictions [about investigating causes and effects] have been necessarily disregarded in … numerous papers. Woolhouse (1873, p. 39).
**244.** Such statisticians who observe without thinking about the *why* or the *how*, who make most involved computations without understanding where all their multiplications and divisions might and will lead them, are extremely numerous. And statistics has to thank them for its ill fame … Chuprov (1903/1960, p. 42).
**245.** [Just as in history it is neccesary] nicht nur das Pourquoi, selbst auch das Pourqoi von dem Pourquoi zu erforschen, so wird es auch bey der Statistik … den gegenwärtigen Zustand eines Staat aus dem vergangenen Zuständen begreiflich zu machen. Gatterer (1775, p. 15).
**Comment.** The author apparently borrowed the beginning of his phrase from Leibniz (Sheynin 2001a, p. 90).
**246.** [Statistics can establish] règles de conduite pour l'avenir; [should] estimer le degré de prospérité de … population, sa force, ses besoins, et jusqu'à un certain point se faire des idées justes sur son avenir. Quetelet (1829, p. 9).

**247.** La statistique a donc pour objet de nous présenter l'exposé fidèle d'un état, à une époque déterminée (p. 264). … Je pense que la définition que je propose, et qui du reste s'écarte peu de celle donnée par plusieurs savants modernes, circonscrit suffisamment les attributions de la statistique, pour qu'on ne puisse pas la confondre avec les sciences historiques ou les autres sciences politiques et morales qui s'en rapprochent le plus. Elle ne s'occupe d'un état que pour une époque déterminée; elle ne réunit que les éléments qui se rattachent à la vie de cet état, s'applique à les rendre comparables et les combine de la manière la plus avantageuse pour reconnaître tous les faits qu'ils peuvent nous révéler (pp. 268 – 269). Quetelet (1846, pp. 264, 268 – 269).

**248.** En l'absence de documents exacts, recueillis par la science, il doit se borner bien souvent à prendre des données voisines de la vérité. Ibidem, p. 73.

**249.** Il convient de ne jamais oublier qu'un document statistique n'est pas un document certain, mais probable seulement; c'est dans l'estimation de cette probabilité que consiste l'importance du résultat que l'on considère, et que réside en général toute l'utilité des calculs statistiques. Quetelet & Heuschling (1865, p. LXV).

**Comment.** Here is Quetelet's pertinent remark (1869, t. 1, p. 112):

*Peut-on appliquer des corrections mathématiques à des nombres, quand on est persuadé que ces corrections sont dépassées de beaucoup par les erreurs qu'on néglige?*

But what exactly had he attributed to mathematical corrections? In essence, Quetelet's opinion was apparently caused by the existence of systematic corruptions, incompleteness of data and forgeries.

**250.** La statistique a la mission d'apprécier la valeur des documents qu'elle rassemble et d'en déduire des conclusions. Quetelet (1869, t. 1, pp. 102 – 103).

**Comment.** This is too narrow, cf. NNo. 257, 392.

**251.** L'heureuse influence qu'exercent nécessairement sur les individus et sur la société des doctrines vraies, de bonnes lois, de sages institutions, ne se trouve pas seulement démontrée par le raisonnement et par la logique, elle se démontre aussi par l'expérience. Par conséquent, la Statistique offre un moyen en quelque sorte infaillible de juger si une doctrine est vraie ou fausse, saine au dépravée, si une institution est utile ou nuisible aux intérêts d'un peuple et à son bonheur. Il est peut-être à regretter que ce moyen ne soit pas plus souvent mis en œuvre avec toute la rigueur qu'exige la solution des problèmes; il suffirait à jeter une grande lumière sur des vérités obscurcies par les passions; il suffirait à détruire bien des erreurs. Cauchy (1845/1896, p. 242).

**Comment.** Cf., however, No. 453.

**252.** The time may not be very remote when it will be understood that for complete initiation as an efficient citizen of one of the new great complex worldwide states that are now developing, it is as necessary to be able to compute, to think in averages and maxima and minima, as it is now to be able to read and write. Wells (1904, pp. 191 – 192) as quoted by Moritz (1964, No. 1201).
**Comment.** Wells (1931, vol. 1, pp. 432 – 433), as quoted by Tee (1980) had more to say, for example: "The movement of the last hundred years is all in favour of the statistician".
**253.** It is to the statistician that the present age turns for what is more essential in all its more important activities. Fisher (1953, p. 2).
**Comment.** The scope of modern statistics can be revealed by studying the *Handbook of Social Indicators* (1989, published by the United Nations) which lists several hundred indicators separated under 13 heads; see also De Vries (2001). Several pertinent papers are included in the *International Statistical Review*, vol. 71, No. 1, 2003. However, it also necessary to take into account the requirements of natural sciences.
**254.** Lies, damned lies and statistics.
**Comment.** More than a hundred years ago, Mark Twain indicated that this pithy saying had been attributed to Disraeli, Lord Beaconsfield. However, in 2005 P. M. Lee stated, in the *Newsletter* of the Royal Statistical Society, that the real author was Lord L. H. Courtney and referred to Baines (1896). There, on p. 87, both Lord C. and the saying itself (treated as generally known) were indeed mentioned, but no definite source was cited.
**255.** No study is less alluring or more dry and tedious than statistics, unless the mind and imagination are set to or that the person studying is particularly interested in the subject; which last can seldom be the case with young men in any rank of life. Playfair (1801, p. 16) as quoted by Gaither et al (1996).
**256.** For the community as a whole, there is nothing so extravagantly expensive as ignorance. Shaw (1926/1942, p. v).
**257.** [Statistics is] die Wissenschaft der Kunst statistische Data zu erkennen und zu würdigen, solche zu sammeln und zu ordnen. Butte (1808, p. xi).
**258.** Le but principal de la statistique [is to render different materials comparable]. Quetelet (1845, p. 225).
**259.** Quand il s'agit de deux pays différents, il semble qu'on ait pris plaisir à rendre toute espèce de rapprochement impossible. Quetelet (1846, p. 364).
**Comment.** Quetelet is meritorious for his attempts to standardize statistics on an international level and Karl Pearson (1914 – 1930, vol. 2, p. 420) highly praised his achievements *in organizing official statistics in Belgium and unifying international statistics*.
**260.** Les unes voudraient tout réduire à des nombres et faire consister la science dans un vaste recueil de tableaux; d'autres, au contraire, semblent craindre les nombres et ne les regardent que comme donnant des idées incomplètes et superficielles des choses. Ces deux excès seraient également nuisibles … Quetelet (1846, p. 432.
**Comment.** *Les unes* evidently were the partisans of the Staatswissenschaft (see NNo. 235 and 236) and *d'autres*, those keeping to the so-called numerical method (No. 310).

**261.** Pour que la statistique mérite le nom de science, elle ne doit pas consister simplement dans une compilation de faits et de chiffres: elle doit avoir sa théorie, ses règles, ses principes. Cournot (1843, §105.

**262.** Allgemein anerkannte Prinzipien, an denen die Richtigkeit der Schlüsse und die Zweckmäßigkeit der angewendeten Methoden geprüft werden könnten, gibt es gegenwärtig in der Statistik gar keine. Chuprov (1905, p. 422).

**263.** Le but essentiel du statisticien, comme de tout autre observateur, est de pénétrer autant que possible dans la connaissance de la chose en soi. Cournot (1843, § 106; roughly repeated in § 120).

**264.** Years ago a statistician might have claimed that statistics deals with the processing of data. … Today's statistician will be more likely to say that statistics is concerned with decision making in the face of uncertainty. Chernoff & Moses (1959, p. 1).

**265.** More and more, the word *statistics* signifies *automatic treatment of data* … Ekeland (1991/1993, p. 167).

**266.** [Bismarck] für die Statistik sehr wenig übrig hatte und sie eigentlich für entbehrlich hielt. Saenger (1935, p. 452).

**Comment.** Saenger had not justified his statement which deserves to be investigated.

**267.** They [the poets] will sometime array the coarse statistical prose in verses because numbers reveal power, rule, human weakness, the course of history and many other … aspects of the world. Mendeleev (1888/1949, p. 54).

**268.** Isolating [from statistical methodology], first, that which relates to the properties of judgements and concepts, i. e., to logic, and then to the properties of quantitative images upon which it [?] is operating, i. e., of mathematics, we nevertheless obtain some remainder for which no acknowledged sanctuary is in existence, which remains uncoordinated and homeless until we perceive its special theoretical essence and provide it with the missing unity in the system of judgements fully deserving the name of theoretical statistics. Slutsky (1916, p. 110/2009, p. 94).

**Comment.** I would replace *logic* by *philosophy*. This definition seems to be too complicated.

**269.** What is the present outlook for statistics? Some hints of future possibilities are suggested by the foregoing outline of its development. … I have already mentioned the immense progress which has been made of late. It might perhaps seem difficult to make further advances, but the movement has been so rapid, that we can hardly think it will at once come to a standstill, and … some opportunities … are evident. We need only to follow the lead of old political arithmeticians. An essential feature in those early days was the frequent use of representative statistics … This method was frequently used too naively, but its principle was not bad. Laplace pointed in the right direction … What remains to be done, is simply to develop a theory of these representative enumerations … In fact in many cases it will be practically impossible to do without representative statistics. Westergaard (1916, pp. 237 – 238).

**Comment.** Chuprov (1912) devoted a popular report to sampling.

**270.** Die Zukunft der "materiellen Statistik" in einer scharfen Hervorkehrung des idiographischen Momentes, – in einer planmässigen Rückkehr zur alten beschreibenden Staatenkunde [Staatswissenschaft], – liegt. Das ehrwürdige Geschöpf der deutschen "Universitätsstatistik" darf natürlich nicht aus dem Grabe, in welchem es ein Jahrhundert lang geschlummert hat, so hervorgeholt werden, wie es zur Ruhe gebracht wurde. … Die Statistik, als selbstständige Wissenschaft, wird die Form … einnehmen, wobei zahlenmässig characterisierte Massenerscheinungen … in den Vordergrund treten, aber doch nicht zur Alleinherrschaft gelangen werden … Chuprov (1922b, p. 339).
**Comment.** The Staatswissenschaft did not die, at least not in Germany, and does not at all abstain from figures. Ideography is a philosophical term, completely forgotten in statistics, signifying description rather than study. I think that Chuprov was mistaken.
**271.** Les mesures géodésiques, les observations relatives aux températures et à l'état de l'atmosphère, aux maladies communes, à la salubrité de l'air, des alimens et des eaux, l'exposition des procédés des arts, les descriptions minéralogiques appartiennent sans doute à la statistique …mais cette science n'a point pour but de perfectionner les théories … Delambre (1819, p. LXX).
**Comment.** Strange as it seems, this evidently means that before Quetelet quite alien from statistics subjects, as we would now say, had been ascribed to it. As to the results of geodetic measurements, we would rather consider them data for the theory of errors.
   Is this partly what Chuprov (No. 270) had in view?

### 3.2. Population and Moral Statistics, Insurance, Games of Chance

**272.** Instead of murmuring at what we call an untimely Death, we ought with Patience and unconcern to submit to that Dissolution which is the necessary Condition of our perishable Materials, and of our nice and frail Structure and Composition. Halley (1694, p. 655).
**273.** Ce sont donc deux choses différentes que l'espérance ou la valeur de l'aage futur d'une personne et l'aage auquel il y a égale apparence qu'il parviendra ou ne parviendra pas. Le premier est pour régler les rentes a vie, et l'autre pour les gageures. Huygens, letter of 1669; Huygens (1888 – 1950, t. 6, p. 537).
**Comment.** Huygens had not explained expectation quite clearly.
**274.** A *Census* … once established, and repeated at proper intervals, would furnish to our Governours, and to ourselves, much important instruction of which we are now in a great measure destitute: Especially if the whole was distributed into the proper *Classes* of *married* and *unmarried*, *industrious* and *chargeable* Poor, *Artificers* of every kind, *Manufacturers*, &c. and if this was done in each County, City and Borough, separately; that particular useful conclusions might thence be readily deduced; as well as the general state of the Nation discovered; and the Rate according to which *human Life* is wasting from year to year (p. 348). The years that are multiples of 10 are generally overloaded (p. 347). De Moivre (1724/1756, pp. 348, 347).

**Comment.** Unlike Graunt (No. 229), De Moivre already believes that everyone requires statistics.

**275.** Das in Geburt, Vermehrung, Fortpflanzung, im Leben, Tode und in den Ursachen des Todes eine beständige, allgemeine, grosse, vollkommene und schöne Ordnung herrschte: dieses ist die Sache, die in dieser Abhandlung soll erwiesen werden (§ 12).

Es hat … mein hochgeschätzer Freund und College, der Herr Professor Euler, … ausser dem mir bey der Berechnung der Verdoppelung geleisteten Beystande, die Durchlesung der abgedruckten Bogen gütigst übernommen, und es kann mich desselben bezeugte Zufriedenheit und freundschaftliches, jedoch unpartheyisches Urtheil in etwas beruhigen … (Vorrede).

Der weiseste Schöpfer und Regierer der Welt lässet das zahlreiche Heer des menschlichen Geschlechts durch die Zeugung aus seinem Nichts hervorgehen, so viel er derselben zum Leben geordnet hat. Der Ewige lässet uns in der Zeit gleichsam vor seinem Angesichte vorbey gehen, bis wir, nach Erreichung des einem jeden gesteckten Zieles, wiederum von diesem Schauplatze abtreten. …

Wen nun das zahlreiche Heer der Menschen gleichsam im Marsche betrachtet wird; so kann man es sich als in verschiedene Züge abgetheilt vorstellen. Die in jedem Jahre des menschlichen Alters Lebende machen einen Zug oder eine Abtheilung aus. Man kann sie sich auch nach grössern Zügen, deren jeder 5 oder sogar 10 Jahre in sich begreift, vorstellen. Hier ist nun zwar kein Zug so groß, wie der andere; aber es hat doch allezeit ein jeder seine richtige Proportion gegen das ganze Heer, und wird dadurch bestimmet. Süssmilch (1761 – 1762, Tl. 1, § 12, Vorrede, § 14).

**Comment.** As an ordained chaplain, Süssmilch participated in military operations which explains his picturing mankind as a marching regiment. Euler participated in preparing the second edition (1761 – 1762) of Süssmilch's contribution, was coauthor of at least one of its chapters (Chapter 8). In particular, he compiled a table showing that increase under quite arbitrary and simplified conditions and noted that each 24 years the number of living increased approximately threefold. Gumbel (1917) proved that the numbers of births and deaths and of the living in Euler's table were approaching a geometric progression and noted that many authors have thought that such was indeed the law of the increase in population under specified conditions.

**276.** Euler was in essence the first in [mathematical demography] to lay quite a distinct mathematical foundation for a number of the main demographic concepts such as the sequence of extinction (mortality tables), increase in population, the period of its doubling, etc. In addition, Euler formulated, as lucid as possible, the guiding rules for establishing the institution of private insurance (of life insurance in all its forms). Paevsky (1935, p. 103/2004, p. 8).

**277.** La manière de former les Tables de mortalité est fort simple. On prend dans les registres civils un grand nombre d'individus dont la naissance et la mort soient indiquées. On détermine combien de ces individus son morts dans la première année de leur âge, combien dans la seconde année, et ainsi de suite. On en conclut le nombre d'individus vivants au commencement de chaque année, et l'on écrit ce nombre dans la table à côté de celui qui indique l'année (pp. XCIX – C). … Tant de causes variables influent sur la mortalité, que les tables qui la représentent doivent changer suivant les lieux et les temps (p. CII). Laplace (1814/1886, pp. XCIX – C, CII).

**Comment.** This explanation seems superficial. In particular, if mortality tables ought to be changed *suivant les lieux et les temps*, how then is it possible to compile the data for many decades?

**278.** Ce n'est même que dans ces derniers temps que l'on a commencé a introduire dans les tables de mortalité la distinction des sexes. Quetelet & Smits (1832, p. 33).

**Comment.** The first to distinguish the sexes, although only for annuitants, was apparently Nicolaas Struick in 1740 (Hald 1990, p. 395).

**279.** [The origin of moral statistics] is usually attributed to 1741 – 1742, to the time when the first edition of Süssmilch's *Göttliche Ordnung* had appeared. Issues belonging to moral statistics were being touched even earlier. Süssmilch himself dealt above all with demography, but his right to be called the father of moral statistics cannot be questioned. … It was mostly the general reanimation of theoretical thinking, the first impetus to which was given by the brilliant school of French mathematicians, that led to the rebirth of moral statistics in the 1820s. The appearance of rich materials pertaining to some spheres of moral statistics, chiefly to criminal statistics, no less influenced the situation. France showed the way: from 1825 onwards, the *Compte général de l'administration de la justice* is being published there. Chuprov (1897, pp. 403, 404/2004, pp. 12, 13).

**280.** Das numerische Verhältnis der jährlichen Trauungen, der unehelichen Geburten, oder der Findelkinder varirt fast überall weniger als das der Sterblichkeit; die Sexualproportion der Gebornen, oder die jährliche Zahl der Zwillingsgeburten schwanken selbst bei ziemlich ansehnlichen Massen oft mehr als die der Trauungen oder der Verbrecher, die doch als freiwillige Handlungen zu betrachten sind. Chr. Bernoulli (1842/1849, p. 17).

**281.** Il est un budget qu'on paie avec une régularité effrayante, c'est celui des prisons, des bagnes et des échafauds … et, chaque année, les nombres sont venus confirmer mes prévisions, a tel point, que j'aurais pu dire, peut-être avec plus exactitude: Il est un tribut que l'homme acquitte avec plus de régularité que celui qu'il doit à la nature ou au trésor de l'État, c'est celui qu'il paie au crime! … Nous pouvons énumérer d'avance combien [crimes of each kind will be committed] à peu près comme on peut énumérer d'avance les naissances et les décès qui doivent avoir lieu.

La société renferme en elle les germes de tous les crimes qui vont se commettre, en même temps que les facilités nécessaires à leur développement. C'est elle, en quelque sorte, qui prépare ces crimes, et le coupable n'est que l'instrument qui les exécute. Tout état social suppose donc un certain nombre et un certain ordre de délits qui résultent comme conséquence nécessaire de son organisation. Quetelet (1836, t. 1, p. 12).

**Comment.** Quetelet became the father of moral statistics which at his time studied criminality, suicides and marriages. Süssmilch and even Graunt ought to be credited with its elements. Quetelet should have mentioned constancy relative to the number of adult (and not very old) men in the population. His last phrase meant that constancy may only be expected under invariable social conditions.

**282.** Statt verklingender Moralpredingten [about suicides] … habe ich sprechende Thatsachen angeführt, die ich der Aufmerksamkeit der hohen Regierungen dringend empfehle. Casper (1825, p. xi).

**283.** Sehen wir also zu, wie es um diese wundersame Constanz und Regelmäßigkeit in einer Reihe viel celebrirter Parade-Beispiele derselben in Wahrheit bestellt ist (p. 47). Da war es doch wahrlich reichlich früh und Quetelet unmotivirt rasch bei der Hand mit der Behauptung von der grausenerregenden Exactitüde in der Wiederkehr der Zahlen (p. 52). Rehnisch (1876, pp. 47, 52).

**Comment.** Quetelet had not studied the data for different years attentively enough, see Sheynin (1986, p. 300).

**284.** A Lottery … is properly a Tax upon unfortunate self-conceited fools. Petty (1662/1899, p. 64).

**285.** Une loterie est un impôt volontaire, parce qu'il est payé par ceux qui veulent jouer à ce jeu. Condorcet (1788/1847, p. 388).

**286.** On dit encore que cet impôt est volontaire. Sans doute il est volontaire pour chaque individu; mais, pour l'ensemble des individus, il est nécessaire; comme les marriages, les naissanceas et tous les effets variables sont nécessaires, et les mêmes à peu pres, chaque année, lorsqu'ils sont en grand nombre; en sorte que le revenu de la loterie est au moins aussi constant que les produits de l'agriculture. Laplace (1819/1912).

**287.** [On regularities in the sphere of moral statistics in general and on the stability of the revenue from lotteries.] Tels sont aussi les produits de la loterie de France, dont la suppression n'est pas moins désirable. Poisson (1837, p. 11).

**288.** Die Pest gab die Natur dem Oriente,/ unbillig ist sie nie;/ dafür gab sie dem Okzidente/ die Zahlenlotterie. Wahrschauer (1885, p. 12) as quoted by Bierman (1957, p. 651).

**Comment.** Bierman explains that that was a popular verse in Berlin (apparently for some time in the second half of the 18th century).

**289.** Die Interessante in den moral-statistischen Zahlen … ist überhaupt nicht ihre Stabilität, sondern ihre Veränderlichkeit. Lexis (1903a, p. 98).

**290.** Des Menschen Thaten und Gedanken … sind notwendig wie des Baumes Frucht, sie kann der Zufall gaukelnd nicht verwandeln. Schiller, *Wallenstein* (1800), Wallenstein's Tod, Aufzug 2, Auftritt 3.

**291.** The present value of a large marriage settlement is not [exactly] proportionate to that of a small marriage settlement. For a marriage settlement of 1,000 *zuz* that can be sold for a present value of 100, if the face value were 100, it could not be sold for 10, but rather for less. Maimonides, as quoted by Rabinovitch (1973, p. 164).

**Comment.** These settlements only stipulated payments for widows and divorced wives. An embryo of the notion of expectation is seen here.

**292.** Whereas we [in England] began to be busy in this direction about the middle of the 16$^{th}$ century (1560), Italy practised this civilising art of insurance [against accidents and against sickness] as early as the end of the 12$^{th}$ century. Guy (1885, p. 74).

**293.** Wide research concedes that Life Insurance came into its own, not by a front-door entrance, but by the marine insurance porthole. O'Donnell (1936, p. 78).

**294.** And whereas it hathe bene tyme out of mynde an usage amongste merchants, both of this realme and of forraine nacyons, when they make any great adventure … to give some consideracion of money to other persons … to have from them assurance made of their goodes, merchandizes, ships, and things adventured, … whiche course of dealinge is commonly termed a police of assurance. Publicke Acte No. 12 – An Acte Conc'ninge Matters of Assurances amongste Marchantes (1601). *Statutes* (1819, pp. 978 – 979).

**295.** I sacrificed my pay for his portrait [and the bronze spent] and he arranged that he would keep my money at 15 per cent during my natural life. Cellini (1965, § 80, p. 385).

**Comment.** Benvenuto Cellini was born in 1500; his autobiography covered the period up to 1562.

**296.** Eine päpstliche Bulle von 1423 erklärt schliesslich [after about a century of prohibition] den Rentenkauf für erlaubt. Du Pasquier (1910, pp. 484 – 485).

**297.** On peut envisager un peuple libre comme une grande association dont les membres se garantissent mutuellement leurs propriétés, en supportant proportionnellement les charges de cette garantie (p. CX). … Autant le jeu est immoral, autant ces établissements sont avantageux aux mœurs, en favorisant les plus doux penchants de la nature (p. CXI). Laplace (1814/1886, pp. CX, CXI).

**298.** Des tontines exercent deux penchants funestes: l'un est la disposition à attendre du hasard ce qui devrait être le fruit d'une industrie profitable à tous, ou le résultat ordinaire des institutions; l'autre est le désir d'augmenter ses jouissances personnelles en s'isolant du reste de la société. Fourier et al (1826/1890, p. 619).

**299.** No cases that, under the guise of insuring life, conceal deals in paying out some moneys upon the occurrence of a stipulated event not inflicting a [pecuniary] loss on the insured, – deals which do not restore actual damage, – should be attributed to insurance. They abuse the idea of insurance. Nikolsky (1895, p. 243).

**300.** [Actuarial science] owed more to the early workers in the Calculus of Finite Differences, especially to Newton, than to the early work in the Theory of Probability; I have often regretted that courses on the computational techniques of the Calculus of Finite Differences should be so seldom given to others than actuarial students. Fisher (1953, p. 1).

**301.** It was my actuarial work that brought me into contact with probability problems and gave a new turn to my mathematical interests. From 1918 on, I tried to get acquainted with the available literature on probability theory, and to do some work on certain problems connected with the mathematical aspect of insurance risk. These problems, although of a special kind, were in fact closely related to parts of probability theory that would come to occupy a central position in the future development … Cramér (1976, p. 511).

**302.** [The game pharao] où le banquier n'est qu'un fripon avoué, & le ponte une dupe … Buffon (1777/1954, p. 463).

**303.** You say that NN had won the main prize and all at once gained a fortune – so why won't I also win the same prize? You are unwittingly comparing yourself with the winner, and do not want to admit that it is much more natural to place yourself among the losers since they are by far more numerous. Ostrogradsky (1847/1961).

**304.** [La roulette] n'a ni conscience ni mémoire. Bertrand (1888a, p. XXII).

**305.** Parmi les joueurs, les unes choisissaient des numéros parce qu'ils n'étaient pas sortis depuis long-temps, d'autres choisissaient, au contraire, ceux qui sortaient le plus souvent. Ces deux préférences étaient également mal fondées. … Poisson (1837, pp. 69 – 70).

**Comment.** Montmort (1708, Préface) and De Moivre (1756, Preface) treated gamblers' prejudices in detail and Laplace (1814/1886, p. CXIIff) largely repeated their remarks. See also NNo. 303, 304.

### 3.3. Medical Statistics

**306.** Diese Art von Wahrscheinlichkeit, da wir das Verhältnis der Fälle selbst, erst durch einen wahrscheinlichen Schluss suchen müssen, nennet Rüdiger die medicinische Wahrscheinlichkeit, weil man in der Heilungskunst aus dem Verhältnisse derer, die an einer gegebenen Krankheit gestorben, oder durch ein gewisses Arzeneymittel genesen sind, zu der Zahl derjenigen bey welchen dieses nicht erfolgt ist, auf die Wahrscheinlichkeit in einzelnen vorkommenden Fällen schließt. Mendelsohn (not later than 1761, p. 204).

**Comment.** Not clear enough, but the idea is understandable. A reference to Rüdiger (Ridigeri; Leipzig, 1741) is supplied.

**307.** Except for the few high priests … the rest of the Esculapian train are nearly as ignorant [of statistics] as the ancients. Black (1788, p. 235).

**308.** Observations suivies et multipliées … peuvent nous apprendre des vérites utiles [about the way of life and avoidance of epidemics]. Condorcet (1795/1988, p. 542).

**309.** [On D'Alembert (1767b)] Une des choses qui m'ont le plus enchanté, c'est votre Mémoire sur l'inoculation. Il est plein de vues et de réflexions très-fines et très-exactes qui avaient échappé à tous ceux qui avaient déjà traité cette matière et qui la rendent tout à fait neuve et intéressante. A l'égard de vos difficultés sur le Calcul des probabilités, je conviens qu'elles ont quelque chose de fort spécieux qui mérite l'attention des philosophes plus encore que celle des géomètres, puisque de votre aveu même la théorie ordinaire est exacte dans la rigueur mathématique. Lagrange (1778/1882, pp. 87 – 88).

**Comment.** D'Alembert published his first memoir on inoculation back in 1761, but hardly Lagrange would have thought about it.

Daniel Bernoulli opened up that field but his memoir was only published in 1766. For D'Alembert, however, it was not difficult to acquaint himself with that memoir beforehand and to publish his criticism in advance, which disgusted Daniel.

**310.** La médecine systématique me parait … un vrai fléau du genre humain. Des observations bien multipliées, bien détaillées, bien rapprochées les unes des autres, voilà … à quoi les raisonnemens en médecine devraient se réduire (p. 163).

[A physician is] un Aveugle armé d'un baton … il lève son baton sans savoir où il frappe; s'il attrape la Maladie, il tue la Maladie; s'il attrape la Nature, il tue la Nature... Le médecin le plus digne d'être consulté, était celui qui croyant le moins à la medicine (p. 167). D'Alembert (1759/1821, pp. 163 and 167).

**Comment.** At best, this is an appeal to the numerical method of the first half of the 19th century, see No. 314.

**311.** En matière de statistique, c'est-à-dire dans les divers essais d'appréciation numérique des faits, le premier soin avant tout c'est de perdre de vue l'homme pris isolément pour ne le considérer que comme une fraction de l'espèce. … En médecine appliquée au contraire, le problème est toujours individuel (p. 173). La statistique mise en pratique, qui est toujours en définitive le mécanisme fonctionnant du calcul des probabilités, appelle nécessairement des masses infinies [?], un nombre illimité des faits non-seulement en vue d'approcher le plus près possible de la vérité, mais aussi afin d'arriver à faire disparaître, à éliminer, autant qu'il est possible, et à l'aide de procédés connus, les nombreuses sources d'erreurs si difficiles à éviter (p. 174). La condition des sciences médicales, à cet égard [of lending themselves to mathematization], n'est pas pire, n'est pas autre que la condition de toutes les sciences physiques et naturelles, de la jurisprudence, des sciences morales et politiques, etc. (p. 176). Poisson et al (1835, pp. 173, 174, 176).

**Comment.** Nowadays it seems strange that physical and natural sciences had been somehow separated.

**312.** Ce n'est qu'après avoir long-temps médité les leçons et les écrits de l'illustre géomètre [Poisson], que nous sommes parvenu à saisir l'étendue de cette question … de régulariser l'application de la méthode expérimentale à l'art de guérir. Gavarret (1840, p. xiii).

**313.** Le premier travail d'un observateur qui constate une différence dans les résultats de deux longues séries d'observations, consiste donc à chercher si l'anomalie n'est qu'apparente, ou si elle est réelle et accuse l'intervention d'une cause perturbatrice; il devra ensuite … chercher à déterminer cette cause … Ibidem, p. 194.

**Comment.** Gavarret was likely the first who clearly described the concept of null hypothesis. He was Poisson's student, graduated from the École Polytechnique, but then took to medicine.

**314.** [Discussions concerning the numerical method were held] uniquement de savoir si on remplacerait, par des rapports numériques les mots souvent, rarement, … etc. La méthode numérique, considérée sous ce point de vue rétréci, ne pouvait s'étendre au-delà d'une simple réforme dans le langage, mais il était impossible d'y voir une question de méthode scientifique et de philosophie générale. Ibidem, p. 10.
**Comment.** The numerical method was not much more than empiricism, see Sheynin (2009a, § 10.9).
**315**. [142 proposed questions about mammals, for example:] number of young at birth, number of pulsations per minute … temperature, average duration of life, proportion of males to females produced. [And about man:] quantity of air consumed per hour, of food necessary for daily support, average proportion of sickness amongst working classes. Babbage (1857, p. 294).
**Comment.** Babbage was a fervent collector of all kinds of statistical data. He even communicated to Quetelet (1869, t. 2, p. 226) how many drunkards did the London police force collect each month during the year 1832.
**316.** Quelques hommes distingués luttaient avec persévérance pour faire adopter l'emploi de la statistique en médecine. C'était, disaient-ils, le seul moyen de recueillir l'expérience des siècles en thérapeutique (p. x). [The members of the Paris Académie Royale de Médecine] n'avaient qu'une idée fort imparfaite de l'emploi du calcul en médecine (p. xiv). Gavarret (1840, pp. x, xiv).
**317.** Der erste, der die mathematische Statistik wesentlich anwandte, war, 1840, ein Mediziner, J. Gavarret. … hat Gavarret eine neue Epoche in der mathematischen Statistik eingeleitet. Vor Gavarret hatte man sich auf statistische Schlüsse mit einer Zuverlässigkeit von 99,999 …% beschränkt; man hatte sich darauf versteift, mit statistischen Argumenten Sicherheit zu erlangen. Die mathematische Statistik konnte sich erst recht entfalten, als man, seit Gavarret, lernte, sich mit einem vernünftigen Grad von Zuverlässigkeit zu begnügen. Freudenthal & Steiner (1966, pp. 181, 182).
**Comment.** I was unable to find anything of the sort in Gavarret (1840).
**318.** Wer kann sich übrigens wundern über die häufigen Widersprüche und den Wechsel der Heilmethoden, so lange man nicht numerische Belege von der relativen Wirksamkeit der Heilmittel in gleich vielen ganz analogen Fällen hat? (p. 7) … Nun gibt es deren freilich, die, obgleich bis in den höchsten Kreisen spukend, jede statistische Prüfung zur Lächerlichkeit machten (p. 7). Der Statistik und hiemit allen anthropologischen [anthropometric] Studien auf statistischem Wege liegt aber noch ob, und unstreitig ist diese Aufgabe die wichtigste und interessanteste, zu erforschen, welche Erscheinungen als konstante oder normale zu betrachten sind, welche Abweichungen als zufällige, und innert welchen Grenzen sie Schwankungen erleiden … (p. 13). Chr. Bernoulli (1842/1849, pp. 7, 7, 13).
**319.** Schon häufig ist uns praktischen Ärzten von theoretischer Seite in unwiderleglicher Weise gezeigt worden, dass alle unsere Schlüsse über Vorzüge oder Nachtheile einzelner Behandlungsmethoden, insofern sie auf die Statistik des thatsächlichen Erfolgs gegründet sind, vollständig in der Luft stehen, so lange wir nicht die strengen Regeln der Wahrscheinlichkeitsrechnung anwenden (p. 935).

Wenn die Ärzte bisher so selten von der Wahrscheinlichkeitsrechnung Gebrauch gemacht haben, so ist die Ursache davon weniger darin zu suchen, dass sie etwa dieser Disciplin nicht die gebührende Bedeutung beigelegt hätten; es beruht vielmehr hauptsächlich auf dem Umstand, dass bisher der analytische Apparat zu unvollkommen und unbequem war (p. 936). … Haben die Mathematiker gut sagen: Ihr Ärzte müsst, wenn Ihr sichere Schlüsse ziehen wollt, immer mit großen Zahlen arbeiten; Ihr müsst Tausende oder Hunderttausende von Beobachtungen zusammenstellen. Das ist ja eben bei therapeutischer Statistik gewöhnlich nicht möglich (p. 937). Wenn sie aber erfüllt sind, – dann kann es oft fraglich erscheinen, ob dann die Wahrscheinlichkeitsrechnung noch so dringend nöthig sei (pp. 937 – 938).

Gavarret nimmt, wie es Poisson bei einzelnen speciellen Aufgaben gethan hatte, einigermaßen willkürlich eine Wahrscheinlichkeit von 0.9953 oder 212/213, also eine Wahrscheinlichkeit welche 212 gegen 1 wetten lassen würde, als ausreichend an… (p. 939).

Wie nun, wenn dieser Grad der Wahrscheinlichkeit nicht erreicht ist? Sollen denn etwa alle Beobachtungsreihen, bei welchen die Wahrscheinlichkeit für Ausschließung des Zufalls nicht ganz 212 gegen 1 beträgt, vollkommen werthlos sein? … Und würden wir nicht selbst dann, wenn die Wahrscheinlichkeit des Vorzugs der einen Behandlungsweise nach den bis jetzt vorliegenden Erfahrungen nur 10 gegen 1 betrüge, doch für den nächsten Fall lieber diese als die andere anwenden? ... (p. 940). Liebermeister (ca. 1877, pp. 935, 936, 937, 937 – 938, 939, 940).

**Comment.** It is possibly Liebermeister rather than Gavarret (No. 317) who the *neue Epoche … eingeleitet*. Note however that beginning with 1863 (Chauvenet) and even earlier astronomers and geodesists had begun to offer tests for rejecting outliers quite in vein with Liebermeister's later reasoning. Incidentally, he could have strengthened his statement by mentioning previous practitioners who had made plausible inferences on the strength of scarce data (Bull 1959), and, in the first place, by citing Celsus (No. 463). Also see Descartes (No. 19) and J. Bernoulli (No. 12). Finally, I note that N. Bernoulli (1709/1975, p. 302) thought that an absentee ought to be declared dead once his death becomes only twice as probable as his still remaining alive.

See Seneta (1994) for a description of Liebermeister's own stochastic investigation.

**320.** Until recently, famous practitioners did not recognize medical statistics. Erisman (1887, vol. 1, p. 3).

**Comment.** Erisman included in his vol. 2 an Appendix on sanitary statistics, 184 pages long.

**321.** This disgraceful state of our Chatham Hospitals is only one more symptom of a system, which, in the Crimea, put to death 16,000 men – the finest experiment modern history has seen upon a large scale, viz., as to what given number may be put to death at will by the sole agency of bad food and bad air. Florence Nightingale; Kopf (1916, p. 390/1977, p. 312, without an exact reference).
**Comment.** Nightingale referred to the Crimean war of the mid-19[th] century. Pirogov also participated there and described the unsatisfactory situation in Russian military hospitals.
**322.** Any deaths in a people exceeding 17 in a 1,000 annually are unnatural deaths. If the people were shot, drowned, burnt, poisoned by strichnine, their deaths would not be more unnatural than the deaths wrought clandestinely by disease in excess of … seventeen in 1,000 living. Farr (ca. 1857; 1885, p. 148).
**Comment.** Farr added that 17 deaths in 1,000 were also too many because sanitary conditions were unsatisfactory even in localities in which mortality was lower.
It is hardly amiss to add Budd's statement (1849, p. 27) concerning cholera even if it does not bear on statistics:

*By reason of our common humanity, we are all the more nearly related here that we are apt to think. … And he that was never yet connected with his poorer neighbour by deeds of Charity or Love, may one day find, when it is too late, that he is connected with him by a bond which may bring them both, at once, to a common grave.*

**323.** Farr is rightly regarded as the founder of the English national system of vital statistics. For over forty years he supervised the actual compilation of English vital statistics, introduced methods of tabulation which have stood the test of time and a classification of causes of death which has been the basis of all subsequent methods. [His] life tables are still used in actuarial calculations and he formulated practical lessons as to the causation and prevention of disease … Newsholme (1931/1954).
**324.** Jahrhunderte lang glaubte man, und jetzt noch glauben Viele, an die besondere Gefährlichkeit der so genannten klimakterischen Jahre, an die vermehrte Sterblichkeit des Weibes zwischen 45 und 50 Jahren u. a. m. Genügend aber hat bereits die Statistik die Unstatthaftigkeit dieser Besorgnisse dargethan. Und ebenso ist vielfach schon die stereotype Meinung von der eminenten Schädlichkeit gewisser Umgebungen, der Frequenz mancher Krankheiten bei gewissen Beschäftigungen und Berufsarten u. a. zweifelhaft gemacht worden. Chr. Bernoulli (1842/1849, p. 7).
**325.** [The application of statistics] is in complete agreement with the spirit of surgery because the diseases included in the province of this science depend incomparably less on individual influences or modifications. N. I. Pirogov (1849a, p. 6).
**Comment.** Apart from epidemiology and public hygiene, surgery had indeed been the first medical science to yield to statistics.
**326.** [In general, statistical problems in surgery are] d'un abord plus facile que celles qui concernent la médecine. Quetelet (1846, p. 347).
**Comment.** Quetelet had not mentioned statistical studies of surgical diseases. And I cannot comment on the implied separation of surgery and medicine.

**327.** Die geheiligten Grundsätze der alten Schule, deren Ansichten im ersten Decennium dieses Jahrhunderts vorherrschen, sind durch die Statistik erschüttert – dass muss man ihr lassen, mit neuen Grundsätzen hat sie aber die alten nicht ersetzt … N. I. Pirogov (1864, p. 690).

**Comment.** For medicine, the new *Grundsätze* only appeared with the development of mathematical statistics. Poisson and Gavarret had taken the first steps in this direction even before the time of Pirogov, see NNo. 306, 312 – 314.

**328.** Even a slightest oversight, inaccuracy or arbitrariness makes [the data] far less reliable than the figures only founded on a general impression with which one is left after a mere but sensible observation of cases. … I … did not yet [in 1849] imagine all the blind alleys into which figures sometimes lead. N. I. Pirogov (1865 – 1866/1961, p. 20).

**Comment.** It follows that Pirogov was compelled to adhere to the first stage of the statistical method (NNo. 461, 463). Incidentally, he was unable to talk things over with statisticians.

**329.** Ich neige zur Ansicht, dass die Letalität jedes pathologischen Prozesses, jeder Verwundung und jeder Operation, die verschiedensten Verhältnisse ungeachtet, doch im Ganzen etwas Konstantes und bestimmtes sein muß. Pirogov (1864, p. 5).

**330.** Während meines Aufenthalts im Ausland hatte ich mich zur Genüge davon überzeugen können, dass die wissenschaftliche Wahrheit bei weitem nicht der Hauptzweck selbst berühmter Kliniker und Chirurgen ist. Ich hatte mich hinlänglich davon überzeugt, dass in berühmten klinischen Anstalten gar oft Massregeln nicht zur Enthüllung, sondern zur Verdeckung der wissenschaftlichen Wahrheit getroffen wurden. Pirogov (1884 – 1885/1894, p. 453).

**331.** Wovon hängen die Erfolge der Behandlung oder die Mortalitätsverminderung in den Armeen ab? Doch gewiss nicht von der Therapie und der Chirurgie an und für sich. … Für die Massen steht von der Therapie und Chirurgie – ohne eine tüchtige Administration – auch in Friedenzeiten wenig Nutzen zu erwarten; um so weniger also noch bei solchen Katastrophen wie ein Krieg. N. I. Pirogov (1871, p. 439).

**Comment.** Pirogov himself was a really good administrator of military medicine. Administration on a large scale is now studied by operational research.

**332.** Ich in der Schule wohl der beste Schüler in Geschichte und russischer Literatur, nicht aber in der Mathematik gewesen bin. Dabei glaube ich aber wohl zu behaupten zu dürfen, dass in mir etwas von mathematischer Ader steckte; sie entwickelte sich aber, glaube ich, nur sehr langsam mit meinem fortschreitenden Alter, und als ich das sogar sehr lebhafte Verlangen fühlte, etwas von Mathematik zu verstehen, da war es schon zu spät dazu. N. I. Pirogov (1884 – 1885/1894, p. 144).

**333.** The development of medical statistics had begun exactly after mankind had to convince itself too clearly and too bitterly of the utter helplessness of medicine against such of its evils as cholera, [typhoid?] fever, etc. Peskov (1874, p. 10).
**Comment.** Surgery forestalled epidemiology.
**334.** Medical statistics should have at its disposal … mean values as accurate as those used in meteorology so as to make possible the construction of lines of equal sickness, mortality, etc. and thus to discover the laws of sickness. Peskov (1874, p. 89).
**335.** Approximately half the articles published in medical journals that use statistical methods use them incorrectly. Glantz (1980, p. 1) as quoted by Gaither et al (1996).
**336.** La probabilité des mathématiciens … n'est guère que la théorie du hasard. Invoquer la probabilité … c'est donc invoquer le hasard; c'est renoncer à toute certitude médicale, à toute règle rationnelle tirée des faits propres de la science. … La médecine ne sera plus un art, mais une loterie (p. 14). [To apply the theory of probability] aux faits réels du monde physique et moral [is] inutile ou illusoire (p. 15). [Probability provided] solutions ou nuisibles, ou insuffisantes, ou trompeuses. … Son importation en médecine est anti-scientifique, abolissant … la véritable observation (p. 31). D'Amador (1837, pp. 14, 15, 31).
**Comment.** It is seen here and in many other instances (below) that some authors did not understand either the dialectical connection between the random and the necessary (No. 792), or that statistical conclusions were not applicable to individuals.
**337.** Je me suis convaincu par expérience, comme les résultats sont différents entre les opérations faites dans les petits établissements cliniques, et les opérations exécutées dans les grands hôpitaux; et même, combien la différence est grande dans les résultats obtenues par les opérations dans les différents hôpitaux de la même ville, exécutées dans les conditions exactement semblables en apparence. [An important feature was] la constitution générale d'un hôpital. N. I. Pirogov (1849c, p. 191).
**338.** The origin and spread of fever in a hospital, or the appearance and spread of hospital gangrene, epysipelas and pyaemia generally are much better tests of the defective sanitary state of a hospital than its mortality returns. Fl. Nightingale (1859/1963, p. 10).
**Comment.** [If a hospital only admits serious cases] A large mortality may even be considered as a presumption of an hospital being well conducted … Blane (1813/1822, p. 140).
**339.** Cette prétendue application de ce qu'on appelle la statistique à la médecine, dont plusieurs savants attendent des merveilles, et qui pourtant ne saurait aboutir, par sa nature, qu'à une profonde dégénération directe de l'art médical, dès lors réduit à d'aveugles dénombrements. Une telle méthode, s'il est permis de lui accorder ce nom, ne serait réellement autre chose que l'empirisme absolu, déguisé sous de frivoles apparences mathématiques. Comte (1830 – 1842/1877, t. 3, No. 40, p. 329).
**Comment.** In part, and without recognizing it, Comte criticized the numerical method.

**340.** The only problem was to compare the probabilities of a floating kidney, a chronic catarrh and appendicitis. The problem was not about Ivan Ilyich's life. Tolstoy (1884 – 1886/2003, p. 27).
**Comment.** Again a criticism of the same method.

### 3.4. Administration of Justice

**341.** The witness [in law-suits pertaining to loans] to whom, within seven days after he has given evidence, happens [a misfortune through] sickness, a fire, or the death of a relative, shall be made to pay the debt and a fine. Bühler (1886, § 108).
**Comment.** This is an early statement actually concerning errors of the first and second kind. The words in brackets were added by Bühler.

**342.** Laplace rejeté les résultats de Condorcet, Poisson n'a pas accepté ceux de Laplace; ni l'un ni l'autre n'a pu soumettre au calcul ce qui y échappe essentiellement: les chances d'erreur d'un esprit plus ou moins éclairé, devant des faits mal connues et des droits imparfaitement définis. Bertrand (1888a, p. 319).
**Comment.** Since he mentioned *droits*, Bertrand apparently had in mind civil cases. Suchlike entirely negative, and sometimes even wrong statements (as in this instance) were characteristic of Bertrand, see NNo. 120, 365, 737.

**343.** Any one of us would prefer to pass a sentence acquitting a wrong-doer rather than condemning as guilty one who is innocent … Aristotle, Problemata 951b0.

**344.** Viel besser ist, zehen schuldige zu befreuen, als einen unschuldigen zum Tode zu condemniren. Russia's *Kriegs-Reglement* (1716/1830, p. 403).

**345.** [In accord with the idea of Laplace [1816/1886, p. 521], as Poisson remarked, there should be] plus de danger pour la sûreté publiques à l'acquittement d'un coupable, que de crainte de la condemnations d'un innocent … Poisson (1837, p. 6).

**346.** Il est manifeste que les conditions de majorité, de pluralité imposées aux décisions d'un corps judiciaire ou d'une assemblée délibérante, doivent avoir des relations avec la théorie mathématique des chances. Cournot (1843, § 192).

**347.** Il y a … causes d'erreur qui sont de nature à influer en même temps sur tous ceux qui prendront connaissance de l'affaire … Ibidem, § 206.
**Comment.** Cournot devoted a chapter of the same contribution to an attempt to study legal proceedings by taking into account the dependence between the decisions made by judges (jurors) and he at least traced some elementary steps in this direction.

**348.** L'accusé, quand il arrive à la cour d'assises, a déjà été l'objet d'un arrêt de prévention et d'un arrêt d'accusation, qui établissent contre lui une probabilité plus grande que 1/2, qu'il est coupable … Poisson (1837, p. 4).
**Comment.** The introduction of a measure of prior culpability into stochastic formulas was indeed reasonable, but Poisson should have also stressed that the presumption of innocence of any individual had still persisted.

**349.** [The personality of the accused and the probability of his being convicted; extract. Quetelet (1836, t. 2, p. 313 and in earlier contributions of 1832 and 1833)]

**Table.** Meaning of columns explained above

| 1 | Ayant une instruction supérieure | 0.400 |
| 3 | Accusé de crimes contre les personnes | 0.477 |
| 5 | Étant femme | 0.576 |
| 6 | Ayant plus de 30 ans | 0.586 |
| 8 | Sans désignation aucune | 0.614 |
| 9 | Étant homme | 0.622 |
| 10 | Ne sachant ni lire ni écrire | 0.627 |
| 11 | Ayant moins de 30 ans | 0.630 |
| 12 | Accusé de crimes contre les propriétés | 0.655 |

**Comment.** The probabilities in lines 8, 9, 10 and 11 are almost the same and it is doubtful whether they were really different.
**350.** Les tables de criminalité, pour les différents âges, méritent au moins autant de confiance que les tables de mortalité. Quetelet (1842, p. 14).
**351.** Especially surprising … is … that the idea of analysing, for example, 'repeated offenders' completely eluded [Quetelet]. … It seems likely that [this was] due …at least in part to the inadequacies of the data available at the time. Landau & Lazarsfeld (1978, p. 831).
**352.** [The study of testimonies is] the weakest section of the theory of probability. Markov (1915, p. 32).
**353.** Les lois du hasard ne s'appliquent pas à ces questions. … Des causes multiples entrent en action, elles troublent les hommes, les entraînent à droite et à gauche, mais il y a une chose qu'elles ne peuvent détruire, ce sont leurs habitudes de moutons de Panurge. Poincaré (1896/1912, p. 22).

**Comment.** Without explaining their restriction, Condorcet, Laplace and Poisson studied the ideal case of independent witnesses and judges (jurors) although Laplace (1816/1886, p. 523) once, and only in passing, mentioned this supposition. Nevertheless, Gauss' opinion (No. 356) should not be forgotten. And, anyway, the efforts of these scholars turned public attention to the problems inherent in legal proceedings. The interest in the application of probability to this subject has been rekindled (Heyde & Seneta 1977, p. 34); to those authors who are mentioned in that source I can now add Zabell (1988), Gastwirth (2000) and Dawid (2005). The last-mentioned author stressed the importance of the interpretation of indirect stochastic information.

**354.** L'application du calcul des probabilités aux sciences morales est, comme l'a dit je ne sais plus qui, le scandale des mathématiques parce que Laplace et Condorcet, qui calculaient bien, eux, sont arrivés à des résultats dénuées de sens commun! Rien de tout cela n'a de caractère scientifique … Poincaré, letter of ca. 1899; *Procès* (1900, t. 3, p. 325).

**Comment.** Poincaré wrote this letter in connection with the notorious Dreyfus case. His next statement (No. 355) was also caused by the same case. Condorcet certainly preceded Laplace.

**355.** L'application du calcul des probabilités à ces matières [to studying handwritings] n'est pas légitime. Poincaré et al (1906, p. 245).

**356.** Kann die Wahrscheinlichkeitsrechnung dem Gesetzgeber eine Richtschnur für die Bestimmung der Zahl der Zeugen und der Ritter geben, wenn sie auch für den einzelnen Fall nichts lehrt. Gauss, in 1841, as reported by W. E. Weber; Gauss, *Werke*, Bd. 12, pp. 201 – 204.

**357.** Äusserst nützlich ist die Abschätzung von Wahrscheinlichkeiten, wiewohl in Beispielen aus dem Rechts- und Staatswesen nicht so sehr eine subtile Rechnungsart nötig ist als eine genaue Aufzählung aller Umstände. Leibniz, letter to J. Bernoulli 3 Dec. 1703, in Latin. Gini (1946, p. 405).

**358.** A very slight improvement in the data by better observations, or by taking into fuller consideration the special circumstances of the case, is of more use than the most elaborate application of the calculus of probabilities founded on the data in their previous state of inferiority. The neglect of this obvious reflection has given rise to misapplications of the calculus of probability which have made it the real opprobrium of mathematics. It is sufficient to refer to the applications made of it to the credibility of witnesses, and to the correctness of the verdicts of juries. Mill (1843/1886, p. 353).

**Comment.** The second half of this statement is a one-sided opinion, its first half is weak: no sound judgement can be based on inferior data. Concerning the consideration of circumstances see also comment to No. 353.

### 3.5. Quetelet (apart from His Appearance Elsewhere)

**359.** Les plantes et les animaux sont restés tel qu'ils sortis de la main du créateur. Quelques espèces, à la vérité, ont disparu, et d'autres se sont montrées successivement. Quetelet (1846, p. 259).
**Comment.** Quetelet never ever referred to Darwin.
**360.** Die … Pflanzenformen seien nicht ursprünglich determinirt und festgestellt … ihnen sei vielmehr … eine glückliche Mobilität und Biegsamkeit verliehen, um in so viele Bedingungen … sich zu fügen und darnach bilden und umbilden zu können. Goethe (1831/1891, p. 120).
**361.** Mais tout d'espèces d'animaux éteintes dont M. Cuvier a su reconnaître … l'organisation dans les nombreux ossements fossiles qu'il a décrits, n'indiquent-elles pas dans la nature une tendance à changer les choses mêmes les plus fixes en apparence? Laplace (1796/1884, p. 480).
**362.** Je ne crains pas de suivre l'exemple … de Humboldt; il m'ofrit, par un mot emprunté à la langue grecque, le titre de mon ouvrage *Anthropométrie* [1871] … Quetelet (1870, p. 671).
**363.** L'homme que je considère ici est dans la société l'analogue du centre de gravité dans les corps; il est la moyenne autour de laquelle oscillent les éléments sociaux. Ce sera … un être fictif … Quetelet (1832a, p. 4).
**364.** Si l'homme moyen était déterminé pour une nation, il présenterait le type de cette nation … Quetelet (1832b, p. 1).
**Comment.** Quetelet set too high a store by his Average man and he defined his concept not clearly enough. By average he sometimes meant the arithmetic mean, and in other cases, the median. This, however, is the mildest possible reproach (Sheynin 2009a, § 10.5).
**365.** Dans le corps de l'homme moyen, l'auteur belge place une âme moyenne. L'homme type sera donc sans passions et sans vices, ni fou ni sage, ni ignorant ni savant, souvent assoupi: c'est la moyenne entre la veille et la sommeil, ne répondant ni oui ni non; médiocre en tout. Après avoir mangé pendant trente-huit ans la ration moyenne d'un soldat bien portant, il mourrait, non de vieillesse, mais d'un maladie moyenne que la Statistique révélerait pour lui. Bertrand (1888a, p. XLIII).
**Comment.** It was absolutely wrong that the Average man was *sans passions et sans vices*: while introducing mean inclinations to crime and to marriage, Quetelet attributed these concepts to him. Furthermore, the Average man is convenient at least as the average producer and consumer. Fréchet (1949) replaced him by the closely related *typical* man.
**366.** [Concerning Quetelet] Einen gedankenreichen, aber unmethodischen und daher auch unphilosophischen Geist. Knapp (1872, p. 124).
**367.** Quetelet's promises and hopes and his achievements in 1835 – 1836 remained in statu quo up to the last edition of [*the Physique sociale*] in 1869. He achieved nothing hardly of real value in all those 33 years. Galton, letter of 1891; Karl Pearson (1914 – 1930, 1924, p. 420).
**368.** Er [Quetelet] hat großes Verdienst in der Vergangenheit, indem er nachwies, wie selbst die scheinbaren Zufälle des sozialen Lebens durch ihre periodische Rekurrenz und ihre periodischen Durchschnittszahlen eine innere Notwendigkeit besitzen. Aber die Interpretation dieser Notwendigkeit ist ihm nie gelingen. Er hat auch keine Fortschritte gemacht, nur das Material seiner Berechnung ausgedehnt. Er ist heut nicht weiter, als er vor 1830 war. Marx, letter of 1869; Marx (1952, pp. 81 – 82).

**369.** Er was een statistisch bureau, er waren statistieken en statistici, maar er war nog geen [but there was no] statistiek [before Quetelet]. Freudenthal (1966, p. 7).

**370.** Il est peu de pays, croyons-nous, où le calcul des probabilités tienne une place aussi considérable dans l'enseigné supérieur qu'en Belgique. … En France … il n'est enseigné qu'accidentellement, comme accessoire du cours de physique mathématique à la Sorbonne: à l'École Polytechnique, on ne lui consacre que quelques leçons des cours d'analyse et d'astronomie. En Allemagne, …la théorie de la compensation des erreurs d'observation fait souvent l'objet d'une Vorlesung spéciale, mais rarement le calcul des probabilités y est exposé dans toute son étendue. Mansion (1905, p. 3).

**Comment.** The author explained the situation by Quetelet's lasting influence. Concerning mathematical physics see No. 208, Comment.

**371.** You are indeed a worthy example in activity & power to all workers in science and if I cannot imitate your example, I can at least appreciate & value it. Faraday, Letter No. 1367 of 1841 to Quetelet; Faraday (1996, p. 42).

**Comment.** Several other letters written by Faraday (1996 – 1999) are also relevant.

    1) Letter No. 1862 of 1846 (1996, p. 501). F. had been receiving Quetelet's *most important exertions*.

    2) Letter No. 2139 of 1848 (1996, p. 743). F. was *astonished to think what must be* [Quetelet's] *industry*.

    3) Letter No. 2263 of 1850 (1999, p. 126). F. studied Quetelet's *results in atmospheric electricity. They are*, [he] *think*[s], *very admirable; and* [he] *admire*[s] *the truly philosophical spirit in which* [Quetelet had] *been content to give them, without any addition of imagination or hypothesis*.

    4) Letter No. 2379 of 1851, this time to Richard Taylor, the then Editor of the *Lond., Edinb. and Glasgow Phil. Mag.* (1999, p. 270). Quetelet's observations of *phenomena of atmospheric electricity* are *so important …*

    5) Letter No. 2412 of 1851 (1999, p. 281). Quetelet's observations of atmospheric electricity are *most imteresting*.

**372.** Wenn mir die Statistik sagt, dass ich im Laufe des nächsten Jahres mit einer Wahrscheinlichkeit von 1 zu 49 sterben, mit einer noch größeren Wahrscheinlichkeit schmerzliche Lücken in dem Kreis mir theurer Personen zu beklagen haben werde, so muss ich mich unter den Ernst dieser Wahrheit in Demuth beugen; wenn sie aber, auf ähnliche Durchschnittszahlen gestützt, mir sagen wollte, dass mit einer Wahrscheinlichkeit von 1 zu so und so viel [I shall commit a crime] so dürfte ich ihr unbedenklich antworten: ne sutor ultra crepidam! [Cobbler! Stick to your last!]. Rümelin (1867, p. 25).

**373.** [Quetelet's worshippers' *zealous beyond reasoning*] naïve admiration for 'statistical laws'; their idolizing of stable statistical figures; and their absurd teaching that regarded everyone as possessing the same 'mean inclinations' to crime, suicide, and marriage, undoubtedly provoked protests. Regrettably, however, the protests were hardly made in a scientific manner. Chuprov (1909/1959, p. 23).

### 3.6. Statistics and Mathematics

**374.** Chaque branche des connaissances humaines a besoin de la méthode statistique, qui, dans le fait, n'est que la méthode numérique (p. 333). [This science] consiste à savoir réunir les chiffres, les combiner et les calculer, de la manière la plus propre à conduire à des résultats certaines. Mais ceci n'est, à proprement parler, qu'une branche des mathématiques (p. 334). Alph. De Candolle (1833, pp. 333, 334).

**375.** [About writings on social sciences before Quetelet] in de ene die werken, waarin te veel formules en te weinig (of in't geheel geen) [if any] empirische cijfers voorkomen; de andere [works] met veel cijfers en weinig of geen [if any] wiskunde … Freudenthal (1966, p. 7).

**Comment.** Poisson's and partly Cournot's contributions might have been attributed to the works of the first type. It seems that Quetelet (No. 285) had appraised the then existing situation no less correctly.

**376.** Les sciences statistiques ne feront de véritables progrès, que lorsqu'elles seront confiées à ceux qui ont approfondi les théories mathématiques. Fourier, letter of unknown date to Quetelet; Quetelet (1826, p. 177).

**Comment.** Quetelet is known to have greatly influenced the work of the *Congress international de statistique*, and here is a pertinent statement (*Congrès* 1868, p. 6): "statistical questions" should have their "scientific base" in mathematics, and should be studied "in direct connection with the theory of probability".

**377.** [Poisson] exprimait parfois dans sa correspondance [with Quetelet] avec une sévérité narquoise et peu rassurante pour les statisticiens qui prétendaient subsister leurs fantaisies aux véritables principes de la science. Quetelet (1869, t. 1, p. 103).

**378.** The most sublime problems of the arithmétique sociale can be only resolved with the help of the theory of probability. Libri-Carrucci et al (1834, p. 535).

**Comment.** The forgotten term meant demography, medical statistics and actuarial science. I was unable to locate that source and do not have its original French text anymore.

**379.** La théorie des probabilités prit naissance presque en même temps que la statistique, sa sœur puînée, dont elle devait devenir la compagne la plus sûre et la plus indispensable. Cette concordance n'est point accidentelle, mais l'une de ces sciences interroge en quelque sorte par ses calculs et coordonne ce que l'autre obtient par ses observations et ses expériences. Quetelet (1869, t. 1, p. 134).

**380.** La plupart des observateurs, les meilleurs même, ne connaissent que très-vaguement, je ne dirai pas la théorie analytique des probabilités, mais la partie de cette théorie qui concerne l'appréciation des moyennes. Quetelet (1846, p. 63).

**381.** [Statistics] steht in gar keiner inneren Beziehung [with mathematics]. Die Mathematik beruht eben auf der Deduktion, die Statistik auf Induktion. Haushofer (1872, pp. 107 – 108).

**Comment.** Haushofer had thus explained his refusal to apply probability to statistics. However, statistics is only partly based on induction, and, in addition, the law of large numbers is a bridge between deduction and induction. Nevertheless, Knapp (1872, pp. 116 – 117) declared that for statistics that law *ist von geringerer Bedeutung* because statisticians always make only one observation (as when counting the population of a city).

**382.** Man braucht mehr als nur die Urnen des Laplace mit bunten Kugeln zu füllen um eine wissenschaftliche Statistik herauszuschütteln. Knapp (1872, p. 117).

**Comment.** Indeed, statistics is not reduced to probability theory. Knapp, however, opposed the latter (Sheynin 1990b/1996, p. 35).

**383.** Die mathematische Statistik ist kein Automat, in den man nur das statistische Material hineinzustecken hat, um nach einigen mechanischen Manipulationen das Resultat wie an einer Rechenmaschine abzulesen. Charlier (1920, p. 3).

**384.** Statistik spielende Mathematiker können nur durch mathematisch ausgerüstete Statistiker überwunden werden. Chuprov (1922b, p. 143).

**385.** [Statistics is] the application of mathematical theory to the interpretation of mass observations. K. Pearson (1978, p. 3).

**386.** Mathematics may be compared to a mill … but … what you get depends on what you put in. … pages of formulae will not get a definite result out of loose data. T. H. Huxley (1869, p. L).

**Comment.** Huxley could have mentioned statistics as well.

**387.** A statistician ought to become a mathematician because his science is a mathematical science. Slutsky (1912, p. 3/2009, p. 15).

**388.** The mathematician, the statistician, and the philosopher do different things with a theory of probability. The mathematician develops its formal consequences, the statistician applies the work of the mathematician and the philosopher describes in general terms what this application consists in. The mathematician develops symbolic tools without worrying overmuch what the tools are for; the statistician uses them; the philosopher talks about them. Each does his work better if he knows something about the work of the other two. Good (1959, p. 443).

**389.** The science of statistics is essentially a branch of Applied Mathematics and may be regarded as mathematics applied to observational data.

Statistics may be regarded as (i) the study of populations; (ii) as the study of variation; (iii) as the study of methods of the reduction of data. Fisher (1925/1990, pp. 1 – 2).

**390.** If mathematical statistics, as justifiably accepted, is a science of the mathematical methods of studying mass phenomena, then the theory of probability should be considered its organic part (p. 216/249). The system of the main concepts of theoretical statistics is still … in the making. Only gradually does this discipline cease to be the *applied theory of probability* in the sense of a collection of separate stochastic problems unconnected with each other by general ideas (p. 218/251). Kolmogorov (1948b, pp. 216, 218/2005, pp. 249, 251).

**391.** We have for a long time cultivated a wrong belief in the existence, in addition to mathematical statistics and statistics as a socio-economic science, of something like yet another non-mathematical, although universal *general* theory of statistics which essentially comes to mathematical statistics and some technical methods of collecting and treating statistical data. Accordingly, mathematical statistics was declared a part of this *general theory of statistics*. Such views … are wrong. …

All that which is common in the statistical methodology of the natural and social sciences, all that which is here indifferent to the specific character of natural or social phenomena, belongs to … mathematical statistics. A. N. Kolmogorov (Anonymous 1954, p. 47/Sheynin (2005a, pp. 254, 255).

**392.** *Statistics, mathematical*, the branch of mathematics devoted to the mathematical methods for the systematization, analysis, and use of statistical data for the drawing of scientific and practical inferences. Here, the term *statistical data* denotes information about the number of objects in some more or less general set which possess certain attributes … Kolmogorov & Prokhorov (1974/1977, p. 721).

**Comment.** Knies (1850, p. 163) was apparently the first to introduce the term *mathematical statistics* (in German).

**393.** [(The notion of) mathematical statistics] is a grotesque phenomenon. Anscombe (1967, p. 3 note).

**394.** In a number of cases mathematical statistics becomes a purely formal doctrine, entirely torn away from concrete reality (p. 64). … The theory of probability cannot be considered as the foundation of the statistical science (p. 77). Nekrash (1947, pp. 64, 77).

**395.** Its only justification lies in the help it can give in solving a problem. … statistical theory is not a branch of mathematics. Statistics, like engineering, requires all the help it can receive from mathematics; but … mathematical statistics as a separate discipline cannot simply exist.

The aim of statistics is to reach a decision on a probabilistic basis, on available evidence. If the problem is of theoretical nature, statistics supplies a valid method for drawing general conclusions from particular experience. If the problem is a practical one, statistics supplies the basis for choosing a particular course of action … by balancing the risks of gain and loss. Mahalanobis, 1950, as quoted by Rao (1993, p. 339).

**396.** Apart from his work in India, Mahalanobis will be remembered as one of the pioneers who along with Karl Pearson, R. A. Fisher, J. Neyman and A. Wald laid the foundations of statistics as a separate discipline. Ibidem, p. 337.

**397.** Diejenigen statistischen Grössen, die sich in die Schema der Wahrscheinlichkeitsrechnung am ehesten fügen lassen, zugleich solche sind, die meistens kein grosses materiell-statistisches Interesse in Anspruch nehmen können. Bortkiewicz (1894 – 1896, p. 356).
**398.** In Wirklichkeit operiert in der Statistik auch der grimmigste Feind der Analogie mit den Zufallsspielen mit Vorstellungen, die gerade diesem Erscheinungsgebiete entstammen. … Er wendet also die Wahrscheinlichkeitstheorie an, ohne es zu wollen und ohne es zu wissen, und darum in unmethodischer Weise, nach der rohen Art des reinen Empirikers. Bortkiewicz (1904, pp. 251 – 252).
**399.** Such procedures as the construction of curves of densities, "smoothing" of series, interpolation not only do not help in revealing the real nature of the studied phenomena, but, to the contrary, can provide a distorted notion about them. … The so-called correlation method … adds nothing in essence to the results of an elementary analysis. Kaufman (1922, p. 152).
**Comment.** The posthumous edition of this source (Moscow, 1928) with Kaufman's name still indicated on the title-page carries a contrary statement written by V. I. Romanovsky: The theory of correlation is "one of the most important and amazing chapters of modern statistics".
**400.** [The] Schema der Wahrscheinlichkeitsrechnung [is] auch die höchste wissenschaftliche Form in welcher die Statistik ihren Stoff fassen kann. Lexis (read 1874, 1903c, p. 241).
**401.** [The only aim of applying probability is to obtain] ein verständliches Schema für die Verteilung der Fälle und … einen Massstab für die Stabilität der statistischen Verhältniszahlen zu bieten. Lexis (1903b, p. 230).
**402.** Like you [N. S. Chetverikov], I even now, just as [in 1909], … see no possibility of throwing a formal logical bridge across the crack separating frequency from probability … Chuprov, letter dated 1923; Sheynin (1990b/1996, p. 97).
**Comment.** Chuprov could have recalled the strong law of large numbers.

### 3.7. Application in Science; General Considerations

**403.** Her [Nightingale's] statistics were more than a study, they were indeed her religion. For her Quetelet was the hero as scientist, and the presentation copy of his *Physique sociale* is annotated by her on every page. [She] believed – and in all actions of her life acted upon the belief – that the administrator could only be successful if he were guided by statistical knowledge. The legislator – to say nothing of the politician – too often failed for want of this knowledge. Nay, she went further; she held that the universe – including human communities – were evolving in accordance with a divine plan; that it was man's business to endeavour to understand this plan and guide his actions in sympathy with it. But to understand God's thoughts, she held, we must study statistics, for these are the measure of His purpose. Thus, the study of statistics was for her a religious duty, pp. 414 – 415.

[Statistics is] the most important subject in the whole world: for upon it depends the practical application of every other science and of every art: the one science essential to all political and social administration, all education, all organisation based on experience, for it only gives the results of our experience. Fl. Nightingale, in a letter written after 1856 (Pearson 1914 – 1930, vol. 2, pp. 414 – 415, 415), as quoted by Tee (2002, p. 4).

**404.** Il y a dans chaque science, art ou objet d'étude, une partie statistique, soit numérique. On la retrouve en agriculture, en médecine, dans toutes les branches de l'administration, dans les sciences physiques, naturelles, et jusque dans les sciences morales. Elle occupe une très grande place dans la géographie botanique. Pour moi, j'en conviens, j'aime les chiffres autant que d'autres les détestent; mais ce qui me plaît, ce n'est pas d'accumuler des chiffres, c'est de montrer à quel degré il est nécessaire de choisir convenablement les valeurs, de les discuter, en d'autres termes, de les subordonner aux lois de la logique et du bon sens. Alph. De Candolle (1855, p. xvi).

**405.** La statistique a pris un développement en quelque sorte exubérant; et l'on n'a plus qu'à se mettre en garde contre les applications prématurées et abusives qui pourraient la discréditer pour un temps, et retarder l'époque si désirable où les données de l'expérience serviront de bases certaines à toutes les théories qui ont pour objet les diverses parties de l'organisation sociale. Cournot (1843, § 103).

**406.** It is beyond hope to rectify bad observations having a large variance by increasing their number. Tiurin (1975, p. 59).

**407.** In den Anwendungen der Wahrscheinlichkeitsrechnung sehr gefehlt werden könne, wenn man nur auf die Zahlen bauet, welche wiederholte Beobachtungen geben, und nicht jeder andern Kenntnis, die man sich von der Natur der Sache und deren Verhältnissen verschaffen kann; ihr Recht widerfahren lässt, so schwer dies oft auch sei. Gauss as reported by W. E. Weber, 1841; Gauss, *Werke*, Bd. 12, pp. 201 – 204.

**408.** The data that I have adduced … have been objected to on the ground that they are collected from too many different hospitals, and too many sources. But … I believe all our highest statistical authorities will hold that this very circumstance renders them more, instead of less, trustworthy. J. W. Simpson (1847 – 1848/1871, p. 102).

**Comment.** This was an unfortunate example of applying the statistical method in the presence of statistical data, cf. No. 409, which originated with Graunt in statistics and Tycho in astronomy. Later Simpson himself coined the still existing term *Hospitalism* and noted that mortality in hospitals increased with the number of beds (1869 – 1870/1871, p. 399).

**409.** Si l'on recueille l'urine d'un homme pendant vingt-quatre heures et qu'on mélange toutes les urines pour avoir l'analyse de l'urine moyenne, on a précisément l'analyse d'une urine qui n'existe pas; car à jeun l'urine diffère de celle de la digestion, et ces différences disparaissent dans le mélange. Le sublime du genre a été imaginé par un physiologiste qui, ayant pris de l'urine dans un urinoir de la gare d'un chemin de fer où passaient des gens de toutes les nations, crut pouvoir donner ainsi l'analyse de l'urine *moyenne* européenne! Bernard (1865/1926, t. 2, pp. 117 – 118).
**Comment.** Bernard denied the application of statistics to medicine, but at least he had not repeated Simpson's mistake (No. 408).
**410.** There is no science which has not sooner or later discovered the absolute necessity of resorting to figures as measures and standards of comparison; nor is there any sufficient reason why physiology and medicine should claim exemption. … On the contrary, they belong in an especial manner to the class of sciences which may hope to derive the greatest benefit from the use of numbers (p. 801). Without statistics a science bears to true science the same sort of relation which tradition bears to history (p. 802). Guy (1852, pp. 801, 802).
**411.** A knowledge of statistics is like a knowledge of foreign languages or algebra; it may prove of use at any time under any circumstances. Bowley (1901/1946, p. 4).
**412.** A theory can be proved by experiment, but no path leads from experiment to the birth of a theory. Einstein (1976) as quoted by Gaither et al (1996).
**413.** We are to admit no more causes of natural things than such as are both true and sufficient to explain their appearances. Newton (1687/1960, Book 3, Rule 1).
**Comment.** Cf. the Talmud (Rabinovitch 1973, p. 120): "What we see we may presume; but we presume not what we see not".
**414.** We ought … to resolve that which is common to all the stars … into a single real motion of the Solar system, as far as that will answer the known facts, and only to attribute to the proper motions of each particular star the deviations from the general law the stars seem to follow … W. Herschel (1783/1912, vol. 1, p. 120).
**Comment.** He added a reference to Newton's Rule, see above.
**415.** If we indulge a fanciful imagination and build worlds of our own, we must not wonder at our going wide from the path of truth and nature. … If we add observation to observation, without attempting to draw not only certain conclusions, but also conjectural views from them, we offend against the very end for which only observations ought to be made. I will endeavour to keep a proper medium; but if I should deviate from that, I could wish not to fall into the latter error. W. Herschel (1785/1912, vol. 1, p. 223).
**416.** Statistical procedure and experimental design are only two different aspects of the same whole, and that whole is the logical requirements of the complete process of adding to natural knowledge by experimentation. Fisher (1935/1990, p. 3).
**417.** No human mind is capable of grasping in its entirety the meaning of any considerable quantity of numerical data. Fisher (1925/1990, p. 6).
**418.** Diagrams prove nothing, but bring outstanding features readily to the eye; they are therefore no substitute for such critical tests as may be applied to the data, but are valuable in suggesting such tests, and in explaining the conclusions founded upon them. Ibidem, p. 24.

**419.** The null hypothesis is never proved or established, but is possibly disproved, in the course of experimentation. Every experiment may be said to exist only in order to give the facts a chance of disproving the null hypothesis. Ibidem, p. 16

**420.** The knowledge of an average value is a meagre piece of information (p. 35). … It is difficult to understand why statisticians commonly limit their inquiries to Averages, and do not revel in more comprehensive views. Their souls seem as dull to the charm of variety as that of a native of one of our flat English counties, whose retrospect of Switzerland was that, if its mountains could be thrown into its lakes, two nuisances would be got rid of at once (p. 62). Galton (1889, pp. 35, 62).

**421.** Von den mehreren Hunderten thätig Vulkane … ist die Kenntniss noch so überaus unvollständig, daß die einzig entscheidende Methode, die der Mittelzahlen, noch nicht ausgewendet werden kann. Humboldt (1845 – 1862, Bd. 3, p. 288).

**422.** The entire aim of statistics consists in constructing mean values and studying deviations from them. Yanson (1871/1963, p. 264).

**423.** En consultant les Tables de mortalité, je vois qu'on en peut déduire, qu'il n'a que dix mille cent quatre-vingt-neuf à parier contre un, qu'un homme de cinquante-six ans, vivra plus d'un jour. Or comme tout homme de cet âge, où la raison a acquis toute sa maturité & l'expérience toute sa force, n'a néanmoins nulle crainte de la mort dans les vingt-quatre heures, quoiqu'il n'y ait que dix mille cent quatre-vingt-neuf à parier contre un, qu'il ne mourra pas dans ce court intervalle de temps; je conclus, que toute probabilité égale ou plus petite, doit être regardée comme nulle, & que toute crainte ou toute espérance qui se trouve au-dessous de dix mille, ne doit ni nous affecter ni même nous occuper un seul instant le cœur ou la tête. Buffon (1777 /1954, p. 459).

**424.** [Probability $10^{-6}$ is insignificant] on the human scale; [$p = 10^{-15}$ is insignificant] on the terrestrial scale. Borel (1943/1962, p. 27).

**Comment.** Borel did not elaborate. *Human scale* apparently had to do with a single individual.

**425.** Let's remember that the *Titanic* couldn't sink; the probability was not even infinitesimal, it was zero. This fact was widely publicized, yet the ship sank on her maiden voyage. Ekeland (1991/1993, p. 142).

**426.** All scientific conclusions drawn from statistical data require a critical investigation of the basis on which they rest. Newcomb (1902/1906, p. 303).

**427.** I like to make everyone happy when I can. Newcomb, an undated note. The Staatsbibliothek zu Berlin, Manuskriptabteilung, Darmstaedter J 1871 (11); Newcomb.

**Comment.** Darmstaedter was a chemist and a collector of autographs some of which were sent to him in response to his request. The year 1871 possibly indicates the date when Darmstaedter received the note.

**428.** I think that we shall have to get accustomed to the idea that we must not look upon science as a 'body of knowledge' but rather as a system of hypotheses; that is to say, as a system of guesses or anticipations which in principle cannot be justified, but with which we work as long as they stand up to tests, and of which we are never justified in saying that we know they are 'true' or 'more or less certain' or even 'probable'. Popper (1935/1959, p. 317).

**429.** L'imagination, impatiente de remonter aux causes, se plait à créer des hypothèses, et souvent elle dénature les faits pour les plier à son ouvrage; alors les hypothèses sont dangereuses. Mais, quand on ne les envisage que comme des moyens de lier entre eux les phénomènes pour en découvrir les lois, lorsqu'en évitant de leur attribuer de la réalité on les rectifie sans cesse par de nouvelles observations, elles peuvent conduire aux véritables causes, ou du moins nous mettre à portée de conclure des phénomènes observés ceux que des circonstances données doivent faire éclore. Laplace (1814/1886, p. CXLIV).

**430.** Telle est la faiblesse de l'esprit humain, qu'il a souvent besoin de s'aider d'hypothèses, pour lier les faits entre eux. En bornant les hypothèses à cet usage, en évitant de leur attribuer une réalité qu'elles n'ont point, et en les rectifiant sans cesse par le nouvelles observations, on parvient enfin aux véritables causes, ou du moins aux lois des phénomènes. L'histoire de la Philosophie nous offre plus d'un exemple des avantages que peuvent ainsi procurer les hypothèses, et des erreurs auxquelles on s'expose en les réalisant. Laplace (1798 – 1825/1878, t. 3, p. xi).

**431.** On doit considérer que, les lois les plus simples devant toujours être préférées, jusqu'à ce que les observations nous forcent de les abandonner … Laplace (1798 – 1825/1878, t. 1, p. 135).

**432.** On le voit clairement, l'induction, l'analogie, des hypothèses fondées sur les faits et vérifiées, rectifiées sans cesse par de nouvelles observations; tels sont les principaux moyens de parvenir à la vérité. Poisson et al (1835, pp. 176 – 177).

**433.** What is consider'd a pestilence? – If, in a city that can furnish five hundred foot-soldiers, three men [of normal health] die one after the other within three successive days [three men in one day would be deemed a mere accident and not an epidemic], this is accounted a pestilence … *Mishnah*, Taanith $3^4$ with explanations by Blackman, the Editor (No. 2).

**Comment.** Rabinovitch (1973, p. 44) concluded that, wherever there was no clear majority, the existence of a doubt was regarded in the Talmud as 'half and half'. Already here an embryo of the principle of insufficient reason, or indifference, is seen. Still, the probability of the death of an inhabitant of a town during three days and in the absence of an epidemic had hardly been considered equal to 1/2. But there is a similar example, also from the *Mishnah* (Sabbath $6^2$): for being approved, an amulet should have healed three patients consecutively. A probability of 1/8 was apparently neglected.

**434.** And the law which must always be made the foundation of the whole theory is the following: When several hypotheses are presented to our mind, which we believe to be mutually exclusive and exhaustive, but about which we know nothing further, we distribute our belief equally amongst them. Donkin (1851, p. 353).
**Comment.** See No. 435.

**435.** Der Wahrscheinlichkeitsbegriff der klassischen Wahrscheinlichkeitsrechnung, der sich auf das so genannte Prinzip des mangelnden Grundes stützt, kann Niemand mehr, als Grundlage der statistischen Anwendungen der Wahrscheinlichkeitsrechnung, befriedigen. Chuprov (1924, p. 435).
**Comment.** It was Kries (1886, p. 6) rather than Donkin (No. 434) who is thought to have introduced that principle. Keynes (1921/1973, p. 44) renamed it the *principle of indifference*. In essence, Chuprov could have added that the law of large numbers enabled statisticians to apply statistical probability.

**436.** We can't get away from determinism. Chase it out the door by postulating total incoherence, and it comes back through the window, in the guise of statistical laws. Ekeland (1991/1993, p. 50).

**437.** [The statistical method is] as yet little known and is not familiar to our minds; [still,] it is the only available method of studying the properties of real bodies. … [It consists in] estimating the average condition of a group of atoms … [in studying] the probable number of bodies in each group [under investigation]. Maxwell (1871a/1890, p. 253; 1877, p. 242)**.**
**Comment.** Note the sudden similarity of Maxwell's definition with modern ideas (No. 392). Bertrand and Poincaré had also applied the term *probable* instead of *mean*.

**438.** We are compelled to adopt … the statistical method of calculation, and to abandon the strict dynamical method in which we follow every motion by the calculus. Maxwell (1871b, p. 309).

**439.** [The results of the statistical method] belong to a different department of knowledge from the domain of exact science. Maxwell, read 1873; Campbell & Garnett (1882/1884, p. 362).
**Comment.** Maxwell did not elaborate.

**440.** We meet with a new kind of regularity, the regularity of averages, which we can depend upon quite sufficiently for all practical purposes, but which can make no claim to that character of absolute precision which belongs to the laws of abstract dynamics. Maxwell (1873/1890, p. 374).

**441.** Cette [kinetic] théorie mérite-t-elle les efforts que les Anglais y ont consacrés? On peut quelquefois se le demander; je doute que, dès à présent, elle puisse rendre compte de tous les faits connus. Mais il ne s'agit pas de savoir si elle est vraie; ce mot, en ce qui concerne une théorie de ce genre, n'a aucun sens. Il s'agit de savoir si sa fécondité est épuisée ou si elle peut encore aider à faire des découvertes. Poincaré (1894/1954, p. 246).
**Comment.** Poincaré was the father of conventionalism whose essence is well enough illustrated by this very passage.

**442.** Peut-être est-ce la théorie cinétique des gaz qui va prendre du développement et servir de modèle aux autres. Alors les faits qui d'abord nous apparaissaient comme simples ne seraient plus que les résultantes d'un très grand nombre de faites élémentaires que les lois seules du hasard feraient concourir à un même but. La loi physique alors prendrait un aspect entièrement nouveau; ce ne serait plus seulement une équation différentielle, elle prendrait le caractère d'une loi statistique. Poincaré (1905/1970, p. 210).

**443.** Toute loi n'est qu'un énoncé imparfait et provisoire … Il me semble que la théorie cinétique des gaz va nous fournir un exemple frappant. Ibidem, pp. 251 – 252.

**444.** Man möchte die Regelmäßigkeit jener Durchschnittswerte mit der bewunderungswürdigen Konstanz der von der Statistik gelieferten Durchschnittszahlen vergleichen . … Die Bestimmung von Durchschnittswerten ist Aufgabe der Wahrscheinlichkeitsrechnung. Die Probleme der mechanischen Wärmetheorie sind daher Probleme der Wahrscheinlichkeitsrechnung. Boltzmann (1872/1909, p. 316).

**Comment.** Boltzmann likely acquainted himself with Quetelet's writings. On the other hand, he never once mentioned Laplace.

**445.** Die moderne Forschung in ihren höchsten Errungenschaften sich je weiter, in um so grosserem Umfange von stochastisch-statistischen Gesichtspunkten leiten lässt. Chuprov (1922a, pp. 65 – 66).

**446.** In the course of the last sixty or eighty years, statistical methods and the calculus of probability have entered one branch of science after another. Independently, to all appearances, they acquired more or less rapidly a central position in biology, physics, chemistry, meteorology, astronomy, let alone such political sciences as national economy etc. At first, that may have seemed incidental: a new theoretical device had become available and was used wherever it could be helpful, just as the microscope, the electric current … or integral equations. But in the case of statistics it was more than this kind of coincidence.

   On its first appearance the new weapon was mostly accompanied by an excuse: it was only to remedy our shortcoming, our ignorance of details or our inability to cope with vast observational material. … But … the attitude changes. … The individual case is entirely devoid of interest. … The working of the statistical mechanism itself is what we are really interested in. … The first scientific man aware of the vital role of statistics was Darwin. His theory hinges on the law of big [!] numbers. Schrödinger (1944).

**Comment.** The penetration of the statistical method in chemistry is hardly studied. On the other hand, the author forgot medicine. The first who acknowledged the part of statistics in natural sciences had been, to mention only a few scholars, Lamarck, Humboldt and Poisson (Sheynin 2009a, § 10.9 and NNo. 481, 421).

**447.** The statistical method will undoubtedly possess its own merits until science acquires other, more accurate methods. However, statistical inferences cannot be considered on a par with laws of nature, because, for one thing, these inferences have only to do with specific phenomena. Statistics has no power to discover the true causes of the rules which it determines. As any other rules, these are liable to exceptions and [even] radical changes. Beketov (1892, p. 45).

**Comment.** Statistics still *possesses its merits*; it reveals *true causes* when statisticians work together with appropriate specialists and statistical laws of nature do exist.

**448.** Appliquons aux sciences politiques et morales la méthode fondée sur l'observation et sur le calcul, méthode qui nous a si bien servi dans les sciences naturelles. Laplace (1814/1886, p. LXXVIII).
**Comment.** Cf. NNo. 354, 355.

**449.** [Concerning a questionnaire:] Eléments que les sciences naturelles doivent fournir à la statistique pour que celle-ci puisse représenter de la manière la plus complète les diverses manifestations de la vie sociale. *Congrès* (1858, p. 390).
**Comment.** The idea itself was of course timely, but who exactly would have provided the necessary data?

**450.** Le calcul des probabilités a presque toujours été néglige des naturalistes; et cependant des applications nombreuses pourraient en être faites à presque toutes les branches des sciences naturelles (p. 645). … Si nous exceptons la minéralogie, le calcul des probabilités est, comme on l'a vu, la seule branche des sciences mathématiques qui soit applicable aux sciences naturelles (p. 646). Geoffroy Saint-Hilaire (1857, pp. 645, 646).
**Comment.** Mineralogy, or, more precisely, its later isolated branch, petrography, nevertheless does apply statistical methods for appraising the deposits of minerals.

**451.** La condition des sciences médicales, à cet égard [yielding to mathematisation] n'est pas pire, n'est pas autre que le condition de toutes les sciences physiques et naturelles, de la jurisprudence, des sciences morales et politiques, etc. (p. 176). La médecine serait ni une science, ni un art, si elle n'était pas fondée sur de nombreuses observations, et sur le tact et l'expérience du médecin … (p. vi, note). Poisson (1837, pp. 176, vi note).

**452.** Le calcul des probabilités dans les choses morales, telles que les jugements des tribunaux, ou les votes des assemblées, parait à M. Poinsot une fausse application de la science mathématique; il pense qu'on n'en peut tirer aucune conséquence qui puisse servir à perfectionner des décisions des hommes. Suivant M. Laplace lui-même [he quoted the passage here included as No. 101]. Cette idée seule d'un calcul applicable à des choses où se melent les lumières imparfaites, l'ignorance et les passions des hommes, pouvait faire une illusion dangereuse pour quelques esprits, et c'était surtout cette considération qui avait détermine M. Poinsot à prendre un moment la parole sur une question si peu géométrique. Poisson (1836, p. 380; Discussion).

**453.** Je suis loin de prétendre que cette analyse doive suffire à toutes les sciences de raisonnement. Sans doute, dans les sciences qu'on nomme naturelles, la seule méthode qu'on puisse employer avec succès consiste à observer les faits et à soumettre ensuite les observations au calcul. … Cultivons avec ardeur les sciences mathématiques, sans vouloir les étendre au-delà de leur domaine; et n'allons pas nous imaginer qu'on puisse attaquer l'histoire avec les formules, ni donner pour sanction à la morale des théorèmes d'algèbre ou de calcul intégral. Cauchy (1821/1897, p. v).

**Comment.** Cf., however, No. 251.

**454.** Les efforts des géomètres, pour élever le calcul des probabilités au-dessus de ses applications naturelles, n'ont abouti, dans leur partie la plus essentielle et la plus positive, qu'à présenter, relativement à la théorie de la certitude, comme terme d'un long et pénible travail algébrique, quelques propositions presque triviales, dont la justesse est aperçue du premier coup d'œil avec une parfaite évidence par tout homme de bon sens. Comte (1854, p. 120).

**455.** Legionenweise erschienen statistische Angaben in Zahlen und statistische Tabellen voll Zahlen. Lueder (1812, p. 9).

**Comment.** Much earlier similar statements were made concerning natural sciences (Descartes, 1637/1982, p. 63):
*Je remarquais, touchant les expériences, qu'elles sont d'autant plus nécessaires qu'on est plus avancé en connoissance.*

**456.** [After mentioning Marshal Vauban's attempt, at the beginning of the 18$^{th}$ century, to estimate the production of agriculture in France by sampling: this procedure] semble étrange ajourd'hui. Moreau de Jonnès (1847, pp. 53 – 54).

**457.** Enfin, en supposant même que les indications eussent été prises avec le plus grand soin, en dehors de tout esprit de système, et au nombre de trois mille, elles n'auraient pas été en quantité suffisante pour établir un jugement avec le degré de probabilité convenable. Quetelet (1846, p. 8).

**Comment.** Quetelet denied sampling.

**458.** Little experience is sufficient to show that the traditional machinery of statistical processes is wholly unsuited to the needs of practical research. Not only does it take a cannon to shoot a sparrow, but it misses the sparrow! Fisher (1925/1990, pp. viii – ix).

**Comment.** Fisher was here concerned with collection and treatment of data on a large scale.

**459.** A questionnaire is never perfect: some are simply better than others. Deming (1950, p. 31) as quoted by Gaither et al (1996).

**460.** Die Präcision einer statistischer Grösse ist stets als etwas Accessorisches anzusehen Bortkiewicz (1894 – 1896, p. 353).

**Comment.** Until perhaps the 1920s statisticians had not studied estimation of precision (Sheynin 2001a, p. 99).

### 3.8. Application in Natural Science

**461.** Persons who are naturally very fat are apt to die earlier than those who are slender. Hippocrates (1952, No. 44).

**Comment.** This statement illustrates the first stage of the statistical method, when qualitative correlational inferences, which conformed with the qualitative nature of ancient science, were being reached on the basis of a general statistical impression; see Sheynin (2009a, § 1.1.1) for similar examples concerning Aristotle. Also see Sheynin (1984b, § 2.1) about (the not altogether wrong) opinions up to the beginning of the 18$^{th}$ century concerning the influence of the Moon.

**462.** Why do the wicked live, reach old age, and grow mighty in power? … How often is it that the lamp of the wicked is put out? Job 21:7, 21:17.

**463.** Careful men noted what generally answered the better, and then began to prescribe the same for their patients. Thus sprang up the Art of medicine. Celsus (1935, p. 19).

**464.** Aucun de [parts of the body] ne peut changer sans que les autres ne changer aussi … Cuvier (1812/1861, p. 62).

**Comment.** Cuvier discussed (also qualitatively) the *corrélation des formes dans les êtres organisés* and he was apparently one of the first to apply the term itself, *correlation*.

**465.** Co-relation or correlation of structure' is a phrase much used in biology, and not least in that branch of it which refers to heredity, and the idea is even more frequently present than the phrase; but I am not aware of any present attempt to define it clearly, to trace its mode of action in detail, or to show how to measure its degree. Galton (1888, p. 135).

**466.** Its [the correlation theory's] positive side is not significant enough and consists in a simple usage of the method of least squares to discover linear dependences. However, not being satisfied with approximately determining various coefficients, the theory also indicates their probable errors, and enters here the realm of imagination, hypnotism and belief in mathematical formulas that actually have no sound scientific foundation. Markov (1916/1951, p. 533; translation 2004, p. 212).

**Comment.** Revealing dependences (even if only linear) is nevertheless important. In his commentary, Linnik (Markov 1951, p. 670) noted that in those times correlation theory had only been developing so that Markov's criticism made sense. This opinion, however, is hardly true (Hald 1998, § 27.7). In the posthumous edition of his treatise Markov (1924) did not repeat his criticisms, but neither had he considered correlation theory in any detail.

**467.** I ask you to correct and supplement [the review of current literature] without ceremony… our cordial accord is such that I shall sign any of your changes with closed eyes. Chuprov, letter of 22.11.1925, i. e., from a hospital in Rome, to Chetverikov (Sheynin 1990b/1996, p. 57).

**Comment.** There also, on the same page, I quoted Chuprov's letter of 1922 to Chetverikov: "for me, even among my students, no one can replace you". Chetverikov looked after the publication of Chuprov (1926). There, on p. 6, Chuprov expressed his thanks to Chetverikov "without whose friendly assistance and tireless care the Russian edition of this [German] book could not have appeared".

**468.** Excluding biological applications, most of its [correlation theory's] practical usage is based on misunderstanding. Bernstein (1928/1964, p. 231; translation 2005, p. 22).

**469.** To the author the main chain of probability theory lies in the enormous variability of its applications. Few mathematical disciplines have contributed to as wide a spectrum of subjects, a spectrum ranging from number theory to physics, and even fewer have penetrated so decisively the whole of our scientific thinking. Kac (1959, p. ix) as quoted by Gaither et al (1996).

**470.** [The method of gaging *the Heavens, or the Star-Gage*] consists in repeatedly taking the number of stars in ten fields of view … very near each other, and by adding their sums, and cutting off one decimal on the right, a mean of the contents of the heavens, in all the parts which are thus gaged, is obtained. W. Herschel (1784/1912, p. 162).

**Comment.** Herschel believed that the universe was finite and hoped to estimate the number of stars. He later came to realize that his telescope only discovered stars down to a certain magnitude.

**471.** It may be presumed that any star promiscuously chosen … out of such a number [more than 14 thousand], is not likely to differ much from a certain mean size of them all. W. Herschel (1817/1912, p. 579).

**Comment.** With regard to size and physical properties stars extremely differ, they do not belong to a single population at all. The sizes' variance is, practically speaking, infinite and the notion of mean size is meaningless. The theory of probability, like any other scientific discipline, cannot lead to concrete results in the absence of positive knowledge, see NNo. 106, 107.

**472.** Cette inégalité [lunaire] quoique indiquée par les observations, était négligée par le plus grand nombre des astronomes, parce qu'elle ne paraissait pas résulter de la théorie de la pesanteur universelle. Mais, ayant soumis son existence au calcul des probabilités, elle me parut indiqués avec une probabilité si forte, que je crus devoir en rechercher la cause. Laplace (1812/1886, p. 361).

**Comment.** An excellent example of the application of statistical data and probability in natural sciences! Regrettably Laplace did not leave any traces of his stochastic calculations.

**473.** La statistique des astres (s'il est permis de recourir à cette association de mots) doit servir un jour de modèle à toute les autres statistiques. Cournot (1843, § 145).

**Comment.** Actually, stellar statistics had then already been born (Sheynin 1984a, § 6).

**474.** The great Cosmical questions to be answered are not so much what is the precise parallax of this or that particular star, but – What are the average parallaxes of those of the first, second, third and fourth magnitudes respectively, compared with those of lesser magnitude? [And] what connection does there subsist between the parallax of a star and the amount and direction of its proper motion or can it be proved that there is no such connection or relation? Hill & Elkin (1884, p. 191).

**475.** In recent times what we may regard as a new branch of astronomical science is being developed, showing a tendency towards unity of structure throughout the whole domain of the stars. This is what we now call the science of stellar statistics. … In the field of stellar statistics millions of stars are classified as if each taken individually were of no more weight in the scale than a single inhabitant of China in the scale of the sociologist. The statistics of the stars may be said to have commenced with Herschel's gauges of the heavens. … The subject was first opened into an illimitable field of research through a paper presented by Kapteyn to the Amsterdam Academy of Sciences in 1893 (p. 302). The outcome of Kapteyn's conclusions is that we are able to describe the universe as a single object (p. 303). Newcomb (1902/1906, pp. 302, 303).

**476.** Just as the physicist … cannot hope to follow any one particular molecule in its motion, but is still enabled to draw important conclusions as soon as he has determined the mean of the velocities of all the molecules and the frequency of determined deviations of the individual molecules from this mean, so … our main hope will be the determination of means and of frequencies. Kapteyn (1906a, p. 397).

**477.** As in making a contour map, we might take the heights of points at the corners of squares a hundred meters on a side, but we should also take the top of each hill, the bottom of each lake, … and other distinctive points. Kapteyn (1906b, p. 67).

**Comment.** He initiated the study of the starry heaven by stratified sampling.

**478.** Ich schon vor dreißig Jahren angefangen habe, die Unterschiede gleichzeitiger Beobachtungen [made in different localities] an die Stelle der gewöhnliche meteorologischen Constanten zu substituiren, und ich habe jetzt noch die Überzeugung, dass dieser Weg der geeignetste ist, um die Meteorologie als mathematische Disciplin auszubilden. Lamont (1867, p. 245).

**Comment.** Earlier Quetelet (1849, pt. 1, Chapt. 4, p. 53) noted that such differences corresponded to *erreurs accidentelles*. He did not elaborate.

**479.** So groß aber ist … die Herrschaft der Mittel [in meteorology] gewesen, dass Untersuchungen über tägliche und jährliche Veränderungen sich nur dann haben Eingang verschaffen können, wenn sie den Zweck aussprachen, dadurch Methoden anzugeben, durch welche man die mittleren Zustände auf die bequemste und richtigste Weise berechnen könne. Dove (1837, p. 122).

**Comment.** The desire to study the appropriate densities can be seen here. Dove also introduced monthly isotherms.

**480.** Die der Meteorologie gestellte Aufgabe ist dem nach eine dreifache: die Bestimmung der Mittel, die Feststellung der Gesetze der periodischen Veränderungen, und die Angabe der Regeln für die unregelmäßigen. Dove (1839, p. 285).

**481.** Bei allem Beweglichen und Veränderlichen im Raume sind mittlere Zahlenwerthe der letzte Zweck, ja der Ausdruck physischer Gesetze: sie zeigen uns das Stetige in dem Wechsel und in der Flucht der Erscheinungen; so ist, z. B. der Fortschritt der neueren messenden und wägenden Physik vorzugsweise nach Erlangung und Berichtigung der mittleren Werthe gewisser Größen bezeichnet. Humboldt (1845 – 1862, Bd. 1, p. 82).

**Comment.** The study of the laws of density had then only begun. But is it always possible to discuss means *bei allem Beweglichen* etc?

**482.** Tout état de choses dans l'atmosphère … résulte non seulement d'une réunion de causes qui tendent à opérer mais encore de l'influence de l'état de choses qui existait auparavant. Lamarck (1800 – 1811, t. 5, p. 5).

**Comment.** This fact had been suspected even earlier. In 1793 John Dalton applied elementary considerations to study the dependence of weather on previous auroras. In the 19th century, Quetelet (Sheynin 1984b, pp. 78 – 79) described the tendency of the weather to persist by elements of the theory of runs.

**483.** Pour découvrir les lois de la nature, il faut, avant d'examiner les causes des perturbations locales, connaître l'état moyen de l'atmosphère et le type constant de ses variations. Humboldt (1818, p. 190).

**484.** Meine … Linien gleicher Wärme … sind ganz nach Analogie von Halley's isogonischen Curven geformt. Humboldt (1845 – 1862, Bd. 4, p. 59).

**Comment.** In 1701, Halley constructed lines of equal magnetic declination over Northern Atlantic. The introduction of such lines connecting places with equal values of some physical magnitude represented extremely important examples of preliminary treatment of statistical data.

**485.** [A new meteorology is being born.] Little by little [it had begun] to master, synthesize, forecast. Mendeleev (1885/1952, p. 527).

**486.** On a reconnu … que l'abaissement du mercure au-dessous de la moyenne est en général plus grand que son élévation au-dessus de ce terme. Les exemples où la moyenne ne tombe pas à égale distance des deux valeurs limites, et où la courbe de possibilité perd de sa symétrie, sont assez fréquents; ils méritent d'autant plus d'être étudiés, que ce défaut de symétrie tient toujours à des causes plus ou moins curieuses, dont on peut apprécier l'influence. Quetelet (1846, p. 168)**.**

**Comment.** It nevertheless seems that Quetelet largely adhered to usual views even in spite of introducing, in 1848, a mysterious *loi des causes accidentelles* with possibly asymmetric densities (1853, p. 57). Without mentioning him, Liagre (1852, p. 502) did not agree that the density of atmospheric pressure was asymmetric.

**487.** Die Fehlerrechnung ist in der Meteorologie principiell unzulässig. Meyer (1891, p. 32).

**Comment.** Meyer justified his conclusion by the asymmetry of the appropriate densities. Karl Pearson (1898) applied Meyer's data to examine the applicability of his theory of such densities even to the extreme case of antimodal curves.

**488.** Il en est de la géographie des plantes comme de la météorologie; les résultats de ces sciences sont si simples, que de tout temps on a eu des aperçus généraux: mais ce n'est que par des recherches laborieuses et après avoir réuni un grand nombre d'observations précises, que l'on pu parvenir à des résultats numériques et à la connaissance des modifications partielles qu'éprouve la loi de la distribution des [vegetable] formes. Humboldt (1816, p. 228).

**489.** Just as there is political arithmetic, or, in Latino-barbare, statistics, an extremely difficult science, most part of it being moreover conjectural, there is also Arithmetica botanica. Humboldt (1815/1936, p. viii).

**490.** The Calculus of Probabilities … may … serve as a check and a guide to the physical philosopher by pointing out where he may and where he may not employ his study of causes with reasonable hope of a successful result. Spottiswoode (1861, p. 154).

**491.** In short, Statistics reigns and revels in the very heart of Physics. Edgeworth (1913, p. 168/1996, p. 357).

**492.** Hat die Wahrscheinlichkeitsrechnung, nach einer vorübergehenden Periode der Stagnation, infolge des schließlich Sieges der atomistischen Anschauungsweise eine für die Physik ganz grundlegende Bedeutung gewonnen und bildet heute das wichtigste Werkzeug bei Forschungen auf dem Gebiete der modernen Theorien der Materie, der Elektronik, Radioaktivität und Strahlungstheorie. Entspricht doch ihr Wesen durchaus der heute zur Herrschaft gelangten Tendenz, sämtliche Gesetze der Physik* – nach dem Vorbild der kinetischen Gastheorie – auf Statistik verborgener Elementarereignisse zurückzuführen, wobei die *Einfachheit* derselben als sekundäre Folge des Wahrscheinlichkeitsgesetzes *der großen Zahlen* aufgefasst wird.

*Von dieser Tendenz sind bisher nur die Lorentzschen Gleichungen der Elektronentheorie, das Energiegesetz und Relativitätsprinzip unberührt geblieben, aber es ist wohl möglich, dass im Laufe der Zeit auch der exakte Gesetzesformen durch statistische Regelmäßigkeit ersetzt werden dürften. Smoluchowski (1918, p. 253).

**493.** In unserer wissenschaftlichen Erwartung haben wir uns zu Antipoden entwickelt. Du glaubst an den würfelnden Gott und ich an volle Gesetzlichkeit in einer Welt von etwas objektiv Seiendem, das ich auf wild spekulativem Wege zu erhaschen suche. Ich *glaube* fest, aber ich hoffe, dass einer einen mehr realistischen Weg, bezw. eine mehr greifbare Unterlage finden wird, als es mir gegeben ist. Der große anfängliche Erfolg der Quantentheorie kann mich doch nicht zum Glauben an das fundamentale Würfelspiel bringen, wenn ich auch wohl weiß, dass die jüngeren Kollegen dies als Folge der Verkalkung auslegen. Einmal wird's sich ja herausstellen, welche instinktive Haltung die richtige ist. Einstein, letter of 7.9.1944 to Born; Born (1969, p. 204).

**494.** Aber davon bin ich überzeugt, dass man schließlich bei einer Theorie landen wird, deren gesetzmäßig verbundene Dinge nicht Wahrscheinlichkeiten sondern gedachte Tatbestände sind, wie man es bis vor kurzem als selbstverständlich betrachtet hat. Einstein, letter of 3.3.1947 to Born; Ibidem, p. 215.

**495.** Our most precise description of nature *must* be in terms of *probabilities*. There are some people who do not like this way of describing nature. They feel somehow that if they could only tell what is *really* going on with a particle, they could know its speed and position simultaneously. In the early days of the development of quantum mechanics, Einstein was quite worried about this problem. He used to shake his head and say, "But surely God does not throw dice in determining how electrons should go!" He worried about that problem for a long time and he probably never really reconciled himself to the fact that this is the best description of nature that one can give. Feynman (1963/1974, Chapt. 6, p. 15).

**496.** The task of systematizing this science [statistical mechanics], of compiling it into a large book, and of giving it a characteristic name, was executed by one of the greatest American scholars, and in regard to abstract thinking [and] purely theoretic investigation, perhaps to the greatest, Willard Gibbs, the recently deceased professor. Boltzmann (1905, p. 601).

**497.** The statistical ideas … injected [into the kinetic theory] remained however somewhat insulated in a specialized application until Willard Gibbs recognized the fundamental nature of Statistical Mechanics, and blazed the trail for the Quantum Mechanics of the present day. Fisher (1953, p. 4).

**498.** E. B. Wilson, a student and worshipper of Gibbs [in 1953], in a conversation I straightaway wrote down, told me that Gibbs always lectured over the heads of his students and always refused to teach undergraduates at all. … He once told me that in all his years of teaching he had had only six students sufficiently prepared in mathematics and physics to follow him; these included … myself [himself]. Truesdell (1976 – 1981/1984, p. 415).

### 3.9. The Two Currents

**499.** There are periods in the life of each science when the specialists scornfully reject the attempts to analyse theoretically its principles. … Such feelings prevailed in statistics until recently. After the downfall of the Quetelet theoretical system under the blows dealt by German criticisms, years of extreme empiricism ensued in statistics. Chuprov (1909/1959, p. 17).

**500.** Unsere (jüngere) Generation der Statistiker kann sich kaum jenen Sumpf vorstellen, in welchen die statistische Theorie nach dem Zusammenbruch des Queteletschen Systems hineingeraten war und der Ausweg, aus welchem damals nur bei Lexis und Bortkiewicz gefunden werden konnte. Anderson (1932, p. 243/1963, p. 531).

**Comment.** See No. 373.

**501.** Im April 1942 wurde ich an die Universität Kiel berufen. Ich bin ord. Mitglied des Internationalen Statistischen Instituts seit 1937, "Fellow" und "Founder" der internat. "Econometric Society", war Mitglied und Korrespondent für Bulgarien der "International Conference of Agricultural Economists", Mitglied der "Bulgarischen Ökonomischen Gesellschaft", vertrat mein Institut als Mitglied der "Royal Economic Society" (London) usw. Bis 1939 war ich auch assoziiertes Mitglied des "Komitees der statistischen Experten" am Völkerbund in Genf. Ich wurde mehrmals eingeladen, öffentliche Vorträge im Ausland zu halten, so z. B. in Wien (Nationalökonomische Gesellschaft), London (London School of Economics) usw. Selbstverständlich habe ich auch an einer Reihe von wissenschaftlichen Konferenzen und internationalen wissenschaftlichen Expertisen hinzugezogen worden.

1929 – 1942 war ich außerdem Mitherausgeber der Vierteljahresschrift der bulgarischen Generaldirektion der Statistik. Anlässlich meines Ausscheidens aus bulgarischen Diensten verlieh mir der bulgarische König das Kommandeutkreuz des Ordens für zivile Verdienste. Anderson (1946, p. 1).

**502.** His [Lexis'] works, where originality competes with clearness of exposition, opened up the newest epoch in the development of moral statistics. Chuprov (1897, p. 408/2004, p. 23).

**503.** [Lexis had achieved a] Neubegründung einer Theorie der Statistik. Bortkiewicz (1915, p. 119).

**504.** Dormoy published independently and about the same time as Lexis some not dissimilar theories, but subsequent French writers have paid little attention to the work of either. Keynes (1921/1973, p. 431).

**505.** We perpetrate some injustice when discussing the "Lexian theory of dispersion". It would be more correct, more proper … to call it after Dormoy and Lexis. Chuprov (1918/1926). Translated from Russian translation of the Swedish paper (Chuprov 1960, p. 228).

**506.** Ging Dormoy das Verständnis für Anwendungen der Wahrscheinlichkeitsrechnung auf Erfahrungstatsachen in so starkem Maße ab, dass es der historischen Gerechtigkeit nicht entspricht, ihn, so oft von Dispersionstheorie die Rede ist, in eine Reihe mit Lexis zu stellen. Bortkiewicz (1930, p. 53).

**507.** Es ergab sich, daß die bei den untersuchten Reihen gefundenen Schwankungen den Voraussagungen der Theorie fast vollständig entspechen, worin eben das Gesetz der kleinen Zahlen besteht (p. vi). Die Tatsache, dass kleine Ereigniszahlen (bei sehr großen Beobachtungszahlen) einer bestimmten Norm der Schwankungen unterworfen sind bzw. nach einer solchen tendieren, das Gesetz der kleinen Zahlen wohl benannt werden (p. 36). Bortkiewicz (1898, pp. vi, 36).

**Comment.** See comment to No. 508.

**508.** It is difficult to say to what extent the law of small numbers enjoys the recognition of statisticians since it is not known what, strictly speaking, should be called the law of small numbers. Bortkiewicz did not answer my questions formulated on p. [285 of the third] edition of my *Essays* [1909/1959] either in publications or in written form; I did not question him orally at all since he regards criticisms of the law of small numbers very painfully. Chuprov, letter of 1916 to Markov; Sheynin (1990b/1996, p. 68).

**Comment.** Bortkiewicz published his law in 1898, but never formulated it precisely. See my explanation in Bortkevich & Chuprov (2005, pp. 288 – 291). In particular, I repeat that Bortkiewicz distinguished his law, – the "main statistical law" as Romanovsky (1924, No. 4 – 6, p. 15) described it, – from the Poisson distribution. I stress that for almost 60 years, until 1898, that distribution (from which the Bortkiewicz' law, as is now recognized, does not nevertheless differ) had not been applied. I (2008b) consider his law as of no consequence.

**509.** Bortkiewicz had a photographic memory and knew the literature on practically any topic of economics and statistics. … And he had a philosophy of statistics that he had never developed in his writings. For him statistics was not a body of mathematical formulas and techniques but the art of quantitative thinking. An outstanding mathematical statistician, he liked to play with formulas and had published many articles full of algebra, but this was more or less a game; very often the purpose of his mathematical essay was to prove the futility of mathematics. To him the essence was to use measurement to obtain a better understanding of facts of life. … In Germany, he was called the Pope of statistics. … The publishers have stopped asking [him] to review their books [because of his deep and impartial response]. … [He was] probably the best statistician in Europe. … [Bortkiewicz' statement:] I do not review the work of persons I know and I don't care to meet authors whose work I have to appraise. Woytinsky (1961, pp. 451 – 452).

**510.** Für diese [for the Berlin University, Bortkiewicz] beinahe wie ein "Fremdkörper" geblieben ist, welcher sich in den Betrieb und die Traditionen einer deutschen nationalökonomischen Fakultät nicht recht einzuordnen vermochte. Anderson (1932, p. 246/1963, p. 534).

**Comment.** Bortkiewicz worked at the Philosophical Faculty of that University.

**511.** Bortkiewicz schrieb nicht für die weite Öffentlichkeit und war durchaus kein guter Popularisator seiner eigenen Ideen. Er stellte ferner sehr hohe Ansprüche an die Vorbildung und Intelligenz seiner Leser. Mit einer Hartnäckigkeit, die zum Teil durch seine wissenschaftliche Askese bedingt war und zum Teil wohl aus dem "romantischen" Typus seines Forschergeistes entsprang, weigerte er sich, den Rat des "Klassikers" Tschuproff anzunehmen und eine leichtere äußere Form für seine Veröffentlichen zu wählen. Anderson (1932, p. 245/1963, p. 533).

**Comment.** Winkler (1931, p. 1030) quoted a letter from Bortkiewicz (not providing its date) who was glad to find in him one of the five expected readers of his work (which one?).

**512.** You shall have no other gods before me. Exodus 20:3. Quoted by Schumacher (1931, p. 573).

**Comment.** He thus described Bortkiewicz' attitude towards science.

**513.** Zieht er [H. Bruns] zu Gewinnung der in Frage stehenden Resultate die erzeugende Funktion heran, so ist es, als wenn es unternehmen würde, die Gleichung $2x – 3 = 5$ mit Hilfe von Determinante zu lösen. Bortkiewicz (1917, p. iii).

**Comment.** In letter of 1896 No. 7 to Chuprov (Bortkevich & Chuprov 2005), he also mentioned his unwillingness to apply generating functions (and differentiations) when studying the binomial distribution. His opinion was hardly reasonable. And the equation above cannot at all be solved by determinants!

**514.** [Information from the Archive, Humboldt (formerly Friedrich-Wilhelm) University Berlin.] Bortkiewicz was appointed *Ordentlicher Professor* 6.7.1920 (after being 20 years extraordinary professor). Phil. Fak. 1466, Bl. 186).

[An unsigned and undated manuscript, apparently written before that appointment took place, acknowledged his] aufwiegenden Bedeutung als Forscher [but stated that] wir [?] Bedenken … wegen seiner nationalen Zugehörigkeit hegten. [These, however,] wesentlich abgeschwächst worden [after Bortkiewicz] durch eine schriftliche Erklärung sein Bekenntnis zum Deutschtum ablegte. Phil. Fak. 1469, Bl. 63.

**515.** In 1919 [in 1918 – 1919] there appeared a most important and remarkable work of Chuprov devoted to the theory of stability of statistical series (p. 255/168). The revolutionary Essay 3 [part 3] is the most interesting and important (p. 256/171). Romanovsky (1923, pp. 255, 256/2004, pp. 168, 171).

**516.** Chuprov's contributions to science were admired by all … they did much to harmonize the methods of statistical research developed by continental and British workers. *Resolution of Condolence*, Royal Statistical Society (1926).

**517.** In the years before the war [the WWI], it had seemed to me that the continental mathematicians and the Anglo-Saxon statisticians were working without sufficient mutual contact, and that it might be useful to try to join both these lines of research. Cramér (1981, p. 315).

**518.** Chuprov refused to marry, abandoned his homeland, and quit teaching thus giving up his secured existence in order not to divorce himself from science. … [He, Ioffe] knew only one person, Einstein, inspired by his science as much as Chuprov was. Ioffe (1928, p. 315).

**519.** The Laplacian mathematics, although it still holds the field in most text-books, is really obsolete, and ought to be replaced by the very beautiful work which we owe to these three Russians [Chebyshev, Markov, Chuprov] (p. 391). Chuprov … gives by far the best and most lucid general accounts … of the doctrines of the [Continental] school, he alone … [is] writing in a style from which the foreign reader can derive pleasure … (p. 430). The most generalized and, mathematically, by far the most elegant treatment of this problem [concerning different types of frequency] with which I am acquainted, is due to … Chuprov [1918 – 1919, 1919] (p. 380). Keynes (1921/1973, pp. 391, 430, 380).

**520.** To the Most Holy Governing Synod

I have the honour of humbly asking The Most Holy Synod to excommunicate me from the Church. … [In my treatise (1900)] my negative attitude to the legends underlying the Judaic and Christian religions is clearly expressed. … I [also] humbly beg to consider that I do not see any essential difference between ikons and idols … and do not sympathise with religions which, like Orthodoxy, are supported by, and … lend support to fire and sword. Markov, 12 February 1912. Markov Jr (1951, pp. 608 – 609/2004, pp. 247 – 248).

**Comment.** In 1901, Tolstoy was excommunicated from the same Church. In 1910, during his last days, the Synod discussed his case but did not revoke its decision. Markov's request was not granted.

The Synod resolved that he "had seceded from God's church" (Emeliakh 1954, pp. 400 – 401 and 408).

**521.** Markov regarded Pearson, I may say, with contempt. Markov's temper was not better than Pearson's, he could not stand even slightest contradictions either. You can imagine how he took my persistent indications to the considerable scientific importance of Pearson's work. My efforts thus directed were not to no avail as proved by [Markov 1924]. After all, something [Pearsonian] occurred to be included in the field of Markov's scientific interests. Chuprov, letter of 1924; Sheynin (1990b/1996, p. 55).

**522.** Markov … to this day remains an old and hardened sinner in provoking debate. I had understood this long ago, and I believe that the only way to save myself from the trouble of swallowing the provocateur's bait is a refusal to respond to any of his attacks. K. A. Andreev, letter of 1915; Chirikov & Sheynin (1994, p. 132/2004, p. 155).

**Comment.** A similar case (Sheynin 2009a, § 14.3) is Zhukovsky's reproach upon Markov dated 1912.

**523.** Slutsky's book [1912] interests me but does not attract me. I have heard that you have recommended it to students, from which I conclude that you found it worthy of attention. Markov, letter of 1912 to Chuprov; Ondar (1977/1981, p. 53).

**524.** The shortcomings of Pearson's exposition are temporary and of the same kind as the known shortcomings of mathematics in the 17$^{th}$ and 18$^{th}$ centuries. A rigorous basis for the work of the geniuses was built only post factum, and the same will happen with Pearson. I took upon myself to describe what was done [1912]. Sometime Chuprov will set forth the subject of correlation from the philosophical and logical point of view, and describe it as a method of research. An opportunity will present itself to a ripe mathematical mind of a pure mathematician to develop the mathematical basis of the theory. Slutsky, letter of 1912 to Markov; Sheynin (1990b/1996, pp. 45 – 46).

**525.** In Slutsky's person Russian science possesses a serious power especially valuable since in Russia a person occupying himself in the field of social science rarely has mathematical training. Chuprov, no date; Sheynin (1990b/1996, p. 49).

**Comment.** From a reference written in 1916.

**526.** After devoting a year to studying the main works of the English school of mathematical statistics, [Slutsky] published a book [1912] that became a considerable independent contribution; to this day, it remains important and interesting. (p. 68). A refined and witty conversationalist, a connoisseur of literature, a poet and an artist … (p. 72). Kolmogorov (1948a/2002, pp. 68, 72).

**527.** [For a very long time before his death Slutsky remained] almost inaccessible to economists and statisticians outside Russia. … His assistance, or at least personal contacts with him would have been invaluable. Allen (1950, pp. 213 – 214).

**528.** In 1832 [when Darwin obtained the second volume of Lyell's *Principles of Geology*] no book could have been more stimulating for the theories which Darwin later formed. … Lyell determined the order and assigned the successive rock masses the names they now bear by a purely statistical argument. Fisher (1953, p. 2).

**529.** I have no faith in anything short of actual measurement and the Rule of Three. Darwin (1887, vol. 1, p. 411; letter of 1855).

**530.** Favourable variations would tend to be preserved, and unfavourable ones to be destroyed. The result of this would be the formation of new species. Darwin (1859/1958, p. 120).

**531.** The part similar to the main postulate of mechanics, – to the principle of inertia, – is played here by the law that we may call the Darwin law of stationarity. If the existence of some simple trait does not either enhance or lessen the individual's adaptation to life (including fertility and sexual selection), the rate of individuals possessing it persists (in the stochastic sense) from generation to generation. Bernstein (1922/2005, p. 10).

**532.** The simplest reason would tell each individual that he ought to extend his social instincts and sympathies to all the members of the same nation. … This point being once reached, there is only an artificial barrier to prevent his sympathies extending to the men of all nations and races. Darwin (1871/1901, p. 188).

**Comment.** Later history proved him absolutely wrong.

**533.** We looked upon Darwin as our deliverer, the man who had given new meaning to our life and the world we inhabited. K. Pearson (1923, p. 23).

**534.** Some people hate the very name of statistics, but I find them full of beauty and interest. Whenever they are not brutalized, but delicately handled by the higher methods, and are warily interpreted, their power of dealing with complicated phenomena is extraordinary. They are the only tools by which an opening can be cut through the formidable thicket of difficulties that bars the path of those who pursue the Science of man. Galton (1889, pp. 62 – 63).

**535.** [Galton] had written *Hereditary Genius* [1869], one of the most remarkable books of the century. … He became convinced that quantitative, and particularly statistical, methods were needed to consolidate Darwin's ideas, and to give confidence to their practical application. Fisher (1951, p. 34).

**536.** Owing to [my] hereditary bent of mind, … I was well prepared to assimilate the theories of Charles Darwin … [His publications] enlarged the horizon of my ideas. I drew from them the breath of a fuller scientific life, and I owe more of my later scientific impulses to the influence of Charles Darwin than I can easily express.

  The then new doctrine [due to Darwin and Wallace] … burst the enthraldom of the intellect … It gave a sense of freedom … Galton, in 1886 and 1908. K. Pearson (1914 – 1930, vol. 2, p. 201 note).

**537.** The pioneers of the new statistical ideas to a large extent owe their rapid success to the fact that the soil for the propagation of their sermons had been prepared: an authoritative ally, Francis Y. Edgeworth, met those pioneers in the bowels of the Royal Society. For two decades he had been popularizing there the mathematical methods of statistics. His tirelessly penned theoretical studies sparkled with wit and his original and bordering on paradoxes attempts to apply mathematical methods to solve most unusual problems continued to appear persistently. However, those representatives of the statistical science, to whom he had applied, blocked his efforts by a lifeless wall of inertia. His reports were being heard out with respectful attention, his mathematical notes had been published in the Society's periodicals, but his voice had found no response.

Edgeworth is too special in every way: in selecting the subjects of his investigations; in the methods of research; and in the exposition of his papers. Considerable efforts are needed to follow the peculiar train of his ideas and an ingrained habit is required for properly appraising his quaint and aesthetically polished style. Whatever Edgeworth touches, beginning with mathematical psychology and mathematical studies concerning pedagogy (I bear in mind his curious works on the degree of precision of appraising the answers given by those sitting for competitive examinations) and up to theoretical investigations in political economy [economics] and finances, he is able to find such aspects that escape the attention of other researchers; he knows how to do everything in his own way.

But it is this excessive originality that hinders him from becoming widely influential. We may discuss the Pearson school and mention Lexis' students, whereas Edgeworth is a lone figure. Many are those who learn by his writings, but disciples he has none. And we ought to repeat that a single warrior in the field is no warrior [a Russian saying; *Un soldat est une partie de l'armée, et n'est point une armée*; Arnauld & Nicole (1662/1992, p. 295)]. In spite of all his talent, Edgeworth's long-term sermon leaves a decisive trace only after those who independently cultivate the same ideas join him. Nevertheless, the course of subsequent events clearly shows that Edgeworth's activities had not been futile. He was unable to convert his comrades from the Royal Society to his faith, but his authority accustomed them to consider the application of mathematics to statistical investigations as a serious scientific work rather than a hobby pursued by ham dreamers. Under other conditions, the transfer of the new methods developed by statisticians working in biology to the field of the phenomena of social life would have demanded decades, but owing to Edgeworth's previous work this is taking place in only a few years. The triumph of mathematical statistics in England is now a fact come true. Chuprov (1909/1959, pp. 27 – 28).

**538.** Edgeworth's name will stand forever in the history of statistics: I do not mean primarily his work on the Index Numbers … but his work on statistical methods and their foundations that centred in his generalized Law of Error. … Why was this great figure so entirely overshadowed by Marshall? … Edgeworth lacked the force that produces impressive treatises and assembles adherents; amiable and generous, he never asserted himself in any claims of his own; he was oversensitive on the one hand, overmodest on the other; he was content to take a backseat behind Marshall whom he exalted into Achilles; hesitating in conversation, absent-minded to a pathological degree, the worst speaker and lecturer imaginable, he was personally

ineffective – unleaderly is, I think, the word. Schumpeter (1954/1955, p. 831).

**539.** Edgeworth was clearly an original thinker with unusual gifts and an unusual background. But it is hard to see his work in context. Perhaps he had no context; he was a disciple of no school; he founded no school; and Bowley was his only statistical heir. His statistical work is, as a general rule, very tedious to read and hard to follow. Its value was not clear to his contemporaries, including, possibly, Edgeworth himself. Nowadays, I should think, it is not read at all.

Nevertheless he was an important figure in our subject. Apart from the permanent contributions with which his name is associated, he was a great influence on his contemporaries and played a notable part in the development and acceptance of statistics as the subject was then understood. That he did not generate in the social sciences the same surge of interest that Karl Pearson generated in the biological sciences was as much the fault of the soil as of the seed or the sower. Kendall (1968, pp. 274 – 275/1970, pp. 262 – 263).

**540.** The problem of evolution is a problem in statistics. … We must turn to the mathematics of large numbers, to the theory of mass phenomena, to interpret safely our observations. … May we not ask how it came about that the founder of our modern theory of descent made so little appeal to statistics? … The characteristic bent of Darwin's mind led him to establish the theory of descent without mathematical conceptions; even so Faraday's mind worked in the case of electro-magnetism. But as every idea of Faraday allows of mathematical definition, and demands mathematical analysis, … so every idea of Darwin – variation, natural selection … – seems at once to fit itself to mathematical definition and to demand statistical analysis…. The biologist, the mathematician and the statistician have hitherto had widely differentiated field of work. … The day will come … when we shall find mathematicians who are competent biologists, and biologists who are competent mathematicians … Anonymous (1901 – 1902, pp. 3 – 6).

**Comment.** In the mid-19$^{th}$ century, a theory, as opposed to a hypothesis, had to be based on calculations and it seems that Darwin had not therefore originated a *theory*. Faraday was the founder of the science of electromagnetic fields. This statement was published in *Biometrika*; Pearson named the new discipline *biometry*, a term first applied by Chr. Bernoulli (1841, p. 389) in another sense: *Biometrie*, as he called it, represented mass observations concerning *eine ganze Klass oder Gattung* of men.

**541.** I have learned from experience with biologists, craniologists, meteorologists, and medical men (who now occasionally visit the biometricians by night!) that the first introduction of modern statistical methods into an old science by the layman is met with characteristic scorn; but I have lived to see many of them tacitly adopting the very processes they began by condemning. K. Pearson (1907, p. 613).

**542.** [The aims of the Biometric school:] To make statistics a branch of applied mathematics … to extend, discard or justify the meagre processes of the older school of political and social statisticians, and, in general, to convert statistics in this country from being the playing field of dilettanti and controversialists into a serious branch of science. … Inadequate and even erroneous processes in medicine, in anthropology [anthropometry], in craniometry, in psychology, in criminology, in biology, in sociology, had to be criticized… with the aim of providing those sciences with a new and stronger technique. Karl Pearson, written in 1920; E. S. Pearson (1936 – 1937, vol. 29, p. 164).

**543.** Do I therefore call for less human sympathy, for more limited charity, and for sterner treatment of the weak? Not for a moment. Karl Pearson (1909, p. 25) as quoted by Mackenzie (1981, p. 86).

**Comment.** A statement made in connection with his studies of eugenics.

**544.** Shall those who are diseased, shall those who are nighest to the brute, have the power to reproduce their like? Shall the reckless, the idle, be they poor or wealthy, those who follow mere instinct without reason, be the parents of the future generation? Shall the phtisical father not be socially branded when he hands down misery to his offspring, and inefficient citizens to the state? It is difficult to conceive any greater crime against the race. K. Pearson (1887/1901, p. 375).

**545.** "Student" ist nicht ein Fachmann, doch glaub' ich dass glaubt Sie zu vieles Gewicht an seine Wörter liegen. Pearson, letter of ca. 1914 to Anderson; Sheynin (1990b/1996, p. 124).

**Comment.** By that time Student had published six papers, five of them in *Biometrika*!

**546.** [Student was] one of the most original minds in contemporary science. (p. 1). It is doubtful if [he] ever realized the full importance of his contribution to the Theory of Errors (p. 5). Fisher (1939, pp. 1, 5).

**547.** When Karl Pearson and G. Udny Yule began to develop the mathematical theory of correlation in the 1890s, they found that much of the mathematical machinery that Gauss devised for finding "best values" for the parameters of empirical formulas by the method of least squares was immediately applicable in correlation analysis, in spite of the fact that the aims of correlation analysis are the very antithesis of those of the theory of errors. Eisenhart (1978, p. 382).

**548.** Durch die Publikation von … *Grammar of Science* … habe ich einen Forscher kennen gelernt, mit dessen erkenntnisskritischen Ansichten ich mich in allen wesentlichen Punkten in Übereinstimmung befinde, und welcher ausserwissenschaftlichen Tendenzen in der Wissenschaft frei und muthig entgegenzutreten weiss. Mach (1897, Vorwort).

**Commentary.** Mach published these lines in the first after 1897 edition of his book.

**549.** [In 1916 Neyman read Pearson's *Grammar of Science* (1892) on the advice of his teacher at Kharkov University, S. N. Bernstein] and the book greatly impressed us [?]. E. S. Pearson (1936 – 1937, 1936, p. 213).

**550.** It is impossible to write of [Pearson's long article of 1901] without expressing a sense of the extraordinary effort which has gone to its production and of the ingenuity it displays. W. Bateson, letter of 1901 to Pearson (E. S. Pearson 1936 – 1937, vol. 28, p. 230). I respect you as an honest man, and perhaps the ablest and hardest worker I have met. … I have thought for a long time that you are probably the only Englishman I know at this moment whose first thought is to get at the truth in these problems … Bateson, letter of 1902 to Pearson (Ibidem, p. 204).

**551.** [Letter to K. Pearson of 1903.] You are the one living author whose production I nearly always read when I have time and can get at them, and with whom I hold imaginary interviews while I am reading. Newcomb; Sheynin (2002, §7.1).

**552.** Evidently [Pearson] badly takes my remarks about the form of his researches with which I qualify the recognition, – a sufficiently complete recognition, as it seems, – of their scientific value. However, I am one of the most fervent of his apostles among theoreticians of statistics on the Continent. It seems that Pearson is unaware of the extent to which the mathematical forms of his researches hamper an appropriate appraisal of his contributions by scientists who did not go through the English school. Because of Pearson's insufficiently rigorous, to their taste, approaches to mathematical problems Continental mathematicians look down on him to such an extent that they do not even bother to study his works. How many lances did I have occasion to break because of Pearson while substantiating the considerable scientific importance of his oeuvre to the most prominent representatives of the Continental work in related areas! Chuprov, undated letter (ca. 1920); Sheynin (1990b/1996, p. 55).

**553.** Die englische statistische Schule vernachlässigt in ihrem Untersuchungen ein Verfahren, das von russischen und deutschen Gelehrten oft angewandt wird … und geben großer Strenge und Exaktheit noch den Vorzug hat recht elementar zu sein – die Methode der mathematischen Erwartung nämlich. Anderson (1914, p. 269).

**554.** Viel Unheil ist … durch die Abneigung der englischen Forscher gegen die Begriffe der mathematischen Wahrscheinlichkeit und der mathematischen Erwartung gestiftet worden: der Verzicht auf den Gebrauch dieser Grundbegriffe hat der Klarheit der stochastischen Fragestellungen ungemein geschadet, – ja gelegentlich die Lösungsversuche auf Irrwege gelenkt. Legt man jedoch das, dem kontinentalen Auge nicht zusagende Gewand ab und holt, wo dies erforderlich ist, das Versäumte nach, so zeigt such deutlich dass Pearson und Lexis [are akin]. … Nicht *Lexis gegen Pearson*, sondern *Pearson durch Lexis geläutert, Lexis durch Pearson bereichert* sollte gegenwärtig die Parole derer lauten, die, von der geistlosen Empirie der nachqueteletischen Statistik unbefriedigt, sich nach einer rationellen Theorie der Statistik sehnen. Chuprov (1918 – 1919, 1919, pp. 132 – 133).

**555.** While having a look, in connection with my [own] works, at the investigations made by Pearson and his students, I came across quite a number of mistakes – some of them being simply blunders but the rest being rather large methodological mistakes pertaining to numerical solutions. Upon discovering a particular mistake I wrote to Pearson. Unexpectedly the result was spectacular: … Pearson published a long editorial [1919] … where gratitude to me was expressed and the errors were systematically corrected. I was very pleased that the tactic which I had chosen – private information by letter instead of public derision in articles – gave such fruit! Chuprov, letter of 1921; Sheynin (1990b/1996, p. 45).

**556.** Those who have dealt with the sciences have either been empirics or dogmatists. The empirics, in the manner of the ant, only store up and use things; the rationalists, in the manner of spiders, spin webs from their own entrails; but the bee takes the middle path: it collects its material from the flowers of field and garden, but its special gift is to convert and digest it.

The true job of philosophy is not much different, for it depends not only or mainly on the powers of the mind, nor does it take the material gathered from natural history and mechanical experiments and store it unaltered in the memory but lays it up in the intellect changed and elaborated. Therefore from a closer and purer alliance (not so far achieved) of these two faculties (the experimental and the rational) we should have good hopes. Bacon (1620/2004, Book 1, § 95, p. 153).

**557.** From the practical viewpoint the Pearsonian British school is occupying the most considerable place in this field ["theory of probability which comprises the theories of distribution and of the general non-normal correlation"]. Pearson fulfilled an enormous work in managing statistics; he also has great merits, especially since he introduced a large number of new concepts and opened up practically important paths of scientific research. The justification and criticism of his ideas is one of the central problems of current mathematical statistics. Bernstein (1928/1964, p. 228; translation, 2005, p. 19).

**558.** I came in touch with [Pearson] only for a few months, but I have always looked upon him as my master, and myself, as one of his humble disciples. Mahalanobis, letter of 1936; Ghosh (1994, p. 96).

**559.** The modern period in the development of mathematical statistics began with the fundamental works of … K. Pearson, Student, Fisher … Only in the contributions of the British school did the application of probability theory to statistics cease to be a collection of separate isolated problems and became a general theory of statistical testing of stochastic hypotheses … Kolmogorov (1947, p. 63/2005, pp. 79 – 80).

**Comment.** See, however, NNo. 561 and 568.

**560.** The investigations made by Fisher, the founder of the modern British mathematical statistics, were not irreproachable from the standpoint of logic. The ensuing vagueness in his concepts was so considerable, that their just criticism led many scientists (in the Soviet Union, Bernstein) to deny entirely the very direction of his research. Kolmogorov (1947, p. 64/2005, p. 80).

**Comment.** Bernstein (1941/1964, p. 386n; translation 2005, p. 121n) remarked, however, that he did not reject Fisher's work altogether.

**561.** The main weakness of the [Biometric] school when [Slutsky, in 1912] began his research were: 1. Rigorous results on the proximity of empirical sample characteristics to the theoretical ones existed only for independent trials. 2. Notions of the logical structure of the theory of probability, which underlies all the methods of mathematical statistics, remained at the level of eighteenth century results. 3. In spite of the immense work of compiling statistical tables …, in the intermediate cases between 'small' and 'large' samples their auxiliary techniques proved highly imperfect. Kolmogorov (1948a/2002, p. 68).

**562.** Between 1892 and 1911 [Pearson] created his own kingdom of mathematical statistics and biometry in which he reigned supremely, defending its ever expanding frontiers against attacks. Hald (1998, p. 651).

**Comment.** Fisher published only one paper in *Biometrika* (in 1915).

**563.** It is customary to apply the same name, *mean, standard deviation, correlation coefficient*, etc., both to the true value which we should like to know, but can only estimate, and to the particular value at which we happen to arrive by our methods of estimation; so also in applying the term probable error, writers sometimes would appear to suggest that the former quantity, and not merely the latter, is subject to error.

It is this last confusion, in the writer's opinion, more than any other, which has led to the survival to the present day of the fundamental paradox of inverse probability, which like an impenetrable jungle arrests progress towards precision of statistical concepts. Fisher (1922, p. 311).

**Comment.** In the same contribution Fisher (pp. 309 – 310) introduced the notions of consistency, efficiency and sufficiency of statistical estimators, but he is still using the term *true value*. He could have understood it like Fourier did (No. 693). That term is really needed even now, and even in statistics.

**564.** Five or ten years ago I had quite an uphill job to free English opinion from the authoritarian opinion that all enlightenment was to be found in the immaculate gospel of Pearson, and I see no need or profit in undertaking the same thankless task for Thiele. … Thiele and Pearson were quite content to use the same words for what they were estimating and for their estimates of it. … I do *not* mean that his [Thiele's] work is rotten with mathematical errors, only that he had no more glimmer than Pearson of some of the ideas we now use. Fisher, letter of 1931; Fisher (1990b, p. 313).

**565.** [Karl Pearson's] plea of comparability [between the methods of moments and maximum likelihood] is … only an excuse for falsifying the comparison. Fisher (1937, p. 306).

**Comment.** I have not seen any pertinent comment, either by Pearson or anyone else.

**566.** The admiration felt for [Pearson] as a 'giant of the Victorian age' must surely be due to the earlier writings on the philosophy of science, his attempt to develop a consistently rational materialism … I think it would be generally agreed that by far his most important positive contribution to statistical method, the subject on which he wrote most, does lie in the chi-squared goodness of fit test, and that he was singularly unreceptive to and often antagonistic to contemporary advances made by others in this field. The fact that he really achieved, in the eyes of his own circle, his latent ambition of High-priesthood did make his antagonism often effective and deleterious. But for this the work of Edgeworth and of Student, to name only two, would have borne fruit earlier. Fisher, letter of 1946; Edwards (1994, p. 100).

**567.** [Pearson] may be regarded as ploughing the ground in preparation for later developments. That the huge mass of his writings have now little value must be ascribed to two circumstances – first, that his mathematics were on the whole clumsy and lacking in penetration and, second, that without the power of self-criticism he was unable and unwilling to correct his numerous errors or to appreciate the work of others, which would have been of the greatest assistance. He seems to have regarded observational material principally as a means of illustrating his *a priori* concepts, not as a means of correcting them or as providing problems of interpretation in which statistical methods might be of service. He seems to have felt a contempt for the work previously done in the theory of errors and to have known little about it. Fisher (1951, pp. 36 – 37).

**568.** The terrible weakness of his mathematical and scientific work flowed from his incapacity in self-criticism, and his unwillingness to admit the possibility that he had anything to learn from others, even in biology, of which he knew very little. His mathematics, though always vigorous, were usually clumsy, and often misleading. In controversy, to which he was much addicted, he constantly showed himself without a sense of justice. In his dispute with Bateson on the validity of Mendelian inheritance he was the bull to a skilful matador. … His activities have a real place in the history of a greater movement. Fisher (1956/1990, p. 3), about Pearson.

**569.** [Fisher is criticized for his] touchiness and involvement in quarrels [which led to] impediment to communication (p. 446). He never had sufficient respect for the work of [the Neyman – Pearson] school to read it attentively… (p. 448). Savage (1976, pp. 446, 448).

**570.** My new department will … only be slowly organised. I want to give Students of Eugenics … an opportunity to acquaint themselves thoroughly with modern genetical knowledge in animal material. … All the equipment here was very old and bad. The department of Statistics has been separated from the Galton Laboratory … I am afraid I am not an experienced lecturer … Fisher, letter to Romanovsky 5 Febr. 1934 (Sheynin 2008a, pp. 377 – 378).

**571.** I am sometimes accused of intuition as a crime. Fisher, in a conversation with Cramér, 1939; Cramér (1976, p. 529).

**572.** The early statisticians of the present century were competent at mathematics, but they were not great creative mathematicians. Karl Pearson was trained in mathematics, but Edgeworth was a classical scholar and Yule an engineer by training. Fisher, who *was* a creative mathematician, criticized his predecessors for the clumsiness of their style; but even he wrote in the tradition of English mathematics, which does not care much about extreme generalization or extreme rigor as long as it gets the right answer to its problems. The consequence was that, with few exceptions, theoretical statistics in the forties could be understood by anybody with moderate mathematical attainment, say at the first year undergraduate level. I deeply regret to say that the situation has changed so much for the worse that the journals devoted to mathematical statistics are now completely unreadable. Most statisticians deplore the fact, but there is not much they can do about it. Kendall (1972, p. 205) as quoted by Gaither et al (1996).
**Comment.** Kendall found fault with early statisticians, but then everything became even worse! Note that he employed two different terms respectively: theoretical, and mathematical statistics. These terms are indeed not identical: apart from the intuitive feeling correctly representing some real difference, the former, but not the latter includes collection of data and their exploratory (preliminary) analysis.
**573.** [Pearson is] a conscientious and honest enemy of materialism (p. 190). … one of the most consistent and lucid Machians (p. 274). Lenin (1909/1961, pp. 190 and 274).
**Comment.** No wonder that, to the detriment of their scientific work, Soviet statisticians had to consider Pearson almost as an *enemy of the people*.
**574.** [Petersburg (more correctly, Petrograd)] has now for some inscrutable reason been given the name of the man who practically ruined it. K. Pearson (1978, p. 243).
**575.** En octobre 1917, comme au cours de toute sa vie orageuse, Lénine a eu soif du pouvoir pour le pouvoir sans penser ni à la Russie ni au prolétariat russe: il s'intéressait uniquement à une expérience grandiose pratiquée in corpore vili [on a body of small value], sur la people russe. Chuprov (1919, p. 8/2004, p. 95).
**576.** Pearson is a Machian and his curves are based on a fetishism of numbers, their classification is only mathematical. Although he does not want to subdue the real world as ferociously as Gaus [her spelling] had attempted it, his system only rests on a mathematical foundation and the real world cannot be studied on this basis at all. Smit (1934, pp. 227 – 228).
**577.** [Pearson is the author of] some ideas of a racist nature which for five decades forestalled the Göbbels department. Boiarsky & Zyrlin (1947, p. 74).
**Comment.** Contradicting facts was typical of many years of Soviet power.

### 3.10. The Soviet Cul-de-Sac
**578.** The Conjuncture Institute is now firmly established … Our work conditions are now quite acceptable, the researchers' conscience is not violated even to the slightest degree. … At the Institute you will … come upon a staff whose like any other institution in Russia can hardly boast, upon researchers working with all their hearts on [various] problems, at times most agonizing ones. … Minuses: everything here is always "in danger", there is nothing immovable. There are no institutions, no staffs, no plans whose existence are beyond the province of stochastics. But, to repeat, the positions of the … Institute and of Kondratiev are *relatively* very stable. … According to the opinion of an extremely experienced person … "it would be better … to wait at least a year". N. S. Chetverikov, letter to Chuprov, beginning of 1926. Sheynin (1990b/1996, p. 29).

**Comment.** Chetverikov communicated to Chuprov the invitation to join the staff of the Institute from its Director, N. D. Kondratiev. Chuprov would have certainly refused, but, anyway, he was already terminally ill and soon died. And the "extremely experienced person" was absolutely right …

**579.** There is [here] in Paris a friend of mine who was some years ago a lecturer of political economy at the University of Tashkend and now, as an emigrant, allready two years, lives in Paris. He is an able scientist who has published two books. … These books were much praised as original and novel in views and working out of data. They are written … quite not in an orthodox Marxian manner. The name of the author is Alexander [Petrovich] Demidoff. … He was working … on … the present economical state of England, its development and future. Perhaps you and Prof. Hotelling could help him to receive the Rock[e]feller stipend … there is yet a very important point. … Please do not remember [repeat] at all my name, for it can end with my emprisonment by GPU … the most dreadfull and mightfull organisation in the present Russia. It is a crime, and a very heavy one, from its point of view my endeavouring to help an emigrant … Romanovsky; letter of 1929 to Fisher from Paris (Sheynin 2008a, pp. 373 – 374).

**Comment.** OGPU = United Chief Political Administration of Soviet Russia, more commonly remembered by its previous name, GPU (the same, but without the *United*) and the predecessor of the KGB. See also No. 461.

**580.** [An appeal to statisticians] Become the OGPU of scientific thought in statistics and its application to national economy. Smit (1930, p. 168).

**581.** The crowds of arrested saboteurs are full of statisticians. Smit (1931, p. 4), literal translation.

**Comment.** She herself probably informed on the *enemies of the people* still at large and her grammatical merit is obvious. She held a leading position at *Vestnik Statistiki*, the only Soviet statistical periodical (Bortkevich & Chuprov 2005, p. 304), was unbelievably ignorant (see also next entry) but (or perhaps because of that) became Corresponding Member of the Soviet Academy of Sciences. See Sheynin (1998, pp. 534, 536, 539, 540, 543).

**582.** Könnte ich … eine ganze Reihe von in Russland früher sehr geschätzten Statistikern und viel versprechenden jüngeren Schülern … Tschuprows aufzählen, deren Namen nach 1930 aus der sowjet-russischen wissenschaftlichen Literatur plötzlich ganz verschwanden. Anderson (1959, p. 294).

**583.** In justifying its decision to dissolve the [Demographic] Institute [of the Soviet Academy of Sciences], the Academy's Presidium stated [in 1934] that the attempts to introduce social-economic elements [Marxist ideology and conformity?] into the work of the Institute had failed. Tipolt (1972, p. 98).

**584.** Feller and I attended [the Conference of 1937 in Genéva devoted to the theory of probability] from Stockholm, and it was quite exciting to see such a large group of eminent probabilists assembled. Among new acquaintances, there were Steinhaus from Poland, Hopf from Germany, and Jerzy Neyman, at that time still working in England. Several colleagues from the Soviet Union had accepted the invitation and announced lectures, but to our great disappointment none of them turned up. Cramér (1976, p. 528).

**Comment.** All the participants of that Conference are known since they signed a Congratulatory address to Max Born on the occasion of his birthday. It is kept at the Staatsbibliothek zu Berlin, Manuskriptabt., Nachlass Born 129, Gumbel. (Gumbel participated there.) The *disappointment* becomes really *great*.

**585.** It is a fact that Soviet scientists were charged with delivering the leading plenary reports on the theory of probability at the two latest international congresses of mathematicians (in 1932 and 1936). … If even before the revolution the Russian theory of probability, owing to its great weight, might have by right claimed to be a leading force in world science, the Soviet theory of probability has not only totally confirmed this right and justified it even better. … The power necessary for this is certainly drawn exclusively from the inexhaustible source of cheerfulness contained in our new socialist culture. Khinchin (1937, p. 46/2005, p. 54).

**Comment.** Published during the Great Terror! Incidentally, Khinchin could have well been one of those who were forbidden to attend the conference in Genéva, see No 584 (mainly because of the Terror).

   In 1985 or 1986, while having working relations with a few modern mathematicians connected with Kolmogorov's former student Prokhorov, I clearly saw their dislike of the memory of the late Khinchin. And his paper (1961) only appeared posthumously, and then in a philosophical journal.

**586.** The extent or the aftermath of the terror [against statisticians] are not yet studied. Orlov (1990, p. 67).

**587.** Der nunmehr den Staat beherrschende Kommunismus war, seinen Wesen gemäss, die Statistik in hochstem Masse zugeneigt (p. 357). Der Zentralapparat steht jetzt höher, als je. … Aber die statistische Zentralbehörde steht oft vor geradezu unüberwindlichen Schwierigkeiten angesichts der fortschreitenden Desorganisation des Landes. … Sehr schlimm ist mit der Statistik der Bewegung der Bevölkerung bestellt (p. 358). Ebenso trostlos ist der Zustand der Erntestatistik. Die Bauern haben begriffen, dass die statistischen Aufnahmen nicht nur wissenschaftliche Zwecke verfolgen, und verweigern hartnäckig jede Auskunft (p. 359). Chuprov (1922c, pp. 357, 358, 359).

**588.** It is indeed a crime, pure and simple, when contrary to the true situation a mill reports to me that their labour productivity has exceeded the pre-war [1913] level by 12%. And in this respect we are negligent to the utmost, we are entirely unconcerned. Dzerzhinsky (1924, p. 37).
**Comment.** At that time, Dzerzhinsky headed both the OGPU (see No. 602, comment) and the Supreme Council of National economy (VSNKh).
**589.** The editorial staff does not share either the main suppositions of Fisher, who belongs to the Anglo-American empiricists' school, or Romanovsky's attitude to Fisher's constructions. Comment on Romanovsky (1927) on p. 224 of his paper.
**590.** The role of statistics is reduced to measuring the regularities revealed by the specific analysis of the pertinent discipline [of Marxism]. Boiarsky et al (1930, p. 4).
**591.** Market spontaneity is ousted by the activity of organizations directly subordinated to a planned management. … The statistical method begins to retreat in the face of the method of direct accounting. Osinsky (1932, pp. 6 – 7).
**Comment.** He added without explanation that statistics still had a certain role to play. In 1937 Osinsky was arrested and died in 1938. At least during 1930 – 1938 statistics had been all but suppressed (Sheynin 1998, p. 539).
**592.** Under our conditions the method of the theory of probability as applied to mathematical statistics is not sufficient since we have to do with a planned economy and we cannot use the patently unfit means employed by the bourgeois statistics. Yanovskaia (1931)
**Comment.** The Russian term is *theory of probabilities* whereas Yanovskaia, the future celebrated specialist in mathematical logic, used the singular number. In ca. 1963 a much older geodesist told me how, perhaps in 1920, an ignorant *commissar* urged him to measure a baseline in winter by *invar wires*, declaring that the bourgeois restrictions were not valid anymore. At best, this would have been useless. He somehow managed to postpone the measurement, but it was mortally dangerous to refuse flatly.
**593.** Being unable to reveal the regularities in animated nature, they [the geneticists] have resorted to the theory of probability … Physics and chemistry have rid themselves of the accidental. Therefore, they have become exact sciences. … Science is the enemy of the accidental. Lysenko (1948, p. 520).
**Comment.** From ignorant concluding remarks at a conference of 1948 in Moscow. A primitive "researcher" and the so-called *People's academician*, he was indeed squeezed into the Soviet Academy of Sciences.
**594.** Under the impulsion of [Lysenko's] attacks many Russian geneticists, and those among the most distinguished, have been put to death either with or without pre-treatment in a concentration camp. The reward he is so eagerly grasping is Power, power for himself, power to threaten, power to torture, power to kill. Fisher (1948/1974, p. 61).

**595.** Evidently, none of the traditional sciences busies itself about the accidental. Aristotle, *Metaphysics* 1026b.
**Comment.** Aristotle could not have known, but Lysenko, the sham scientist, was obliged to know, that the theory of probability only reveals regularities in *mass* random phenomena.
**596.** For many decades determinism had been preached to students ("science is the enemy of the accidental") (p. 67). We cannot reconcile ourselves anymore to the situation in statistics brought about by the legacy of Stalinism. The split in statistics, and the lack of necessary knowledge typical of many "specialists" are leading to an ever increasing lag behind the advanced nations with respect to the mass application of modern statistical methods (p. 65). Orlov (1990, pp. 67, 65).
**597.** Throughout the discussion [before and during the Conference, see No. 593], Lysenko and his followers treat neo-Mendelism (or Morgano-Mendelism or whatever other title they apply to modern genetics) as a mere theory, in the sense of a hypothesis, not in the usual sense in which the term *theory* is used in science, of a set of conceptions tying together a vast body of experimental results and established laws; what is more, they treat it as a theory inspired primarily by philosophical and political principles hostile to Marxism and Communism, and not by the desire to find the simplest explanation of the facts of nature (p. 83). Michurinites are, without knowing it, simply utilizing the basic Mendelian principles of segregation and recombination, followed by the principles of Darwinian selection, but undiluted by any trace of Lamarckism (p. 93). J. Huxley (1949, pp. 83, 93).
**Comment.** An honest hard-working practitioner and the Soviet counterpart of the American horticulturist Burbank, Ivan Vladimirovich Michurin (1855 – 1935). Lysenko and his followers unfairly contrasted him with geneticists.
**598.** Some hotheads, without gaining an understanding of what did Academician [!] Lysenko really say, and being ignorant of the theory of probability, decided to declare a campaign against it. Gnedenko (1950b, p. 8).
**599.** The chromosome theory of heredity has become part of the gold fund of human knowledge. … I am in a position to verify this theory from the point of view … of statistics. And it also conforms to my ideas. Nemchinov (1948).
**600.** The Conference resolutely condemns the speech of V. S. Nemchinov … for his attempt statistically to "justify" the reactionary Weismann theories. Objectively [he] adhered to the Machian Anglo-American school which endows statistics with the role of arbiter situated over the other sciences, a role for which it is not suited (p. 313/181 – 182). The Conference accepts with satisfaction the statement of … Romanovsky, who confessed to having made ideological mistakes in some of his early works (p. 314/182). The methods of bourgeois statistics were not always critically interpreted; sometimes they had been propagandized and applied (p. 314/182).

[One of the goals:] A construction of a consistent system of mathematical statistics embracing all of its newest ramifications and based on the principles of the Marxist dialectical method (p. 315/183). *Resolutsia* (1948, pp. 313, 314, 315/2005, pp. 181 – 182, 182, 183).

**Comment.** From the Resolution of the Second All-Union Conference on Mathematical Statistics. Much more can be cited. Only the Party had been the arbiter. On Nemchinov see the English edition of the *Great Sov. Enc.*, vol. 17, 1978.

**601.** [Nemchinov made an early attempt to introduce econometrics in the Sovet Union. This is the official comment: He endeavoured to] transfer the arch-bourgeois mathematical school of economics to the Soviet soil …Anonymous (1948, p. 82).

**602.** The fabulous statistics continued to pour out of the telescreen. As compared with last year, there was more food, more clothes, more houses, more furniture, more cooking pots, more fuel, more ships, more helicopters, more books, more babies – more of everything except disease, crime, and insanity. Orwell (1949/1969, p. 59).

**603.** There is practically no statistics in Russia and it is a surprising feature that a country which is so strong in probability theory has made practically no contribution to mathematical statistics. Obviously the political atmosphere is very unfavourable to that type of application. W. Feller, letter of 1948 or 1949 to J. Huxley; J. Huxley (1949, p. 170).

**Comment.** Feller likely had in mind the application of statistics to various sides of social life, and indeed there are no pertinent examples in the generally known source, Gnedenko (1950a). As to mathematical statistics, Feller was wrong.

**604.** The Moscow school [of probability] may … be reproached for not sufficiently or not altogether systematically studying the issues of practical statistics. The main problem next in turn … consists in revising and systematizing the statistical methods that still include many primitive and archaic elements. This subject is permanently included in the plans of the Moscow school and just as invariably postponed. Khinchin (1937, p. 43/2005, pp. 50 – 51).

**605.** Being a theoretician of modern bourgeois statistics, bourgeois narrow-mindedness and formal viewpoints are in his nature. … In actual fact, he completely ignores the qualitative aspect of phenomena. Fisher (1925/1958; Publisher's Preface).

**Comment.** A Russian translation of his book should have appeared 30 years earlier. The qualitative aspect once again!

**606.** [Statistics] is not subordinated to other sciences. M. V. Ptukha; Anonymous (1954, p. 44).

**607.** [Statistics] is an independent science. S. G. Strumilin; Ibidem, p. 41.

**608.** Only the revolutionary Marxist theory had been a robust foundation for the development of statistics as a social science. A. M. Vostrikova; Ibidem, p. 41.

**609.** Lenin had completely subordinated [adapted] the statistical methods of research … to the problem of the class analysis of the rural population. [It is impossible to maintain that] the same methods of research are applied [in astronomy and economics]. K. V. Ostrovitianov; Ibidem, p. 82.

**Comment.** I have quoted four participants of the Conference on the problems of statistics held in 1954 in Moscow. Ostrovitianov's opinion meant that he, the Vice-President of the Soviet Academy of Sciences, had no understanding of statistics. Even more strange were the utterances of V. A. Sobol (p. 61) and S. P. Partigul (p. 74). The former stated that statistics did not study mass random phenomena, and the latter denied any regularities of such phenomena.

**610.** Statistics is an independent social science. It studies the quantitative aspect of mass social phenomena in an indissoluble connection with their qualitative side. Strumilin (1969, p. 87).

**Comment.** This was the definition of statistics adopted at the Conference. The *qualitative aspect* (= the Party directive of the day) had invariably been stressed up to the collapse of the Soviet Union. I am not however denying it, cf. No. 407, provided the statistician works together with the appropriate specialist rather than relies on directives.

**611.** Anyone who does not examine the meaning of the numbers with which he performs some calculations is not a mathematician. Buniakovsky (1866, p. 154).

**Comment.** Mathematicians were not less clever than statisticians! Cf. No. 407.

**612.** The statistician cannot evade the responsibility for understanding the process he applies or recommends. Fisher (1935/1990, p. 1).

**613.** Statisticians are barely acquainted with the theory of probability whereas mathematicians, who applied mathematical calculations for analysing numerical information concerning social phenomena, considered it as abstract quantities, did not allow for special properties of these phenomena and had therefore been arriving at nearly absurd inferences. Yanson (1887/1963, p. 291).

**Comment.** This only means that mathematicians should work together with statisticians.

**614.** Plebiscience is an intermediate stage [in history of science]. The next and last is prolescience. … Its function is to confirm and comfort the proletariat in all that will by then have been ordered to believe. Of course, that will be mainly social science. Truesdell (read 1979, 1981/1984, p. 117).

**Comment.** Truesdell did not study Soviet statistics …

**615.** [The aims of statistics] Collecting reports for the state machinery and one-sided propaganda … [Work under conditions of] secrecy, colouring the truth. Methodological niceties became useless. Yasin (1989, p. 774).

**616.** Only during the perestroika [Gorbachev's doomed attempt to reform the Communist regime] the veil of secrecy began to open slightly. And methods of falsifying statistical data, which made it possible to create a semblance of well-being, were at once revealed (p. 65). We reject the decisions of the All-Union conference of 1954 as impeding the perestroika. The mistaken attribution of statistics to the social sciences considerably delayed the development of national economy. A barrier had been erected between modern theoretical (mathematical) statistics and the agencies [of government statistics] whose activities were almost reduced to registration (p. 67). Vast trustworthy statistics was not needed, [it was] even dangerous for the Soviet system (Ibidem). Orlov (1990, pp. 65, 67).
**Comment.** Neither the national economy nor the structure of the state could have been reformed within the bounds of the existing system of power whose invariability was in the vital interests of the ruling clique at large. The perestroika failed.

### 3.11. Econometry

**617.** Nor need we doubt that the study of Statistics will, ere long, rescue Political Economy from all the uncertainty in which it is now developed. Portlock (1839).

**618.** But still higher tribute [than to Jevons] is due to Cournot who, without encouragement or lead, in what was then a most uncongenial environment, in 1838 fully anticipated the econometric program by his *Recherches*, one of the most striking achievements of true genius, to which we pay respect to this day by nearly always starting out from them. Schumpeter (1933, p. 8).

**619.** Monsieur le Président,
Le calcul des probabilités a été l'objet de ses [Laurent's] constantes préoccupations; ses travaux dans cette science lui ont valu une grande autorité à l'Institut des Actuaires Français, institution très vivante et très utile. Il a résolu habilement diverses questions de hautes Mathématiques relatives à la Théorie des Assurances sur la vie [1894]. C'est aussi au Calcul des Probabilités que se rattache un mémoire sur la méthode des moindres carrés.

   Enfin un traité d'économie politique mathématique [1902]. Cette science nouvelle crée par Walras et ses disciples, et de nature à préciser des conceptions des économistes et si l'on peut craindre que la précision du langage mathématique ne convienne pas toujours à un objet quelquefois un plus vague, les idées nouvelles ont certainement rendu des services en faisant justice de certains raisonnements qui n'avaient que l'apparence de la rigueur alors qu'ils ne pouvaient même servir de première approximation. Laurent dans ses dernières années s'est beaucoup occupé de cette science nouvelle qu'il a professée à la Sorbonne dans un cours libre. Poincaré; undated and unpublished letter kept in his Dossier at the Paris Academy of Sciences.

**620.** The Econometric Society is an international society for the advancement of economic theory in its relation to statistics and mathematics. … Its main object shall be to promote studies that aim at a unification of the theoretical-quantitative and the empirical-quantitative approach to economic problems and that are penetrated by constructive and rigorous thinking similar to that which has come to dominate the natural sciences. Frisch (1933, p. 1).
**Comment.** From the signed Editorial in the first issue of the *Econometrica*. Also see NNo. 220 and 601.

**621.** The main, difficult but necessary aim is to express the desired optimal state of affairs in the national economy by a single indicator. A. N. Kolmogorov, statement at a conference of 1960; Birman (1960, p. 44).
**Comment.** This was an indirect criticism of the law of value.
**622.** In the 42$^{nd}$ year of the existence of the socialist state our economic science does not know precisely what the law of value means in a socialist society or how it should be applied. It does not know what socialist rent is, or whether in general there ought to be some calculations of the effectiveness of capital investment. … We are offered as the latest discovery in the field of economics, for example, the proposition that the law of value does not govern, but only influences … Kantorovich (1959, p. 60).
**Comment.** No one dared study the quantitative side of economics. See also R. W. Campbell (1961).

## 4. Mathematical Treatment of Observations

### 4.1. Ancient Astronomers

**623.** Observations [in antiquity] were more qualitative than quantitative; "when angles are equal" may be decided fairly well on an instrument but not "how large are the angles", says Ptolemy … with respect to the lunar and solar diameters. Neugebauer (1948/1983, p. 101).

**624.** The characteristic type of measurement [in antique astronomy] depended not on instrumental perfection, but on the correct choice of crucial phenomena. Aaboe & De Solla Price (1964, p. 3).

**625.** [Planetary observations that] are most likely to be reliable are those in which there is observed actual contact or very close approach to a star or the moon and especially those made by means of the astrolabe instruments. Ptolemy (1984, IX 2, p. 423; H 213).

**626.** We have hardly anything from Ptolemy that we could not with good reason call into question prior to its becoming of use to us in arriving at the requisite degree of accuracy. Kepler (1609/1992, p. 642).

**627.** Ptolémée, à qui nous devons principalement la connaissance de ses [Hipparchus'] travaux et qui s'appuie sans cesse sur ses observations et sur ses theories, le qualifie avec justice "d'astronome d'une grande adresse, d'une sagacité rare, et sincere ami de la vérite". … Ses Tables du Soleil, malgré leur imperfection, sont un monument durable de son génie, que Ptolémée respecta au point d'y assujettir ses propres observations (p. 413). Même qu'il est entièrement détruit, l'*Almageste*, considéré comme le dépôt des anciennes observations, est un des plus précieux monuments de l'antiquité (p. 422). Laplace (1796/1884, pp. 413, 422).

**628.** All of Ptolemy's *Almagest* [his main work] seems to me to breathe an air of perfect sincerity. Newcomb (1878, p. 20).

**629.** Ptolemy, like astronomers of today, undoubtedly built his edifice on a great array of traditional materials, rejecting, adjusting, or incorporating them as he saw fit, and molding them into a new theoretical framework … Gingerich (1983, p. 151).

**630.** [Ptolemy was an opportunist ready] to simplify and to fudge. Wilson (1984, p. 43).

**631.** [Ptolemy's system] held mankind in spiritual bondage for fourteen centuries. Chebotarev (1958, p. 579).

**Comment.** The author was a worthy Soviet lackey, and had nothing in common with astronomy!

**632.** One of the most admirable features of ancient astronomy [was] that all efforts were concentrated upon reducing to a minimum the influence of the inaccuracy of individual observations with crude instruments by developing to the farthest possible limits the mathematical consequences of very few elements [of optimal circumstances] (p. 250).

The 'doctoring' of numbers for the sake of easier computation is evident in innumerable examples of Greek and Babylonian astronomy. Rounding-off in partial results as well as in important parameters can be observed frequently, often depriving us of any hope of reconstructing the original data accurately (p. 252). Neugebauer (1950/1983, pp. 250, 252).

**633.** In all ancient astronomy, direct measurements and theoretical considerations are … inextricably intertwined. … Ever present numerical inaccuracies and arbitrary roundings … repeatedly have the same order of magnitude as the effects under consideration. Neugebauer (1975/2004, p. 107).

**634.** In former centuries the astronomer selected from among his observations those that seemed the best … [An example concerning Tycho Brahe follows]. In the 17th century scientists like Huygens and Picard realized that the average of a number of equivalent measurements would be better than one of a couple of selected from them, and in the 18th century this averaging came more and more into use … A new attitude was brought into being, typical of the 19th century scientist towards his material: it was no longer a mass of data from which he selected what he wanted, but it was a protocol of an examination of nature … Pannekoek (1961, pp. 339 – 340).

**Comment.** This conclusion is only partly correct. On ancient astronomy see above; Huygens is only mentioned, and no pertinent particulars seem to be known; and already Kepler followed "the new attitude".

**635.** [Ptolemy is] the most successful fraud in the history of science. R. R. Newton (1977, p. 379).

**Comment.** An absolutely wrong statement, see previous items.

**636.** The scientific literature of the seventeenth century – and not only of the seventeenth century – is full of these fictitious experiments. Koyré (1956, p. 150).

**Comment.** This is a very important testimony since many commentators even before R. R. Newton (No. 635) have expressed themselves in the opposite sense. Ptolemy had indeed, without explaining his decisions, treated observations as he saw fit. This, however, was in the spirit of the time as was his likely borrowing without any acknowledgement the Hipparchus observations. However, the history of astronomy knows some examples of deceit committed by classics. Thus, Donahue (Kepler 1609/1992, p. 3, note 7) noted that "the entire table at the end of Chapter 53 [of that source], for example [!] is based upon computed longitudes presented as observations". We ought to realize that *Quod licet Jovi, non licet bovi*!

**637.** Throughout the ages down to [up to] at least AD 1600 the roundings off of fractions is usually done at the computer's pleasure… Hartner (1977, p. 3).

**638.** In all probability I-Hsing [one of the astronomers] thought it very undesirable to admit … a mass of raw data showing considerable scatter, and not being able to assess it statistically, he used it only to satisfy himself that his calculated values came about were they should – indeed, he probably believed that they were much more reliable than most of the observations. Needham (1962, p. 51).

**Comment.** An opinion about the meridian arc measurement in China in the 8[th] century.

**639.** Observers of an eclipse should obtain all its times [phases] so that every one of these, in one of the two towns, can be related to the corresponding time in the other. Also, from every pair of opposite times, that of the middle of the eclipse must be obtained (p. 155). The use of sines engenders errors which become appreciable if they are added to errors caused by the use of small instruments and errors made by human observers (p. 237). Al-Biruni (1967, pp. 155, 237).

**Comment.** Observations of lunar eclipses were being used for determining the longitudinal difference of two places. The described procedure allowed to estimate qualitatively the possible discrepancies. The reasoning about errors was possibly the first pertinent recorded statement.

**640.** Many medieval maps may well have been made from general knowledge of the countryside without any sort of measurement or estimation of the land by the 'surveyor'. D. J. Price (1955, p. 6).

### 4.2. Seventeenth and Eighteenth Centuries

**641.** Since the divine benevolence has vouchsafed us Tycho Brahe, a most diligent observer, from whose observations the $8'$ error in this Ptolemaic computation is shown, it is fitting that we with thankful mind both acknowledge and honour this benefit of God. … Now, because they could not have been ignored, these eight minutes alone will have led the way to the reformation of all of astronomy, and have constituted the material for a great part of the present work. Kepler (1609/1992, p. 286).

**Comment.** I believe that Kepler applied the elements of the minimax method and that the eight minutes became an impossible least value of the maximal deviation of the Ptolemaic theory from the Tychonian observations.

Another point that Kepler (Ibidem, p. 60ff) made after establishing the *New Astronomy* was to explain the known and apparently opposing statement in the Old Testament.

**642.** One might, however, hold suspect such license in making small changes in the data, thinking that by taking the same liberty in changing whatever we don't like in the observations, the full Tychonic eccentricity might also at last be obtained. Anyone who thinks this way should make a try at it, and, comparing his changes with ours, he should judge whether the changes remain within the limits of observational precision. He also needs to beware lest, elated by the results of one such iteration he give himself much more hideous problems in what follows, because of the very divergent apogees found for the sun (p. 334). If this wearisome method has filled you with loathing, it should more properly fill you with compassion for me, as I have gone through it at least seventy times at the expense of a great deal of time (p. 256). Ibidem, pp. 334, 256.

**Comment.** Kepler apparently applied the elements of statistical simulation ("let another quantity be taken at will", p. 256) in which case a successful result was only possible had he allowed for the properties of *usual* random observational errors.

**643.** So cörperlich und so greifflich gehet es nicht zu, dass Himmel und Erde einander anrühreten, wie die Räder in einer Uhr. Kepler (1610/1941, § 57, p. 200).

**Comment.** This is Kepler the astrologer explaining the qualitative correlation between heaven and Earth. See No. 761.

**644.** Newton in fact (but not in explicit statement) had a precise understanding of the difference between random and structurally 'inbuilt' errors. He was certainly, himself, absorbed by the second type of 'inbuilt' error, and many theoretical models of differing types of physical, optical and astronomical phenomena were all consciously contrived so that these structural errors should be minimized. At the same time, he did, in his astronomical practice, also make suitable adjustment for 'random' errors in observation. D. T. Whiteside, private communication of 1972.

**645.** Experiments ought to be estimated by their value, not their number; … a single experiment … may as well deserve an entire treatise … As one of those large and oriental pearls may outvalue a very great number of those little … pearls, that are to be bought by the ounce … Boyle (1772/1999, p. 376).

**646.** He [Flamsteed] does not appear to have taken the mean of several observations for a more correct result … Where more than one observation of a star has been reduced, he has generally assumed that result which seemed to him most satisfactory at the time, without any regard to the rest. Neither, in fact, did he reduce the whole (nor anything like the whole) of his observations: many day's work having been wholly omitted in his computation-book. And, moreover, many of the results, which have been actually computed … have not been inserted in any of his MS catalogues. Baily (1835, p. 376).

**Comment.** Several considerations ought to be taken into account. Thus, it is known that Flamsteed never hurried to publish his observations. Then, Thoren (1972) formulated a very high opinion of Flamsteed as an observer skilful at "manipulating" data.

**647.** This [discovery of nutation] points out to us the great advantage of cultivating [astronomy] as every other branch of natural knowledge, by a regular series of observations and experiments (p. 1).

   Science … had acquired such extraordinary advancement, that future ages seemed to have little room for making any great improvements. But, in fact, we find the case to be very different; for, as we advance in the means of making more nice inquiries, new points generally offer themselves, that demand our attention (p. 2).

   When several observations have been taken of the same star within a few days of each other, I have either set down the mean result, or that observation which best agrees with it (p. 24). Bradley (1748, pp. 1, 2, 24).

**648.** Let *p* be the place of some object defined by observation, *q, r, s* the places of the same object from subsequent observations. Let there also be weights *P, Q, R, S* reciprocally proportional to the displacements arising from the errors in the single observations, and which are given by the limits of the given errors; and the weights *P, Q, R, S* are conceived as being placed at *p, q, r, s*, and their centre of gravity *Z* is found: I say the point *Z* is the most probable place of the object. Cotes, 1722; Gowing (1983, p. 107).

**Comment.** Cotes appended a figure (perhaps representing a three-dimensional picture) showing nothing except these four points. He did not explain what he understood as *the most probable place*. If the limits of the errors are proportional to the mean square errors (unknown to Cotes), the weights will become inversely proportional to the latter rather than to their squares as indicated by Gauss (1823b/1887, § 7).

The centre of gravity corresponds to the arithmetic mean which beginning at least during Kepler's lifetime had become a universal estimator of the measured constant. This mean occurs in approximate calculations of areas of figures and volumes of bodies so as to compensate the errors of approximate formulas and/or the deviations of the real figures and bodies from their accepted models (Colebrooke 1817, p. 97).

**649.** La règle de Cotes fut suivie par tous les calculateurs. Laplace (1814/1886, p. CL).

Les meilleurs astronomes ont suivi cette méthode, et le success des tables qu'ils ont construites à son moyen en a constaté l'avantage. Laplace (1812/1886, pp. 352 – 353).

**Comment.** This preceding statement also discusses the Cotes method, and here Laplace correctly bears in mind the period after 1750. But what exactly did he understand as that method?

**650.** Tob. Mayer nicht nach einem systematischen Princip, sondern nur nach hausbackenen Combinationen gerechnet hat. Gauss, letter of 24 June 1850 to Schumacher. Gauss (1863/*Werke, Ergänzungsreihe*, Bd. 5/3, 1975, p. 90).

**Comment.** In 1750, Mayer considered 27 equations in three unknowns. He separated them in three reasonably selected disjoint groups of nine equations each and solved these by "combinations". Gauss was discussing Mayer's manuscript, but the method of solution was likely the same. Then, Gauss himself, in an earlier letter of the same year (22 February) published in the same sources (pp. 66 – 67), discussed his own similar procedure. True, he applied it for calibrating an aneroid rather than for treating astronomical observations as Mayer did.

**651.** Euler's work [1749] was, in comparison with Mayer's [1750], … a statistical failure (p. 27). He distrusted the combination of equations, taking the mathematician's view that errors actually increase with aggregation rather than taking the statistician's view that random errors tend to cancel one another (p. 18). Stigler (1986, pp. 27, 28).

**Comment.** Euler applied the elements of the minimax method which is the best possible for checking whether the received theory (about which he had serious doubts) conformed to the observations. Then, not only (pure) mathematicians, whom Stigler had not named, but even Laplace and Legendre really feared an accumulation of errors, see No. 653. In his later book Stigler (1999, p. 318), without mentioning his previous opinion, stated that [in 1778] Euler, by denying the principle of maximum likelihood, *was acting in the grand tradition of mathematical statistics*.

**652.** If we assume that among the observations … there is one that should be almost rejected, whose degree of goodness would accordingly be as small as possible, it is evident that the product [the likelihood function] would in fact be reduced to nothing, so that it could not possibly be considered as a maximum, no matter how great it might be, were that observation omitted. … I do not think that it is necessary … to have recourse to the principle of the maximum, since the undoubted precepts of the art of conjecturing are quite sufficient to resolve all questions of this kind. Euler (1778/1970, p. 168).

**Comment.** The addition of an unworthy observation should not change the result, Euler also reasoned, whereas, when applying *the principle of the maximum*, it leads to a significant change. Euler proposed an estimator practically coinciding with the arithmetic mean and heuristically resembling the Gauss principle of maximum weight (least variance), although his reasoning should have resulted in the choice of the median. In the passage above, I replaced the translator's modernized term, *theory of probability*, by Euler's original expression.

**653.** Il se présentoit trois parties [for calculating the triangulation laid out between two baselines]: le premier, d'employer une des bases à l'exclusion de l'autre, le second, de répartir la différence dont nous venons de parler entre les deux bases proportionnellement, si l'on veut, à leur longueur; le troisième, d'employer la [northern] base dans le calcul de la partie boréale de la méridienne, et [the southern] dans le calcul de la partie australe. C'est le parti auquel nous avons cru devoir nous arrêter … Le second parti est beaucoup plus plausible, et paroît au premier abord fort naturel [but will lead to large corrections]; suivant le troisième parti on laisse les bases intactes, et les côtés de chaque triangle ne sont sujets à d'autre incertitude qu'à celle qui peut naître des petites erreurs qui restent dans la détermination des angles. Enfin, en employant la [northern] base depuis … jusqu'à …, et [the southern baseline] depuis … jusqu'à …, on partage l'étendue de la méridienne en deux parties à peu près égales, on prévient l'accumulation ultérieure de petites erreurs, on s'en tient de plus près à ce que les observations même [of the angles] donnent immédiatement. Van Swinden et al (an VII, pp. 421 – 422).

**Comment.** An VII ≈ 1799. Two baselines measured at the ends of a triangulation chain allow to check its linear scale but lead to a discrepancy (the base condition). In those times, the errors of linear measurements had a much more considerable influence than those of angle measurements so that the corrections should have been mainly applied to the baselines. Laplace (ca. 1819/1886, pp. 590 – 591) later justified the then adopted decision by the lack of a sound method of adjustment. Bru (1988, pp. 225 – 228) described this episode and noted that other scholars (Maupertuis, Bouguer, Condamine) had also experienced similar difficulties.

**654.** Auch die übliche Methode einer Ausgleichung in einem Guss, bei der alle im Netz auftretenden Bedingungsgleichungen nach der [method of least squares], kann in einem so dichten astronomisch-geodätischen Netz infolge des ungeheuren Arbeitsumfanges nicht angewendet werden. Außerdem ist, wie Prof. F. N. Krasowski … zum Ausdruck bringt, zu bezweifeln, ob die Ergebnisse einer solchen Ausgleichung die besten sein werden (p. 438). Diesen Weg für die Ausgleichung des astronomisch-geodätischen Netzes der [Soviet Union] hat Prof. F. N. Krasowski vorgeschlagen, der damit die Methode von Helmert modifizierte, die zu kompliziert und umständlich war, um ein so dichtes … Netz wie das sowjetische ausgleichen zu können (p. 440). Sakatow (1950/1957, pp. 438, 440).

**Comment.** At least the idea of temporarily replacing triangulation chains by the appropriate geodesics is due to Helmert. The system of geodesics can then be subjected to adjustment with a subsequent return to the chains (then to be adjusted separately and finally). Concerning the advisability of managing without the geodesics, About 1950, A. A. Isotov, Krasovsky's leading assistant, explained in a lecture that I attended, that systematic errors will then "freely walk" over the entire network.

**655.** [Unnamed persons have] publicly maintained that one single observation, taken wirh due care, was as much to be relied on, as the mean of a great number … (p. 82). Supposing that the several chances for the different errors that any single observation can admit of, are expressed by the terms of the series … T. Simpson (1756, p. 82, both quotations/1757, p. 64, only in the second case).

**Comment.** This was how Simpson explained the essence of his propositions. He thus understood that an observational error was a particular value of a [random variable] and initiated the future theory of errors.

**656.** There is no objection to regard the actual observations as a mere selection, taken at random … from an infinite sequence … and so as generally yielding results that are representative of the whole on the law of distribution of errors. Sampson (1913, p. 173).

**Comment.** Nothing new, but perhaps indicating that, on the other hand, such statements were not generally recognized.

**657.** [The term *Theorie der Fehler*. Its aims: find the] Verhältnis zwischen den Fehlern, ihren Folgen, den Umstanden der Ausmessung und Güte des Instrumente. [Also introduced was the *Theorie der Folgen* aimed at studying the errors of observation and of their functions.] Lambert (1765a, §321).

**Comment.** As far as the theory of errors is considered, Lambert was the main predecessor of Gauss. Unlike Bessel, neither Laplace, nor Gauss ever applied the new term, but in the mid-19th century it suddenly came into vogue.

**658.** Mais pour prendre ce milieu, tel qu'il ne soit point simplement un milieu arithmétique, mais qu'il soit plié par une certaine loi aux règles des combinaisons fortuites et du calcul des probabilités; nous nous servirons ici d'une problème que j'ai indiqué [in 1757] … où je me suis contenté de donner le résultat de la solution. Voici le problème: étant donné un certain nombre de degrés [arcs of the meridian], trouver la correction qu'il faut faire à chacun d'eux, en observant ces trois conditions, la première [application of the proper theory]: la seconde, que la somme des corrections positives soit égale à la somme des négatives: la troisième, que la somme de toutes les corrections, tant positives que négatives, soit la moindre possible, pour le cas où les deux premières conditions soient remplies. … [The second condition is needed] par un même degré de probabilité, pour les déviations du pendule et les erreurs des Observateurs, dans l'augmentation et la diminution des degrés; la troisième est nécessaire pour se rapprocher autant qu'il se pourra des observations. Maire & Boscovich (1770, p. 501).

**Comment.** Pendulum observations served for determining the flattening of the Earth's ellipsoid of revolution. Boscovich was likely one the first to apply the term *calcul des probabilités*. The first was undoubtedly Niklaus Bernoulli; see his Preface to Jakob Bernoulli's *Ars Conjectandi* (1713).

**659.** Si l'on combine ces quatre parties de la méridienne que nous avons déterminées de toutes les manières possibles, il en résulte, avec l'arc total, dix combinaisons qui fournissent dix degrés moyens pour dix latitudes moyennes, et que ces degrés présentent encore tous le même résultat, à une seule irrégularité près, très-légère et qu'il seroit aisé de faire disparaître. Van Swinden et al (an VII, p. 429).

**Comment.** An VII ≈ 1799. The length of a degree of a meridian can be calculated from each arc measurement, and the parameters of the Earth's ellipsoid of revolution, from two such arcs. The combinations were evidently compared one with another for a qualitative check of their closeness. In the same way Boscovich, in 1757 (Cubranic 1961, p. 46), determined the mean latitudinal difference between the endpoints of his meridian arc measurement. He had four of them and calculated the six partial means before determining the general mean. Again, redundant linear equations in two unknowns were solved similarly. In the 19$^{th}$ century, it was proved that, with an appropriate weighting of the pairs of equations, this procedure led to a least-squares solution.

**660.** We shall … find … that the method of least squares, when applied to a system of observations, in which one of the extreme errors is very great, does not generally give so correct a result as the method proposed by Boscovich. … The reason is, that in the former method, this extreme error [like any other one] affects the result in proportion to the second power of the error; but in the other method, it is as the first power … Bowditch, comment on Laplace (1798 – 1825, t. 2, § 40/translation 1832/1966).

**Comment.** More definitely: the reason was the robustness of the median (to which the Boscovich method led). Bowditch translated Laplace's *Traité de mécanique céleste* into English.

**661.** [The principle of least squares] ist … dem La Place'schen vorzuziehen, nach welchem die Summe jener Differenzen = 0, und die Summe derselben Differenzen, aber sämmtlich positiv genommen, ein Minimum sein soll. Man kann zeigen, dass [dies] nach den Gründen der Probabilitätsrechnung nicht zulässig ist, sondern auf Widerspruche führt (p. 329). Ich im Juni 1798 … zuerst La Place's Methode gesehen, und die Unverträglichkeit desselben mit den Grundsätzen der Wahrscheinlichkeitsrechnung in einem kurzen Notizen-Journal … angezeigt habe (pp. 493 – 494). Gauss, letters of 1807 and 1812 to Olbers (1900/*Werke*, *Ergänzungsreihe*, Bd. 4, 1976, pp. 329, 493 – 494).

**Comment.** Gauss wrongly attributed the Boscovich method to Laplace. He (1809b/1887, § 186) then noted that that method led to zero residual free terms of a certain number of equations (this is an important theorem in linear programming!) and considered this consequence unacceptable.

### 4.3. Treatment of Observations. General Considerations

**662.** Aus seiner [Gauss'] mir vorliegenden Protokollen geht vielmehr hervor, dass er auf jeder Station so lange gemessen hat, bis er meinte, dass jeder Winkel sein recht bekommen habe. Er hat dann … die hervorgehenden Richtungswerthe als gleichgewichtig und von einander unabhängig in die Systemausgleichung eingeführt. Schreiber (1879, p. 141).

**Comment.** This is also seen in Gauss' *Werke* (Bd. 9, pp. 278 – 281). Having been a natural scientist as well as a pure mathematician, he did not quite believe in his own formulas for estimating precision, apparently because of the existence of systematic errors.

**663.** Kommt man … früher oder später an die Grenze, um welche [the results] in gewissen geringen Oscillationen herumschwanken, und überzeugt sich, dass jedes weitere Fortsetzen der Repetitionen nur verlorene Arbeit seyn würde. … So aber habe ich es, nach dem Beispiel von Gauss, regelmäßig immer gemacht. Gerling (1839, pp. 166 – 167).

**Comment.** Cournot (1843, §§ 130 and 138) and several later authors (Sheynin 1994, p. 264) expressed themselves in a similar vein. See also No. 664.

**664.** It is a principle in observing generally, that to repeat the same observation over and over, under precisely the same circumstances, is a mere waste of time. Clarke (1880, p. 18).

**665.** To avoid if possible the application of the method of least squares, I prefer to make a few but precise and repeated measurements at several significantly different pressures. … Amassing observations made at various closely spaced pressures not only presents many difficulties for analysis, but also increases the error in drawing inferences. Mendeleev (1872b/1939, p. 144).

**Comment.** Mendeleev described his study of the refinement of the Boyle – Mariotte law. His attitude to the method of least squares was undoubtedly occasioned by the computational difficulties involved. It is possible that under different pressures the systematic errors somewhat changed their values and thus became partly random. Cf. NNo. 663 and 664.

**666.** I believe that isolated observations will add but little to our knowledge, whereas tabulated results from a very large number of observations, systematically made, would probably throw much light on the sequence and period of development of the several [mental] faculties. Darwin, Letter No. 424 of 1881; Darwin (1903, vol. 2, p. 54).

**667.** The mean number from among the data obtained under differing conditions, by different methods and practitioners, is of little value and is always less probable than the result achieved by precise methods and habitual persons. Mendeleev (1872a/1951, p. 101).

**Comment.** That Mendeleev included *differing conditions* seems dubious or at least not applying to metrology.

**668.** When however one of the numbers provides demonstrably better assurance of precision than the others, it alone should be taken into account, ignoring the numbers that either certainly represent worse experimental or observational conditions or give any cause for doubt. … To consider worse numbers taking them with some [even small] weight is tantamount to deliberately corrupting the best number. Mendeleev (1895/1949, p. 159).

**Comment.** See Sheynin (2009a, § 1.1.4) for a similar attitude of ancient astronomers. Then, here is another statement:

*Immerhin bleibt größtmögliche Genauigkeit das erste Requisit numerischer Angaben, und die erste Aufgabe des Statistikers ist daher, sorgfältig sie zu prüfen und nach ihrer Glaubwürdigkeit zu sichten, und nicht Daten von sehr ungleichem Gehalt zu vermengen.* Chr. Bernoulli (1842/1849, p. 10).

**669.** Die so genannten wahrscheinlichen Fehler wünsche ich eigentlich, als von Hypothese abhängig, ganz proscribirt; man mag sie aber berechnen, indem man die mittleren mit 0.6744897 multiplicirt. Gauss, letter of 1825 to Schumacher; Gauss, *Werke*, Bd. 8, p. 143.

**Comment.** It was Bessel who formally introduced the probable error in 1816. Gauss, in spite of his statement above, sometimes applied it in his letters. Owing to its apparent simplicity, the probable error had been employed at least until the mid-20th century, see Sheynin (1979, p. 39). It is strange that Gauss unconditionally adopted one and the same factor for changing from the mean [square] error to it. The highly excessive precision of that factor conformed to the tradition of the day.

**670.** Freilich hat in dem Falle, wo die Anzahl der Beobachtungen vielmale größer ist als die der unbekannten Größen, diese Unrichtigkeit wenig zu bedeuten; allein theils erfordert die Würde der Wissenschaft, dass man vollständig und bestimmt übersehe, wieviel man hierdurch zu fehlen Gefahr läuft, theils sind auch wirklich öfters nach jenem fehlerhaften Verfahren Rechnungsresultate in wichtigen Fällen aufgestellt, wo jene Voraussetzung nicht stattfand. Gauss (1823a/1887, p. 199).

**Comment.** Gauss' opinion about the correct and wrong expressions (the number of the observations, or that number less the number of the unknowns) for the denominator of the formula for the sample variance.

**671.** Pour estimer l'erreur totale [of the Chéops pyramid] à laquelle on serait exposé, il suffisait de multiplier la limite des erreurs partielles par 14, car le nombre de assises est 203 [and $\sqrt{203} = 14$]. Fourier (1829a/1890, p. 569)

**672.** Most statisticians are apt to rely blindly upon the proposition that the random fluctuations of statistical numbers must decrease when the number of trials increases and that we are thus always able to attain any degree of certainty. Chuprov (1918/1926, 1960, p. 230).

**673.** Wäre z. B. eine oder die andere [observation] zur Bestimmung des Werthes von *V* [the observed magnitude] verwendete Beobachtung auch zur Bestimmung des Werthes von *V'* [another such magnitude] benutzt worden, so würden die Fehler [of these observations] nicht mehr von einander unabhängig … Gauss (1823b, § 18/1887).

**Comment.** This is the notion of dependence admitted in the theory of errors. Kapteyn (1912), without knowing about his predecessor, attempted to quantify it. Also see No. 674.

**674.** Die beobachteten Werthe sind unabhängig von einander, wenn jeder aus besonderen Beobachtungen abgeleitet ist. Sind also eine oder mehrere Beobachtungen, die zur Herleitung von A gedient haben, auch zu der von B benutzt werden, so sind A und B abhängig von einander. Schreiber (1882, p. 134).

**675.** La théorie des erreurs était naturellement mon principal but. Poincaré (1921/1983, p. 343).

**676.** [Without the theory of errors his contribution (1925) on stable laws] n'aurait pas de raison d'être. Lévy (1925, p. vii).

**Comment.** For the theory of errors, this contribution proved absolutely useless (Sheynin (1995, § 5).

**677.** [The] fausseté [of the work of Laplace and Gauss] aurait dû apparaître lorsqu'en 1853 Cauchy attira l'attention [to stable laws and in particular to the "Cauchy" distribution]. Lévy (1924, p. 77).

**Comment.** Lévy conditioned the possibility of plausibly estimating the precision of observations by the existence of a stable law of error; Cauchy did not say anything pertinent.

**678.** The analysis of data and the selection therefrom of reliable figures demand so much work and time, that it is more advantageous to make new observations instead … Until the disadvantageous data have been eliminated by a clear critical appraisal, there is no hope of achieving a realistic result. Mendeleev (1887/1934, p. 82).

**Comment.** Astronomers have to analyse ancient observations, for them the situation is much more difficult.

**679.** Zu einer erfolgreichen Anwendung der Wahrscheinlichkeitsrechnung auf Beobachtungen ist allemal umfassende Sachkenntnis von höchster Wichtigkeit. Wo diese fehlt, ist das Ausschließen [of outliers] wegen größerer Differenz immer misslich, wenn nicht die Anzahl der vorhandenen Beobachtungen sehr groß ist … Halte man es wie man will, mache aber zum Gesetz, nichts zu verschweigen, damit andere nach Gefallen auch anders rechnen können. Gauss, letter of 1827 to Olbers; Gauss, *Werke*, Bd. 8, pp. 152 – 153.

**680.** [Newcomb] thr[ew] aside or corrected [an observation] when it appear[ed] probable that it did not or could not correspond to the general mean. Newcomb (1882, p. 372).

[There had been] a decided bias in each generation of astronomers towards depending upon a few recent observations to the exclusion of past ones. Newcomb (1891, pp. 263 – 264).

[Deviating or doubtful observations should be rejected only in the presence of a] well established cause of systematic error. Newcomb (1895, p. 186).

[Rejection was a] matter which has to be left entirely to the judgement of the invesigator. Newcomb (1897, p. 16).
**Comment.** Newcomb acknowledged that he himself was guilty of the "bias". However, combining observations made over a considerable period by different instruments etc. is extremely difficult, but Newcomb had to do so time and time again.
**681.** The question of … rejection … reduces to a question of common sense. … The judgement [of an experienced observer] can undoubtedly be aided by … one or more tests based on the theory of probability [but hardly if a test "requires an inordinate amount of calculation" or is founded on a "complicated hypothesis"]. Rider (1933, pp. 21 – 22).
**682.** One sufficiently erroneous reading can wreck the whole of a statistical analysis however many observations there are. Anscombe (1960, p. 124).
**683.** When all is said and done, the major problem in outlier study remains the one that faced the very earliest workers on the subject – what is an outlier and how should we deal with it? Barnett & Lewis (1978, p. 360).
**684.** Un principe que je n'accepte pas: les observations sont des témoins; si elles sont, avant l'épreuve, jugées dignes de confiance, leur déclaration, quelle qu'elle soit, doit être recueillie et conservée. Bertrand (1888a, p. 305).
**685.** The erroneous numbers, however, are entered in the tables, that it may not be supposed that I have in any one instance tampered with the results. Darwin (1877, p. 150).
**686.** Die Übereinstimmung [of Bradley's observations and the normal law] … ist überall so gut, wie man dies überhaupt … erwarten kann. Aber die schwersten Fehler, die die gewohnten Grenzen weit überschreiten, sind … ein wenig häufiger, als es der Theorie zu entspechen scheint. … Dennoch ist die Differenz zu unbedeutend, als daß sie nicht auch einer noch nicht ausreichend großen Anzahl von Beobachtungen zugeschrieben werden könnte. Bessel (1818); translation Schneider (1988, p. 279).
**Comment.** The three investigated series contained 300, 300 and 470 observations, so that Bessel's explanation was lame. And since larger errors were more frequent, smaller errors should have been less frequent than indicated by theory.
**687.** Owing to unfavourable circumstances under which observations are frequently made, … the arithmetic mean does not necessarily give the most probable result … Any collection of observations of transits of Mercury must be a mixture of observations with different probable errors … (1882, p. 382). In such a mixed system [see above] the most probable result will be, not the arithmetic mean, but a mean obtained by giving less weight to the more discordant observations (1886, p. 346). [Newcomb] set out to modify the usually accepted law in order that it may be applicable to all cases whatever (1886, p. 351). Newcomb (1882, p. 382; 1886, pp. 346, 351).

**Comment.** Observations, especially those made under differing conditions, can undoubtedly obey different laws of distribution. However, Newcomb's proposal to introduce a mixture of normal distributions inevitably involved a subjective choice of the pertinent variances. Then, provided that the law of distribution is more or less even, giving less weight to "discordant" observations amounts to correction for (small) asymmetry. And cases of asymmetrical frequencies ought to be dealt separately. Finally, Newcomb mistakenly thought that his innovation was appropriate for astronomical observations in general. On this subject and on further academic generalizations see Sheynin (1995, pp. 179 – 183).

**688.** I submit … in the absence of evidence to the contrary, that *non-exponential* laws of error of the kind…

$$y = \frac{c}{(a^2 + x^2)^2}$$

do occur in *rerum natura*, that the "ancient solitary reign" of the exponential law should come to an end… Edgeworth (1883, pp. 305 – 306/1996, pp. 194 – 195).

**689.** If [an] archer … [aiming his arrows at a set mark] makes innumerable shots, all with the utmost possible care, … the hits must be assumed to be thicker and more numerous on any given band the nearer it is to the mark. If all the places … whatever their distance from the mark, were equally liable to be hit, the most skilful shot would have no advantage over a blind man. That, however, is the tacit assertion of those who use the common rule in estimating the value of various discrepant observations, when they treat them all indiscriminately [and choose the arithmetic mean]. (pp. 157 – 158). The mean between the two extreme observations … as a rule for several observations I have found to be less often wrong than I thought before I investigated the matter (p. 161). Daniel Bernoulli (1778/1970, pp. 157 – 158, 161).

**Comment.** Daniel alleges that, when choosing the arithmetic mean, a deviating observation excessively influences the result. However, as Karl Pearson (1978, p. 268) commented, although only qualitatively, since *small errors were more probable than large ones, they would be more frequent and have their due weight in the arithmetic mean.* Without mentioning Lambert, Daniel recommended the principle of maximum likelihood; cf. No. 717. Concerning the mean between extreme observations, nothing is known about his pertinent study. On the former use of that mean in meteorology see Sheynin (1984b, p. 74).

### 4.3. The True Value of a Measured Constant and the Method of Least Squares

**690.** En prenant … un milieu entre un grand nombre d'observations, on court peu de risque de se tromper. Condamine (1751, p. 223).
**691.** [Given an even curve of frequency (in modern terms),] Das Mittel aus mehrern Versuchen dem wahren desto näher kommen müsse, je mehr der Versuch ist wiederholt worden. Lambert (1765a, § 3).
**692.** Si l'on multiplie indéfiniment les observations …, leur résultat moyen converge vers un terme fixe, de manière qu'en prenant de part et d'autre de ce terme un intervalle aussi petit que l'on voudra, la probabilité que le résultat moyen tombera dans cet intervalle finira par ne différer de la certitude que d'une qunatité moindre que toute grandeur assignable. Laplace (1795/1912, p. 161).
**693.** Supposons donc que l'on ait ajouté ensemble un grand nombre de valeurs observées, et que l'on ait divisé la somme par le nombre *m*, ce qui donne la quantité A pour la valeur moyenne; nous avons déjà remarque que l'on trouverait presque exactement cette même valeur A, en employant un très grand nombre d'autres observations. En général, si l'on excepte des cas particuliers et abstraits que nous n'avons point à considérer, la valeur moyenne ainsi déduite d'un nombre immense d'observations ne change point; elle a une grandeur déterminée H, et l'on peut dire que le résultat moyen d'un nombre infini d'observations est une quantité fixe, où il n'entre plus rien de contingent, et qui a un rapport certain avec la nature des faits observés. C'est cette quantité fixe H que nous avons en vue comme le véritable objet de la recherche. Fourier (1826/1890, pp. 533 – 534).
**Comment.** This definition recalling the Mises understanding of probability had been either forgotten or even unnoticed and many authors independently from each other and from Fourier introduced it anew, see NNo. 694 and 695. On occasion, Gauss applied the term *real value*. A few years later Fourier (1829a/1890, p. 583) additionally noted that a new long series of observations will provide the same mean which was thus a "quantité entièrement fixe" and that the difference between it and the constant measured "peut devenir moindre que toute quantité donnée".
**694.** The mass of a mass standard is … specified … to be the mass of the metallic substance of the standard plus the mass of the average volume of air adsorbed upon its surface under standard conditions. Eisenhart (1963/1969, p. 31).
**Comment.** Eisenhart made use of the same definition of the real value. As an unavoidable consequence, it occurred that the residual systematic errors were included in that value.
**695.** Der "wahre" Wert der Beobachtung (d. i. derjenige, der sich als Durchschnitt bei einer ins Unendliche fortgesetzten Beobachtungsreihe ergeben müsste) (p. 40). Der "wahre" Mittelwert … ist nicht anderes als die Größe, die nach der Definition des Wahrscheinlichkeitsbegriffes als arithmetisches Mittel einer ins Unendliche fortgesetzten Ziehungsserie sich erheben müsste (p. 46). Mises (1919/1964, pp. 40, 46).
**696.** Das es … einen "wahren Wert" … gibt, folgt aus der grundlegenden Annahme, wonach die Meßergebnisse ein Kollektiv bilden. *Ein anderer Begriff des "wahren Wertes" als dieser hat in einer auf Wahrscheinlichkeitsrechnung gestützten Theorie der Beobachtungsfehler keinen Raum.* Mises (1931, p. 370).

**697.** Mais il faut surtout faire en sorte que les erreurs extrêmes, sans avoir égard à leurs signes, soient renfermées dans les limites les plus étroites qu'il est possible. De tous les principes qu'on peut proposer pour cet objet, je pense qu'il n'en est pas de plus général, de plus exact, ni d'une application plus facile, que celui dont nous avons fait usage dans les recherches précédentes, et qui consiste à rendre *minimum* la somme des quarrés des erreurs. Par ce moyen il s'établit entre les erreurs une sorte d'équilibre qui, empêchant les extrêmes de prévaloir, est très-propre à faire connaître l'état du système le plus proche de la vérité. Legendre (1805, pp. 72 – 73).

**Comment.** Unlike the minimax principle, Legendre's innovation did not at all lead to a least interval of possible errors (more precisely, of the residual free terms of the initial equations). Legendre's formulation thus involved two mistakes.

**698.** Übrigens ist unser Princip, dessen wir uns schon seit dem Jahre 1795 bedient haben, kürzlich auch von Legendre … aufgestellt worden … Gauss (1809b/1887, § 186).

**Comment.** Later Gauss (1823b/1887, § 17) again claimed the principle of least squares although not as resolutely as before.

**699.** Je ne vous dissimulerai-donc pas, Monsieur, que j'ai éprouvé quelque regret de voir qu'en citant mon mémoire …, vous dites *principum nostrum* … Il n'est aucune découverte qu'on ne puisse s'attribuer en disant qu'on avoit trouvé la même chose quelques années auparavant; mais si on n'en fournit pas la preuve en citant le lieu où on l'a publiée, cette assertion devient sans objet et n'est plus qu'une chose désobligeante pour le véritable auteur de la découverte. En Mathématiques il arrive très souvent qu'on trouve les mêmes choses qui ont été trouvées par d'autres et qui sont bien connues; c'est ce qui m'est arrivé nombre de fois, mais je n'en ai point fait mention et je n'ai jamais appelé *principum nostrum* un pr[incipe] qu'un autre avait publié avant moi. Vous êtes assez riche de [votre] fonds, Monsieur, pour n'avoir rien à envier à personne; et [je suis] bien persuadé au reste que j'ai à me plaindre de l'expression seulement et nullement de l'intention … Legendre, letter of 1809 to Gauss; Gauss, *Werke*, Bd. 10/1, p. 380.

**700.** J'ai fait usage de la méthode des moindre[s] carrés depuis l'an 1795 et je trouve dans mes papiers, que le mois de Juin 1798 est l'époque où je l'ai rapprochée aux principes du calcul des probabilités. … Cependant mes applications *fréquentes* de cette méthode ne datent que les l'année 1802, depuis ce temps j'en fait usage pour ainsi dire tous les jours dans mes calculs astronomique[s] sur les nouvelles planètes. … Je ne me suis hâté d'en publier un morceau détaché, ainsi Mr. Legendre m'est prévenu. Au reste j'avais déjà communiqué cette même méthode, beaucoup avant la publication de l'ouvrage de M. Legendre, à plusieurs personnes, entre autres à Mr. Olbers en 1803 … Ainsi, pouvais je dans ma *théorie* [1809b] parler de la méthode des moindre[s] quarrés, dont j'avais fait depuis 7 ans mille et mille applications, … je dis, pouvais je parler de ce principe, que j'avais annoncé à plusieurs de mes amis déjà en 1803 comme devant faire partie de l'ouvrage que je préparois, – comme d'une méthode *empruntée* de Mr. Legendre? Je n'avait pas l'idée, que Mr. Legendre pouvait attacher tant de prix à une idée aussi simple, qu'on doit plutôt s'étonner qu'on ne l'a pas eue depuis 100 ans, pour se fâcher que je raconte, que je m'en suis servi avant lui? … Mais j'ai cru que tous ceux qui me connaissent le croiraient même sur ma parole, ainsi que je l'aurait cru de tout mon cœur si Mr. Legendre avait avancé, qu'il avait possédé la méthode déjà avant 1795. J'ai dans mes papiers beaucoup de choses, donc peut être je pourrai perdre la priorité de la publication: mais soit, j'aime mieux faire mourir les choses. Gauss, letter of 1812 to Laplace; Ibidem, pp. 373 – 374.

**701.** Mein Wahlspruch ist aut Caesar, aut nihil. Gauss, letter to Olbers 30.7.1806 (1900/Gauss, *Werke, Ergänzungsreihe*, Bd. 4/1, 1976, p. 307).

**702.** When I had received Dr. Olbers's observations till April 17, for curiosity's sake I attempted to apply to them the same method, which I had made use of in my calculations about Ceres Ferdinanden, and which without any hypothetical suppositions yields the true conic section as exactly as the nature of the problem & the precision of the observations will permit. Gauss, letter to Maskelyne dated 19 May 1802, in English; MS, Roy. Greenwich Obs., Code 4/122:2; also Sheynin (1999a, p. 255).

**Comment.** Strictly speaking, this passage does not yet prove that Gauss had then applied the method of least squares.

**703.** Mit Unwillen und Betrübnis habe ich … gelesen, dass dem alten Legendre, der eine Zierde seines Landes und seines Zeitalters ist, die Pension gestrichen hat. Gauss, letter of 1824 to Olbers (1900; *Werke, Ergänzungsreihe*, Bd. 4, Tl. 1, 1976, p. 413).

**704.** [Gauss had] noch immer an der Vervollkommnung der Methode selbst gearbeitet, besonders in dem vorigen Winter, und ihre jetzige Gestalt sieht ihrer ersten fast gar nicht mehr. Gauss, letter of 1806 to Olbers; Gauss, *Werke*, Bd. 6, pp. 275 – 277

**705.** M. Legendre eut l'idée simple de considérer la somme des carrés des erreurs des observations, et de la rendre un minimum, ce qui fournit directement autant d'équations finales, qu'il y a d'éléments à corriger. Ce savant géomètre est le premier qui ait publié cette méthode; mais on doit à M. Gauss la justice d'observer qu'il avait eu, plusieurs années avant cette publication, la même idée dont il faisait un usage habituel, et qu'il avait communiquée à plusieurs astronomes. Laplace (1812/1886, p. 353).

**706.** Was einem normalen Autor verboten ist, einem Gauss wohl gestattet werden muss, zumindest müssen wir seine Gründe respektieren. Biermann (1966, p. 18).

**707.** Gauss cared a great deal for priority. … But to him this meant being first to discover, not first to publish; and he was satisfied to establish his dates by private records, correspondence, cryptic remarks in publications. … Whether he intended it so or not, … he maintained the advantage of secrecy without losing his priority in the eyes of later generations. May (1972, p. 309).

**708.** Gauss bereits im Junius 1803 die Güte hatte, mir diese Methode als längst von ihm gebraucht, mitzuteilen und mich über die Anwendung derselben zu belehren. Olbers (1816, p. 192n)
**Comment.** Already in 1812; (1900/Gauss, *Werke, Ergänzungsreihe*, Bd. 4, 1976, Tl. 1, p. 495) assured Gauss that he will confirm *gern und willig* that he came to know the principle of least squares from him (from Gauss) before 1805. In 1812 – 1815 Olbers had only published a few papers on the observation of comets (*Catalogue of Scientific Literature*, Roy. Soc. London).
**709.** [Bessel also came to know the principle of least squares before 1805] durch eine mündliche Mittheilung von Gauss. Bessel (1832/1848, p. 27).
**710.** Le célèbre Docteur Gauss était déjà depuis 1795 en possession de cette méthode, et il s'en est servi avec avantage dans la détermination des élémens des orbites elliptiques des quatre nouvelles [minor] planètes, comme on peut voir dans son bel ouvrage [of 1809]. Von Zach (1813, p. 98n).
**Comment.** Still, the *Theoria motus* does not directly prove the point mentioned. Von Zach's testimonial is nevertheless important because he had been unjustly accused (not publicly) of unwillingness to confirm Gauss' priority. It occurred in addition that until 1809 he only knew that Gauss had applied some new method of adjusting observations.
**711.** Wie ein Axiom pflegt man nämlich die Hypothese zu behandeln, wenn irgend eine Größe durch mehrere unmittelbare, unter gleichen Umständen und mit gleicher Sorgfalt angestellte Beobachtungen bestimmt worden ist, dass alsdann das arithmetische Mittel zwischen allen beobachteten Werthen, wenn auch nicht mit absoluter Strenge, so doch wenigstens sehr nahe den wahrscheinlichsten Werthe hebe, so dass es immer das sicherste ist, an diesem festzuhalten. Gauss (1809b/1887, § 177).
**712.** Liegen für eine Größe mehren Werte aus Beobachtungen vor, so ist ihr wahrscheinlichsten Wert das arithmetische Mittel alles beobachteten Werte. D. Hilbert, from an unpublished lecture of 1905; Corry (1997, p. 161).
**713.** Tel est le postulatum [of the arithmetic mean] de Gauss. Bertrand (1888a, p. 176).
**Comment.** He also called it an axiom. His term took root.
**714.** La fonction φ (*z*) est donc, en réalité, assimilée à une exponentielle de la forme … Bertrand (1888a, p. 267).
**Comment.** Bertrand noted that, for small values of its argument, an even quadratic function, the sum of the first terms of a power series, was approximately identical with the exponential function of the negative square; in other words, that the second Gauss justification of least squares had not after all done away with the normal law.
**715.** En proposant en 1809 une hypothèse sur la théorie des erreurs, [Gauss] ne prétendait nullement établir la vérité, mais la chercher. Bertrand (1888a, p. XXXIV).
**716.** Wenn die unabhängigen Beobachtungsfehler … einzeln das [Normal] Gesetz befolgen, so unterliegt eine homogene lineare

Funktion … derselben einen Gesetz der gleichen Form … Czuber (1903, p. 23).
**Comment.** He called that proposition the "Hauptsatz" of the theory of errors. It was proved by Bessel (1838/1876, § 7).
**717.** Das ich übrigens die in der *Theoria motus* … angewandte Metaphysik für die Methode der kleinsten Quadrate späterhin habe fallen lassen, ist vorzugsweise auch aus einem Grunde geschehen, den ich selbst öffentlich nicht erwähnt habe. Ich muss es nämlich in alle Wege für weniger wichtig halten, denjenigen Werth einer unbekannten Größe auszumitteln, dessen Wahrscheinlichkeit die größte ist, die ja doch immer nur unendlich klein bleibt, als vielmehr denjenigen, an welchen sich haltend man das am wenigsten nachtheilige Spiel hat; oder wenn *fa* die Wahrscheinlichkeit des Werthes *a* für die Unbekannte *x* bezeichnet, so ist weniger daran gelegen, dass *fa* ein Maximum werde, als daran, das $\int fx \, F(x - a) \, dx$ ausgedehnt durch alle möglichen Werthe des *x* ein Minimum werde, in dem für *F* eine Funktion gewählt wird, die immer positiv und für größere Argumente auf eine schickliche Art immer größer wird. Gauss, letter of 1839 to Bessel; Gauss, *Werke*, Bd. 8, pp. 146 – 147.
**Comment.** Metaphysics had apparently been understood as a speculative principle, perhaps not representing reality well enough; for example, leading to the universality of one single distribution.
**718.** [Markov] often heard that my [his] presentation [of the method of least squares] is not sufficiently clear. Markov, letter of 1910 to Chuprov (Ondar 1977/1981, p. 21).
**Comment.** See also Markov's correspondence with B. M. Koialovich, his former student and future professor of physics (Sheynin 2004, pp. 220 – 226).
**719.** When justifying the method of least squares, Markov effectively introduces … notions equivalent to the concepts of unbiased and effective estimators … Linnik et al (1951, p. 637/2004, p. 238).
**Comment.** Instead of Markov, they should have mentioned Gauss.
**720.** The importance of the work of Markov concerning the best linear estimates consists, I think, chiefly in a clear statement of the problem. Neyman (1934, p. 593).
**Comment.** Even this is wrong. Moreover, on p. 595 he mentioned a *Markov theorem*. Later Neyman (1938, p. 52) made known the forthcoming appearance of David & Neyman (1938) where these authors proved "an extension of the Markov theorem". Again, Neyman (1938, pp. 131 and 133) stated that the Russian scientist had calculated the variance of a linear function of [dependent] arguments and indirectly attributed to him the principles of unbiasedness and least variance also known to Gauss.

Then, finally, Neyman (1938/1952, p. 228) admitted the "confusion" to which he "inadvertently contributed by attributing to Markov the main theorem of the method of least squares". The non-existing Gauss – Markov theorem is still alive (Dodge 1990/2003).
**721.** I consider it [Gauss' second justification of the method of least squares] rational because it does not obscure the conjectural essence of the method … We do not ascribe the ability of providing the most probable, or the most plausible results to the method … Markov (1899b/1951, p. 246; translation 2004, pp. 138 – 139).

**Comment.** Markov had somehow forgotten that Gauss did state that, according to his new substantiation, the method furnished most plausible results. This is all the more strange since Markov, in the same memoir, resolutely defended that second justification.

**722.** For stark clarity of exposition the presentation [of the principle of least squares by Legendre] is unsurpassed; it must be counted as one of the clearest and most elegant introductions of a new statistical method in the history of statistics. Stigler (1986, p. 13).

**Comment.** This is wrong, see No. 697 and my comment there.

**723.** Legendre immediately realized the method's potential … it was not merely applications to the orbits of comets he had in mind (p. 57). There is no indication that [Gauss] saw its [the method's] great general potential before he learned of Legendre's work (p. 146). Ibidem, pp. 57, 146.

**Comment.** The first statement seems likely, but the second one is wrong (and disgusting), see Sheynin (1999a; 1999b). Instead of searching for proof or refutation of a fact, it is much easier to assert the second possibility out of the blue! Incidentally, it is likely that Gauss could have also applied least squares for trial computations and/or in a simplified way (Gauss 1809b/1887, § 185) and this is of course impossible to refute.

**724.** Except for one circumstance, Gauss's argument [in 1809] might have passed relatively unnoticed, to join an accumulating pile of essentially ad hoc constructions, a bit neater than some but less compelling than most. That one circumstance was the reaction it elicited from Laplace. Ibidem, p. 143.

**Comment.** Computations made by Gauss enabled other astronomers to find the first minor planet which had disappeared after its first sightings, and this alone has immortalized the Gauss method. In addition, the Gauss *argument* had since been repeated in hundreds of treatises. So does it occur that Laplace had wrongly appraised the importance of Gauss' deliberations, and that only Stigler explained to all of us the real situation? As to the *accumulating pile*, *ad hoc constructions* and *less compelling*, let all this nonsense remain on Stigler's conscience.

**725.** Gauss solicited reluctant testimony from friends that he had told them of the method before 1805. Ibidem, p. 145.

**Comment.** This is libel, pure and simple, see NNo. 708 and 709.

**726.** Although Gauss may well have been telling the truth about his prior use of the method, he was unsuccessful in whatever attempts he made to communicate it before 1805. Ibidem, p. 146.

**Comment.** The beginning of the phrase is appropriate with respect to a suspected rapist, but not to Gauss. The ending of the phrase is meaningless. Should have Gauss published his finding in a newspaper, or announced it through a public crier? Taken as a whole, Stigler's utterances about Gauss (and Euler, see No. 651) are abominable. Here, however, is the opinion of Hald (1998, p. XVI), one of the most authoritative statisticians and historians of statistics: Stigler's book *is epochal*, – in spite of what I stated above and in spite of its not being a history of statistics (as claimed by its title) but a description of a few of its chapters. Epochal, as we ought to believe, along with the works of Newton and Einstein. I also note Tee's reasonable criticisms (1991) of Stigler.

**727.** I distinguish the viewpoints of Gauss and Laplace [on the method of least squares] by the moment with regard to the experiment. The first one is posterior and the second one is prior. It is more opportune to judge à posteriori because more data are available, but this approach is delaying, it lags behind, drags after the event. P. A. Nekrasov, letter of 1913 to Markov; Sheynin (2009a, § 14.5).

**Comment.** This statement is wrong (and meaningless as well). Laplace issued from the central limit theorem which he proved non-rigorously whereas Gauss studied the finite case. Then, a sound appraisal of the precision attained (*to judge*) is only possible à posteriori. See the next items for a discussion of Nekrasov's personality.

**728.** Nekrasov reasons perhaps deeply but not clearly and he expresses his thoughts still more obscurely. I am only surprised that he is so self-confident. In his situation, with the administrative burden weighing heavily upon him, it is even impossible, as I imagine, to have enough time for calmly considering deep scientific problems, so that it would have been better not to study them at all. K. A. Andreev, letter to Liapunov dated 1901; Gordevsky (1955, pp. 40 – 41).

**Comment.** From 1893, Nekrasov was Rector of Moscow University, and from 1898 onward, a highly ranking official at the Ministry of People's Education. He graduated from an Orthodox seminary and his religious upbringing and administrative duties resulted in a change of his personality. Around 1900, his writings (only on probability and statistics!) became unimaginably verbose, sometimes obscure and/or meaningless, and inseparably linked with ethical, political and religious considerations. In a letter of 1916 to P. A. Florensky, a mathematically oriented religious philosopher, Nekrasov stated that the comparison of "Christian science" and philosophers, such as Florensky and he himself, with Karl Marx, Markov and a hardly known author with a Jewish name "clearly shows the crossroads to which the German-Jewish culure and literature are pushing us" (Sheynin 2003c, p. 343).

**729.** The revolution had come, and Nekrasov decided to direct his entire talent towards serving the proletariat. He definitely attempted to grasp the Marxist system. Uritsky (1924).

**730.** [Nekrasov's economic concepts were] equally hostile both to the capitalist and the socialist principles. A. V. Andreev (1999, pp. 105 – 106).

**731.** I, together with some other freshmen, attended Nekrasov's course "Theory of probability". It was the last year that he delivered lectures. He simply read aloud his book (I do not know which one), and, according to some indications, he did not go into the essence of the material read. His course proved to be absolutely useless. Beskin (1993, pp. 168 – 169).

**Comment.** On his p. 164 Beskin stated that he had entered Moscow University in 1921.

**732.** [During the first half of the 1920s] Nekrasov still attended the meetings of the Moscow Mathematical Society and sometimes even presented papers. A queer shadow of the past, he seemed decrepit – physically and mentally – and it was difficult to understand him. … This pitiful old man was like a shabby owl. Liusternik (1967, p. 222).

**Comment.** See Sheynin (2003c, p. 339) for other telltale episodes proving that for the powers that be Nekrasov remained a *persona non grata*. Perhaps only Sluginov (1927), a petty mathematician, favourably commemorated Nekrasov and called him "Professor at First Moscow State University" (from 1917, as it seems, to 1930 that University consisted of two bodies). This is difficult to understand the more so since apparently, as follows from Sluginov's acount, Nekrasov had been working during his last year(s) in the provincial city of Perm.

**733.** Man sieht, dass diese Formel in dem Falle $n = m$ unter die ganz unbestimmte Form 0/0 tritt, und in der Tat ist in diesem Falle gar kein Schluss auf die Präzision gestattet. Dedekind (1860/1930, p. 97).

**Comment.** Dedekind believed that his remark proved the inadequacy of the Gauss formula for the sample variance. However, he was mistaken: under the stated condition (no redundant observations available), regardless of any formulas, estimation of precision is only possible by means of indirect considerations. Cf. No. 670.

**734.** Gauss' Ansatz ist gekünstelt und wenig überzeugend (p. 177). Die Ausgleichung erreichte ihren theoretischen Abschluss erst mit I. J. Bienaymé (p. 178). Freudenthal & Steiner (1966, pp. 177, 178).

**Comment.** I can only note (Heyde & Seneta 1977, pp. 66 – 71) that the principle of least variance makes more sense when applied to the estimation of all the unknowns at once.

**735.** Gauss' second exposition seems … to be no more satisfactory than his first. In each case he starts from a postulate, plausible but not universally valid, which leads inexorably to the foregone conclusion. Harter (1977, p. 28).

**Comment.** Harter certainly had in mind the variance which Gauss had chosen as the main estimator of precision. Kolmogorov (1946, p. 64) decided that the formula for the sample variance ought to be considered as its definition, and Tsinger (1862, § 3) remarked that it already *concealed* the principle of least squares.

**736.** [Principle of least squares] vielleicht nicht mit Unrecht [should be left unjustified] (p. 44). … würde ich für Zweckmässig halten wenn man zur Begründung der Methode der kleinsten Quadrate die Wahrscheinlichkeitstheorie überhaupt nicht mehr in Anspruch nehmen würde (p. 67). Henke (1868/1894, pp. 44, 67).

**Comment.** The two statements are hardly contradictory because no justification is possible without stochastic considerations.

**737.** L'application du Calcul des probabilités à l'étude des erreurs d'observation repose sur une fiction. Bertrand (1888a, p. 222).

**Comment.** Because, as he mistakenly explained, it was founded on the existence of a certain law of distribution.

**738.** [The determination of the orbits of celestial bodies. Lagrange and Laplace] restricted their attention to the purely mathematical aspect whereas Gauss thoroughly worked out his solution from the point of view of computations and took into account all the conditions of the work of astronomers and [even] their habits. Subbotin (1956, p. 297).

**Comment.** The necessity of adjusting a finite number of observations which Gauss understood perfectly well, was quite in line with the author's conclusion.

**739.** [Gauss had often been led to his discoveries] durch peinlich genaues Rechnen … Wir finden [in his works] ganzen Tafeln deren Herstellung allein die Lebensarbeit manches Rechners vom gewöhnlichen Schlage ausfällen würde. Maennchen (1930, p. 3).

**Comment.** Maenchen did not study Gauss' geodetic computations, likely because at that time mathematicians had not yet been interested in examining solutions of systems of linear algebraic equations.

**740.** Some authors have recently begun to assume the expression [the correct formula for the sample variance] … Chebyshev (1879 – 1880, p. 249/2004, p. 228).

**741.** [In 1919 the author only had] un vague souvenir du fait que les erreurs accidentelles obéissent à la loi de Gauss. Lévy (1970, p. 71).

**Comment.** At best, this is only approximately true. Lévy's statement goes to show that at the time mathematicians had not mastered the theory of errors.

**742.** La loi de Gauss est bien la seule pour laquelle cette méthode [of least squares] s'applique. Lévy (1925, p. 79).

**Comment.** Not *la seule*, but best possible. This time, Lévy justified his statement by a practically meaningless conclusion concerning stable laws. Poisson (1824) was the first to derive the *Cauchy distribution*. However, Lévy (1970, p. 78) "propose d'appeler [it] loi de Cauchy". This also was a mistake.

**743.** In the cases to which it is appropriate this method [of least squares] is a special application of the Method of Maximum Likelihood, from which it may be derived. Fisher (1925/1990, p. 260).

**Comment.** This is only true in the sense of the Gauss first justification (1809b) of the method.

**744.** The probable conclusion … perfectly agrees here with the arithmetic mean, and this indicates that the [observational] errors follow a definite law assumed by the Gauss theory of probability [of errors], i. e., that the observations do not contain large random deviations, but are subject to unavoidable observational errors. Mendeleev (1875/1950, p. 209).

**Comment.** *Probable conclusion* is an imprecise and therefore unfortunate substitute for median. In 1809, Gauss arrived at a *definite* law by assuming that its mode (not median!) coincided with the arithmetic mean. That definite (the normal) law nevertheless allows large errors to happen (with low probabilities).

**745.** [The existence of the second formulation of the method of least squares] seems to be virtually unknown to all [of its] American users … except students of advanced mathematical statistics. Eisenhart (1964, p. 24).

**746.** Gauss and Laplace are representatives of two absolutely different opinions on the meaning of the method of least squares. In Laplace's work we find a rigorous [?] and impartial study of this problem. His analysis shows that the results of the method … only enjoy a more or less substantial probability when the number of observations is large whereas Gauss attempted to attach absolute meaning to this method, using extraneous considerations. If we turn our attention to the fact that all the essence of the Theory of chances is contained in the law of large numbers, and that all the properties of random phenomena only take real importance when the number of trials is large, it would not be difficult to perceive the correctness of the Laplacean inference. However, when the number of observations is limited, we cannot at all reckon upon the mutual cancellation of errors … and … any combination of observations can …lead as much to the increase of errors as to their decrease. Tsinger (1862, p. 1).

**Comment.** The author was ignorant of the second Gaussian justification of the method of least squares; of Gauss' qualification remark (1823b/1887, § 6) about the arbitrariness of his method; and of Gauss' correct decision to restrict his attention to the case of a small number of observations. Finally, both the history of the sciences of observation and of mathematical statistics proved that Tsinger's last lines contradicted reality and theory, respectively.

**747.** Tout le monde y croit cependant, me disait un jour M. Lippmann, car les expérimentateurs s'imaginent que c'est un théorème de mathématiques, et les mathématiciens que c'est un fait expérimental. Poincaré (1896, § 108).

**Comment.** Poincaré discussed the belief in the universality of the normal distribution. Already Bessel (No. 686) noticed that observations deviated from normality but paid no attention to that fact. At the mid-19$^{th}$ century such deviations became sufficiently known, see Sheynin (1995, pp. 174 – 177) and NNo. 486 and 487.

**748.** Our results … justify the method of least squares without postulating a Gaussian error law provided that it is understood that the method is not concerned with "most probable values" (pp. 275 – 276). The proof [of the Gaussian law] from the "principle of the arithmetic mean", which is still sometimes quoted, is altogether fallacious (p. 276 note). No practical physicist has ever set himself to make the most accurate determination of a physical constant; he usually wishes to make a more accurate determination than his predecessors; sometimes his intention is merely to obtain a determination which is *accurate enough* … The calculator may well have similar aims. Eddington (1933, pp. 275 – 276 and 271 – 272).

**Comment.** Did Eddington know about the second Gauss' justification of the method? And why is that proof fallacious? My last quotation partly explains his attitude.

**749.** The theory of errors is the last surviving stronghold of those who would reject plain facts and common sense in favour of remote unverifiable guesses, having no merit other than mathematical tractability. N. R. Campbell, discussion of Eddington (1933), see p. 283. Cf. No. 745.

**Comment.** Gauss rather than Campbell was endowed with common sense (and did not at all reject plain facts). Again, ignorance of the second justification of the method.

**750.** Bertrand die Unsicherheit der [Gauss'sche] Formel nach dem absoluten Betrage ihres mittleren Fehlerquadrates beurtheilt, statt, wie er sein muss, den relativen Betrag … zugrunde zu legen. Czuber (1891, p. 460).

**Comment.** Czuber wrote these lines after discussing the issue with Helmert. Bertrand (1888b) calculated the variance of the variance and of his own estimator for a chosen example and normal law and stated that his estimator was better. He did not notice that, unlike the variance, his estimator was biased. And, when calculating, he forgot about the formula provided by Gauss for the normal distribution. In practice, the precision of observations is measured by the mean square error (rather than by the variance), which is biased. In the same paper, Czuber proved that for the normal distribution the relative variance of the variance was minimal.

**751.** I like very much Bredikhin's rule according to which 'in order to admit the reality of a computed quantity, it should at least twice numerically exceed its probable error'. I do not know, however, who established this rule or whether all experienced calculators recognized it. Markov, 1903; Sheynin (1990a, pp. 453 – 454).

**Comment.** Both Newcomb and Mendeleev had also been applying the same rule (Sheynin 1989, p. 351).

### 4.4. The Determinate Branch of the Theory of Errors

**752.** Wir verdanken Ihnen den größten Theil der heutigen Verfeinerung der Astronomie, nicht nur wegen Ihrer kleinsten Quadrate, sondern auch wegen der Erweckung der Sinns für Feinheit, der seit Bradley's Zeit verschwunden zu sein schien und erst seit 18 Jahren wieder erschien. Wir sind erst jetzt auf den Punkt gekommen, kleinen Fehlern oder Abweichungen außer der Grenzen der Wahrscheinlichkeit mit derselben Aufmerksamkeit nachzuspuren als früher großen … F. W. Bessel, letter of 1818 to Gauss; 1880/Gauss, *Werke*, *Ergänzungsreihe*, Bd. 1, 1975, p. 275.

**Comment.** Bessel's own merit is here also considerable, see however Sheynin (2000). Contrary to the implication made (Biermann 1966), the relations between Gauss and Bessel had not at all been serene. Thus, Bruhns (1869, p. 108) testified that a noisy scientific quarrel between them had taken place in 1825. Then, in 1844 Bessel stressed Legendre's priority in the discovery of the principle of least squares (Sheynin 2001b, p. 168).

**753.** [Investigation of the precision of projected geodetic networks. The aim:] Einen nothwendigen Genauigkeitsgrad mit möglichst wenig Zeit und Geld zu erreichen (p. 1). [Or, to achieve more precise results] bei gleicher Mühe (p. 60). Helmert (1868, pp. 1, 60 of Dissertation).

**754.** Ich glaube nicht zu irren, wenn ich annehme, dass, nach einer Reihe von Jahren, das erste Capitel aller Lehrbücher der auf Erfahrung beruhenden Wissenschaften, der Anwendung der Wahrscheinlichkeits-Rechnung auf die Beobachtungskunst gewidmet sein wird. Bessel (1848a, p. 398).

**755.** Sie haben nie die Verpflichtung anerkannt, durch zeitige Mitteilung eines dem ganzen angemessenen Theils Ihrer Forschungen die gegenwärtige Kenntniss der Gegenstände derselben zu befördern; Sie leben für die Nachwelt. Dieses ist aber ganz gegen meine Ansicht … Bessel, letter of 28 May 1837 to Gauss. Gauss (1880/*Werke, Ergänzungsreihe*, Bd. 1, 1975, pp. 516 – 520).

## 5. Randomness

### 5.1. Previous Attitudes

**756.** A certain man drew his bow at a venture and struck the King of Israel. 1 Kings 22: 34; 2 Chronicles 18:33. Similar examples: 2 Samuel 1: 6 and 20:1.

**757.** Digging a hole for a tree, someone finds a treasure [not a rusty nail!]. Aristotle, *Metaphysics* 1025a.

**Comment.** This is his example of a chance event. Another of his examples (*Physics* 196b30): an unforeseen coming across an acquaintance.

**758.** By chance Aristotle means something that takes place occasionally; has the character of an end; is such that it might have been the object of a natural or of a rational appetite but came into being by accident. Junkersfeld (1945, p. 22).

**759.** [Mistakes] in the operations of nature [give rise to] monstrosities. [The first departure of nature from the] type is that the offspring should become female instead of male; … as it is possible for the male sometimes not to prevail over the female … either through youth or age or some other such cause … Aristotle, *Physics* 199b1.

**Comment.** Mistakes therefore occur very often. How are we to combine this with Aristotle's understanding of chance as a rare event? Laplace (1814/1886, p. L) thought that the study of the sex ratio of newly-born animals deserved attention, and Darwin (Sheynin 1980, p. 346) studied this problem statistically.

**760.** Si l'on veut avoir plus de femelles il faut employer des mâles jeunes et des femelles dans l'âge de la force, et nourrir celles-ci plus abondamment que ceux-la. Il faut faire l'inverse si l'on veut produire plus de mâles. Cuvier (1831, p. clxxxvii).

**Comment.** See also Babbage (1829, p. 91).

**761.** Practically all other horoscopic instruments … are frequently capable of error, the solar instruments by the occasional shifting of their positions or of their gnomons, and the water clocks by stoppages and irregularities in the flow of water from different causes and by mere chance. Ptolemy (1956, III, 2, p. 231).

**Comment.** That book was devoted to astrology. Ptolemy believed that the influence of heaven was a tendency rather than a fatal drive, cf. No. 643.

**762.** Such a distribution of animals, plants and climates as exists, is not the result of design – just as the difference of race, or of language, is not, either – but rather of accident and chance. Strabo (1969, 2.3.7)**.**

**Comment.** A strange statement of an ancient geographer and historian.

**763.** Nichts ist der Berechnung und der Regelmäßigkeit so sehr entgegengesetzt wie der Zufall. Ciceron (1991, Buch 2, § 17, p. 149).

**764.** The crow had no idea that its perch would cause the palm-branch to break, and the palm-branch had no idea that it would be broken by the crow's perch; but it all happened by pure Chance. Belvalkar & Ranade (1927, p. 458).

**765.** When a traveller meets with a shower, the journey had a cause, and the rain had a cause …; but because the journey caused not the rain, nor the rain the journey, we say they were contingent one to another. Hobbes (1646/1840, p. 259).

**766.** Les événements amenées par la combinaison ou la rencontre de phénomènes qui appartiennent à des séries indépendantes, dans l'ordre de la causalité, sont ce qu'on nomme des événements *fortuites* ou des résultats du *hasard*. Cournot (1843, § 40).

**Comment.** Such an understanding of randomness had been typical in particular for ancient scholars, see NNo. 756 and 757.

**767.** Some measure [time] …by continuous motion … and, as a rule, this has been done by the use of water. But it is subject to variation in many repects. For instance, the purity and density depend on its sources. … Also, it is subject to accidental variations … Al-Biruni (1967, pp. 155 – 156).

**768.** In those who are healthy … the body does not alter even from extreme causes; but in old men even the smallest causes produce the greatest change. Galen (1951, p. 202).

**769.** Le nez de Cléopatre s'il eût été plus court toute la face de la terre aurait changé. Pascal (1669/2000, p. 675).

**770.** Chance is nothing real in itself … is merely a negation of a cause … 'Tis essential to [chance] to leave the imagination perfectly indifferent, either to consider the existence or non-existence of that object, which is regarded as contingent. Hume (1739/1969, Book 1, pt. 3, § 11, p. 125).

**771.** Creatures that arise spontaneously are called automatic …because they have their origin from accident, the spontaneous act of nature. Harvey (1651/1952, p. 338).

**Comment.** Spontaneous generation had been generally accepted and randomness was thus thought to be an extremely important agent.

**772.** I am not willing to ascribe this wonderful coincidence in time and place to blind chance, especially since the appearance of a new star, on its own, even without regard to time and place, is no ordinary thing, like throwing dice, but a great wonder, whose like has not been heard of or written about before our own time. Kepler (1604/1977, p. 337).

**Comment.** Kepler considered unusual both the time and the place as well.

**773.** Was aber ist Zufall? Wahrlich, er ist ein höchst abscheulicher Götze und nicht anderes als eine Beschimpfung des höchsten und allmächtigen Gottes und der höchst vollkommemem Welt, die er schuf. Kepler (1606/2006, p. 163).

**Comment.** Nevertheless Kepler was compelled to to state that the eccentricities of the planetary orbits were random. Also see No. 774.

**774.** If the celestial movements were the work of mind, as the ancients believed, then the conclusion that the routes of the planets are perfectly circular would be plausible. … But the celestial movements are … the work of … nature … and this is not proved by anything more validly than by the observation of the astronomers, who … find that the elliptical figure of revolution is left in the real and very true movement of the planet; and the ellipse bears witness to the natural bodily power and to the emanation and magnitude of its form. … Those [natural and animal] faculties … could not attain perfection completely. Kepler (1618 – 1821, 1620/1952, p. 932).

**775.** Woher sind ihre Umläufe nicht vollkommen zirkelrund? … Ist es nicht klar einzusehen, daß diejenige Ursache, welche die Laufbahnen der Himmelskörper gestellt hat, … es nicht völlig hat ausrichten können …Kant (1755/1910, 1. Hauptstücke, p. 269).

**776.** Si le système solaire s'était formé avec une parfaite régularité, les orbites des corps qui le composent seraien des cercles, dont les plans, ainsi que ceux des divers équateurs et des anneaux, coïncideraient avec le plan de l'équateur solaire. Mais on conçoit que les variétés sans nombre qui ont du exister dans la temperature et la densité des diverses parties de ces grandes masses ont produit les excentricités de leurs orbites et les déviations de leurs mouvements du plan de cet équateur. Laplace (1796/1884, p. 504).

**Comment.** Neither Kant, nor even Laplace (in 1796) knew about Newton's explanation of the elliptical motion of the planets!

**777.** [Laplace's *Exposition* (1796)] is an ingenious epitome of the principal discoveries. … If he writes the history of great astronomical discoveries, he becomes a model of elegance and precision. No leading fact ever escapes him. … Whatever he omits, does not deserve to be cited. Fourier (1829b, p. 379).

**778.** Zufällige Dinge [are those] deren vollkommener Beweis jeden endlichen Verstand überschreitet. Leibniz (1686/1960, p. 288).

**779.** Ce serait un défaut de raison de s'imaginer que notre esprit étant fini, il put comprendre jusqu'où peut aller puissance de Dieu, qui est infini. Arnauld & Nicole (1662/1992, p. 317).

**780.** Blind chance could never make all the planets move one and the same way in orbs concentrick, some inconsiderable irregularities excepted, which may have risen from the mutual action of comets and planets upon one another, and which will apt to increase, till this system wants a [divine] reformation. Such a wonderful uniformity in the planetary system must be allowed the effect of choice. And so must the uniformity in the bodies of animals. Newton (1704/1931, Query 31).

**Comment.** Many scholars both before and after Newton thought that a random origin of the system of the world, or the universe was impossible, see No. 615 and Sheynin (1974, pp. 132 – 135).

**781.** *Contingent* (both *free*, if it depends on the free will of a reasonable creature, and *fortuitous* and *casual*, if it depends on fortune or chance) is that which can either exist or not exist at present, in the past or future, – clearly because of remote rather than immediate forces. Indeed, neither does contingency always exclude necessity up to secondary causes. Jakob Bernoulli (1713/2005, p. 15).

**782.** La Nature contient le fonds de toutes ces variétés: mais le hasard ou l'art [in case of domestic animals] les mettent en œuvre. Maupertuis (1745/1756, p. 110).

**783.** [Maupertuis was] the first to apply the laws of probability to the study of heredity … virtually every idea of the Mendelian mechanism of heredity and the classical Darwinian reasoning from natural selection and geographical isolation is … combined [in his works] together with De Vries' theory of mutations as the origin of species. B. Glass (1959, p. 60).

**784.** [In 1751 Maupertuis] provided the first careful and explicit analysis of the transmission of a dominant hereditary trait in man. … He calculated the mathematical probability that the trait would occur coincidentally in the three successive generations of [a certain] family [in Berlin] had it not been inherited. On the basis of this study Maupertuis founded a theory of the formation of the fetus and the nature of heredity that was at least a century ahead of its time. B. Glass (1974, p. 188).

**785.** La nature ne forme rien, elle détruit toujours (§ 920). Toute substance composée tend naturellement à ce détruire … (§ 482). Tous les efforts de la nature … sont perpétuellement dirigés vers ce seul but; savoir, d'opérer la destruction des composés quels qu'ils soient, et de rendre aux élémens qui les constituent, la liberté et leurs qualités naturelles, dont ils sont dépourvus dans leur état de combinaison (§ 813). Lamarck (1794, §§ 920, 982, 813).

**Comment.** This is a heuristic statement concerning the not yet existing thermodynamics (increasing enthropy).

**786.** The word *chance* only signifies ignorance of causes (p. 607). The aims of meteorology would have been absolutely useless, unreliable and groundless had there existed any part of nature which … would not obey invariable rules, if that, which is called chance, might be a reality (p. 632). Lamarck (1810 – 1814/1959, pp. 607 and 632).

**Comment.** Soon Lamarck repeated this statement (1815, p. 329).

**787.** La nature a deux moyens puissans et généraux, qu'elle emploie continuellement … Ces moyens sont 1° L'attraction universelle …; 2° L'action répulsive des fluides subtils, mis en expansion; action qui, sans être jamais nulle, varie sans cesse dans chaque lieu, dans chaque temps, et qui modifie diversement l'état de rapprochement des molécules des corps. De l'équilibre entre ces deux forces opposées … naissent ... les causes des tous les faits que nous observons, et particulièrement de ceux qui concernent l'existence des corps vivans. Lamarck (1815, p. 169).

**788.** Une … cause accidentelle et par conséquent variable a traversé ça et la l'exécution de ce plan [des opérations de la nature, – i. e., the layout of the tree of animal life]. Ibidem, p. 133.

**789.** Cette cause nécessaire et inconnue que l'on nomme hasard, mot dont l'étude de ce calcul peut seul bien faire connaître le veritable sens. Condorcet (1795/1988, p. 251).

**790.** Dans une série d'événements indéfiniment prolongée, l'action des causes régulières et constantes doit l'emporter à la longue sur celle des causes irrégulières. C'est ce qui rend les gains des loteries aussi certaines que les produits de l'Agriculture, les chances qu'elles se réservent leur assurant un bénéfice dans l'ensemble d'un grand nombre de mises. Ainsi, des chances favorables et nombreuses étant constamment attachées à l'observation des principes éternels de raison, de justice et d'humanité qui fondent et maintiennent les sociétés, il y a un grand avantage à ce conformer à ces principes et de graves inconvénients à s'en écarter. Laplace (1814/1886, p. XLVIII).

**791.** The regularity that often escapes us when considering a small number of conclusions reveals itself if a large number of them is taken into account. W. J. 'sGravesand, in Latin, quoted by Chuprov (1905, § 2) without mentioning the exact source. N. S. Chetverikov translated that paper into Russian.

**792.** Der Zufall im Einzelnen nichts desto weniger einer Regel im Ganzen unterworfen ist … Kant (1781/1911, p. 508).

**793.** Das Zufällige ist ein Wirkliches, das zugleich nur als möglich bestimmt … was möglich ist, ist selbst ein Zufälliges. Hegel (1812/1978, p. 383).

**794.** Supposing that Gain and Loss were so fluctuating, as always to be distributed equally, … would it be reasonable … to attribute the Events of Play to Chance alone? I think … it would be quite otherwise, for then there would be more reason to suspect that some unaccountable Fatality did rule in it … De Moivre (1718/1756, p. iv).

**795.** Nous voyons sur une table des caractères d'imprimerie, disposés dans cet ordre, *Constantinople*, et nous jugeons que cet arrangement n'est pas l'effet du hasard, non parce qu'il est moins possible que les autres, puisque si ce mot n'était employé dans aucune langue, nous ne lui soupçonnerions point de cause particulière; mais, ce mot étant en usage parmi nous, il est incomparablement plus probable qu'une personne aura disposé ainsi les caractères précédents qu'il ne l'est que cet arrangement est dû au hasard. Laplace (1776/1891, p. 152; 1814/1886, p. XV).

**Comment.** This example is due to D'Alembert (1767a, pp. 254 – 255) who offered another and not really convincing explanation.

**796.** Le mot *hasard* n'exprime donc que notre ignorance sur les causes des phénomènes que nous voyons arriver & se succéder sans aucun ordre apparent. Laplace (1786/1893, p. 296).

**797.** Le mot hasard sert officieusement à voiler notre ignorance. Quetelet (1846, p. 14).

**798.** Si nous voyons arriver tout autre événement qui présente quelque chose de symétrique, tel que …, nous sommes portés croire que ces événements réguliers ne sont pas l'effet du hasard … Poisson (1837, p. 39).

**Comment.** For more detail see Ibidem, p. 114ff.

**799.** Une intelligence qui, pour un instant donné, connaîtrait toutes les forces dont la nature est animée et la situation respective des êtres qui la composent, si d'ailleurs elle était assez vaste pour soumettre ces données à l'Analyse, embrasserait dans la même formule les mouvements des plus grands corps de l'univers et ceux du plus léger atome: rien ne serait incertain pour elle, et l'avenir, comme le passé, serait présent à ces yeux. Laplace (1814/1886, pp. VI – VII).

**Comment.** In actual fact, Laplace had thus admitted randomness but explained it by ignorance. I believe that his statement therefore does not prove his determinism at all. Then, during the latest decades, physicists and mechanicians have begun to admit for randomness an incomparably greater part to say nothing about the existence of unstable motion. Finally, both Maupertuis (1756, p. 300) and Boscovich (1758, § 385) kept to the "Laplacian determinism". Both mentioned calculations of the past and future ("to infinity on either side", as Boscovich maintained), but, owing to obvious obstacles, disclaimed any such possibilities.

**800.** [The action of selection] absolutely depends on what we in our ignorance call spontaneous or accidental variability. Darwin (1868/1885, vol. 2, p. 236).

**Comment.** A similar understanding of randomness is seen in his other contributions (1859/1958, p. 128) as well.

**801.** Throw up a handful of feathers, and all fall to the ground according to definite laws; but how simple is the problem where each shall fall compared with problems in the evolution of species. Darwin (1859/1958, p. 77).

**Comment.** Darwin had not mentioned randomness.

**802.** I am inclined to look at everything as resulting from designed laws, with the details, whether good or bad, left to the working out of what we may call chance. Darwin, letter of 1860; Darwin (1887, vol. 2, p. 312).

**803.** [The use of chance only in relation to purpose in the origination of species.] This is the only way I have used the word chance. … On the other hand, if we consider the whole universe, the mind refuses to look at it as the outcome of chance, – that is, without design or purpose. Darwin, Letter No. 307 of 1881; Darwin (1903, vol. 1, p. 395).

**804.** Hat ein Philosoph, von dem richtigen Gedanken ausgehend, dass alle für die Menschen erreichbare Wahrheit doch nur durch Worte ausgedrückt werden könne, alle Wörter seiner Sprache … geschrieben, und eine Maschine erfunden, welche diese auf allen Seiten beschriebenen Würfel nicht nur wendete, sondern auch in einander verschob. Nach jeder Handhabung der Maschine … wenn drei oder vier mit einander einen Sinn gaben, wurde diese Wortfolge notirt, um auf diesem Wege zu aller möglichen Weisheit zu gelangen. Baer (1873, p. 6).

**Comment.** Baer, Danilevsky (1885, pt. 1, p. 194), and still earlier John Herschel (1861/1866, p. 63) mentioned the ill-starred philosopher from *Gulliver Travels* borrowed by Swift from Lully (13[th] – 14[th] centuries). The three authors reduced randomness to its uniform case for denying the Darwinian hypothesis. See No. 94 for a criticism of the uniform randomness *in general*.

**805.** [The 19[th] century will be called] das Jahrhundert der mechanischen Naturauffasung, das Jahrhundert Darwins … Boltzmann (1872/1909, p. 316).

**Comment.** Boltzmann twice repeated this idea 15 – 20 years later, see Sheynin (1985, p. 375). Hardly anyone else had called the Darwin hypothesis mechanical. Boltzmann's mistaken opinion was apparently occasioned by his ambiguous attitude towards the part of randomness.

**806.** [About Boltzmann:] Randomness … struggles with mechanics. Mechanical philosophy is still able … to overcome randomness and wins a Pyrrhic victory over it, but recedes undergoing a complete ideological retreat. Rubanovsky (1934, p. 6).

**807.** A special development of [Jakob Bernoulli's law of large numbers] is noticed in Darwin's [*Origin of species*] and in Buckle's *History of Civilization in England* [1857 – 1861]. From these writings we see that life is a struggle in which even insignificant but constant forces produce great effects. Vashchenko-Zacharchenko (1864, p. 69).

**808.** Quelquefois les chances ne sont soumises à aucune loi appréciable, et la courbe de possibilité peut affecter les formes les plus capricieuses. Quetelet (1846, p. 182).

**Comment.** This is a very rare and indirect mention of chaos whose existence De Moivre (1733/1756, p. 251) had simply denied.

**809.** Le hasard est sans vertu: impuissant dans les grandes affaires, il ne trouble que les petites. Mais, pour conduire les faits de nature à un fin assurée et précise, il est, au milieu des agitations et des variétés infinies, le meilleur et le plus simple des mécanismes. Bertrand (1888a, p. L).

**810.** [A variable probability in problems concerning drawings without replacement is] en quelque sorte un régulateur de la proportion normale prévue par le théorème de Bernoulli. Ibidem, p. 166.

**811.** Nous sommes devenus des déterministes absolus. … Les mot hasard n'aurait pas de sens, ou plutôt il n'y aurait pas de hasard. C'est a cause de notre faiblesse et de notre ignorance qu'il y en aurait un pour nous. Poincaré (1896/1912, p. 2).

### 5.2. New Ideas. Modernity

**812.** When an infinitely small variation in the present state may bring about a finite difference in the state of the system in a finite time, the condition of the system is said to be unstable. It is manifest that the existence of unstable conditions renders impossible the prediction of future events … Maxwell (read 1873/1884, pp. 362 – 363).

**Comment.** On p. 364 Maxwell indicated as an example the unstable direction of a ray within a biaxial crystal. He had not mentioned randomness which in similar examples only explicitly appeared in Poincaré's illustrations (No. 813).

**813.** Une cause très petite … détermine un effet considérable … et alors nous disons que cet effet est dû au hasard. Poincaré (1896/1912, p. 4).

**814.** Dans chaque domaine, les lois précises ne décidaient pas de tout, elles traçaient seulement les limites entre lesquelles il était permis au hasard de se mouvoir. Ibidem, p. 1.

**815.** Zufällig nennen wir eine Erscheinung, deren Abhängigkeit von den sie hervorrufenden Ursachen eine so complicirte ist, dass sie mit

Hülfe von nur analytischen Functionen gar nicht ausgedrückt werden kann. N. N. Pirogov (1891, p. 518).
**Comment.** Cf. No. 778.
**816.** Il était bon de sygnaler l'extreme mobilité que le hasard pourrait produire dans la composition de la chambre élective … Poisson (1837, p. 259).
**Comment.** Given proportion 9:10 of the strength of two parties, he proved that there was a very low probability of electing a deputy belonging to the less numerous of them.
**817.** Es seien $N$ Kugeln (z. B. 100) fortlaufend nummeriert über zwei Urnen verteilt. Die Urne A enthalte $P_0$ (z. B. 90), die Urne B also $Q_0 = (N - P_0)$ Kugeln. … Überdies befinden sich in einem sack $N$ Lotteriezettel mit den Nummern $1 - N$. Nach je einer Zeiteneinheit wird ein Zettel gezogen und zurückgelegt. – Jedesmal, wenn eine Nummer gezogen wird, hüpft die Kugel mit dieser Nummer aus der Urne, in der sie gerade liegt, in die andere Urne und bleibt dort so lange liegen, bis gelegentlich wieder ihre Nummer gezogen wird. P., T. Ehrenfest (1907/1959).
**Comment.** This celebrated model is considered as the origin of the history of stochastic processes. See No. 818 and Sheynin (2009a, §§ 6.1.1 and 7.1-3) for the previous results of Daniel Bernoulli and Laplace.
**818.** [Laplace considered the same pattern as in No. 817 with many urns, containing balls of two colours, and even with new urns being added, again with any initial distribution of balls of the two colours, and concluded, apparently too optimistically:] On peut étendre ces resultats à toutes les combinaisons de la nature, dans lesquelles les forces constantes dont leurs éléments sont animés établissent des modes réguliers d'action, propres á faire éclore du sein même du chaos des systèmes régis par des lois admirables. Laplace (1814/1886, p. LV).
**819.** There exist very simple systems whose behaviour is just as unpredictable [as that of complicated systems]. Ekeland (1991/1993, p. 91).
**820.** However paradoxical it is, all mathematical definitions of randomness known at present are formulated in terms of the theory of algorithms. Kolmogorov & Uspensky (1987, p. 425).
**821.** The von Mises' Theory constitutes a great and real advance both in the theory of probability proper and expecially in the philosophy of the subject. It is known but little, there are many regrettable misunderstandings about it and it is certainly in the interest of American science to have here one more competent person who could successfully teach it. Neyman (letter of 1940) as quoted by Siegmund-Schulze (1998, p. 259).
**Comment.** Neyman recommended Hilda Geiringer to the Rockefeller Foundation.
**822.** For a long time I had the following views. (1) The frequency concept based on the notion of *limiting frequency* as the number of trials increases to infinity, does not contribute anything to substantiate the applicability of the results of probability theory to real practical problems where we have always to deal with a finite number of trials. (2) The frequency concept applied to a large but finite number of trials does not admit a rigorous formal exposition within the framework of pure mathematics. …
I still maintain the first of the two theses mentioned above. As regards the second, however, I have come to realize that the concept of random distribution of a property in a large finite population can have a strict formal mathematical exposition. Kolmogorov (1963, p. 369).

**823.** I remember conversations [in the family] which included the name of Gumbel … as well as Kolmogorov who was mentioned as having differing viewpoints but being nevertheless a friend. … I think it quite likely that there was some correspondence [with Kolmogorov, likely after Mises had quit Germany]. My mother kept for many years a … picture postcard [to her] of a red troika from him, probably with New Year's greetings. M. Ticza, daughter of Hilda Geiringer, the wife of R. von Mises; two private communications, both dated 2004.

**824.** When von Mises was persuaded to make his theory mathematical, a collective became a mathematical sequence of ones (heads) and zeros (tails) … This took the intuitive content out of the approach and left only an awkward technique to get by detours to the usual mathematical probability discipline. Doob (1976, p. 203).

**825.** The concept of chance has no part in the Mises theory. Is this obscure and contradictory notion necessary for science, and is it not possible to replace it by other more precise and definite concepts? Romanovsky (1924, No. 4 – 6, p. 33n).

**826.** Until now, it proved impossible to embody Mises' intention in a definition of randomness that was satisfactory from any point of view. Uspensky et al (1990, § 1.3.4).

**827.** Supposons actuellement qu'au lieu de deux événements possibles … il y en ait un nombre donné λ, dont un seul devra arriver à chaque épreuve. Ce cas est celui où l'on considère une chose A d'une nature quelconque, susceptible d'un nombre λ de valeurs, connues ou inconnues, … et parmi lesquelles une seule devra avoir lieu à chaque éprouve … Poisson (1837, pp. 140 – 141).

**Comment.** Also see his p. 254. Note that Poisson called the random variable by a purely provisional term.

**828.** The 1970 draft lottery has not helped to mitigate the doubts of many regarding the quality and fairness of random drawings. Fienberg (1971, p. 261).

**Comment.** The doubts apparently concerned the possibility of ensuring an equal probability of the drawing of all the numbers. This is a very simple example of difficulties encountered in generating random and pseudorandom numbers.

**829.** The peculiar motions of the stars are directed at random, that is, they show no preference for any particular direction. I shall further on refer to this hypothesis as the fundamental hypothesis. Kapteyn (1906a, p. 400).

**830.** Der Weise … sucht das vertraute Gesetz in des Zufalls grausenden Wundern, sucht den ruhenden Pol in der Erscheinungen Flucht. Schiller, *Spaziergang*, 1795.

**831.** The words *chance* and *genius* do not designate anything really existing and cannot therefore be defined. These words only designate a certain degree of understanding the [pertinent] phenomena. Tolstoy (1865 – 1869/1999, vol. 2, Epilogue, p. 528).
**832.** Random errors have all the properties of random quantities. Vasiliev (1885, p. 133).
**Comment.** *Random quantity* still remains the Russian term and Vasiliev was the first to introduce it. The equivalent and better English expression *random variable* appeared in Whitworth (1901) and perhaps in earlier editions of that source. Random errors of observation are practically always corrupted by systematic influences.

# Bibliography

   In many cases I mention more than one edition of the same source. Thus,
**Darwin, C.** (1871), *The Descent of Man*. London, Murray, 1901. [London, Penguin, 2004.]
   This means that that source was first published in 1871; that, in spite of having mentioned Darwin (1871) in the text itself, I refer to the edition of 1901; and that there exists an edition of 2004. When possible, I only mention Russian sources in translation; however, most of my own translations were published in a small number of copies, and in such cases I mention the original Russian sources (and provide the appropriate page numbers) as well.

*Abbreviation*
AHES         = *Arch. Hist. Ex. Sci.*
Hist. Scient. = *Historia Scientiarum* (Tokyo)
IMI          = *Istoriko-Matematich. Issledovania*
JNÖS         = *Jahrbücher f. Nationalökonomie u. Statistik*
OC           = *Oeuvr. Compl.*
R            = in Russian
VS           = *Vestnik Statistiki*

# Index of Names

*The Index does not cover either the Introduction, or the References. Letter C means "Comment". When a name appears both in a certain item and in the accompanying Comment, the C is not entered. I have not included my own name.*